\documentclass[hidelinks,10pt]{amsart}

\usepackage{amsmath,amssymb,stmaryrd,dsfont,bm}
\usepackage[pdftex]{graphicx}
\usepackage{subcaption}
\usepackage[all]{xy}
\usepackage[hyperfootnotes=false]{hyperref}

\input xy
\xyoption{all}

\usepackage{pdfsync}
\usepackage{color}

\usepackage[margin=2cm]{geometry}
\usepackage[font=small]{caption}

\newcommand{\beq}{\begin{equation}}
\newcommand{\ee}{\end{equation}}
\newcommand{\bal}{\begin{aligned}}
\newcommand{\eal}{\end{aligned}}
\newcommand{\bac}{\begin{array}{c}}
\newcommand{\ea}{\end{array}}

\newcommand{\Ec}{\mathcal{E}}

\newcommand{\Ss}{\mathrm{S}}
\newcommand{\Tt}{\mathrm{T}}

\newcommand{\vertiii}[1]{{\left\vert\kern-0.25ex\left\vert\kern-0.25ex\left\vert #1 
    \right\vert\kern-0.25ex\right\vert\kern-0.25ex\right\vert}}

\begin{document}
\newtheorem{theorem}{Theorem}[section]
\newtheorem{lemma}[theorem]{Lemma}
\newtheorem*{conjecture*}{Conjecture}
\newtheorem{observation}[theorem]{Observation}
\newtheorem{problem}{Problem}
\newtheorem{definition}[theorem]{Definition}
\newtheorem{example}[theorem]{Example}
\newtheorem{claim}[theorem]{Claim}
\newtheorem{hypothesis}{Hypothesis}
\newtheorem{remark}[theorem]{Remark}

\title[Regularized semiclassical limits]{ Regularized semiclassical limits: \\linear flows with infinite Lyapunov exponents 
}


\author[Athanassoulis]{Agissilaos  Athanassoulis}
\address[A.  Athanassoulis]{Dept. of Mathematics, University of Leicester, UK}
\email{agis.athanassoulis@le.ac.uk}

\author[Katsaounis]{Theodoros Katsaounis}
\address[Th.  Katsaounis]{Computer  Electrical and Mathematical Sciences \& Engineering, King Abdullah Univ. of Science and Technology (KAUST), Thuwal,  KINGDOM of SAUDI ARABIA,  \& IACM--FORTH, Heraklion, GREECE}
\email{theodoros.katsaounis@kaust.edu.sa}

\author[Kyza]{Irene Kyza}
\address[I.  Kyza]{Division of Mathematics,
University of Dundee,
Dundee, DD1 4HN, 
Scotland, UK \& Institute of Applied and Computational Mathematics--FORTH,
Nikolaou Plastira 100, Vassilika Vouton, Heraklion-Crete, GREECE. }
\email{ikyza@maths.dundee.ac.uk}

\thanks{
I.K.\ was partially supported by the ESF-EU and National Resources of the Greek State within the framework of the Action ``Supporting Postdoctoral Researchers'' of the EdLL Operational Programme. Th.K. and A.A  were partially supported by European Union FP7 program Capacities (Regpot 2009-1), through ACMAC (http://www.acmac.uoc.gr).}

\begin{abstract} 
Semiclassical asymptotics for linear Schr\"odinger equations with non-smooth potentials give rise to ill-posed formal semiclassical limits. These problems
have attracted a lot of attention in the last few years, as a proxy for the treatment of eigenvalue crossings, i.e. general systems. It has recently been shown that the semiclassical limit for conical singularities is in fact well-posed, as long as the Wigner measure stays away from singular saddle points.
In this work we develop a family of refined semiclassical estimates, and use them to derive regularized transport equations for saddle points with infinite Lyapunov exponents, extending the aforementioned recent results. In the process we answer a related question posed by P. L. Lions and T. Paul in 1993. 
If we consider more singular potentials, our rigorous estimates break down. 
To  investigate whether conical saddle points, such as $-|x|$, admit a regularized transport asymptotic approximation, we employ a numerical solver based on posteriori error control. Thus rigorous upper bounds for the asymptotic error in concrete problems are generated. In particular, specific phenomena which render invalid any regularized transport for $-|x|$ are identified and quantified. In that sense our rigorous results are sharp. 
Finally, we use our findings to formulate a precise conjecture for the condition under which conical saddle points admit a regularized transport solution for the WM. 
\end{abstract}

\maketitle

\noindent  \textbf{Keywords:} \textit{semiclassical limit for rough potential, multivalued flow, selection principle, a posteriori error control}



\section{Introduction}
The study of the linear Schr\"odinger equation in the {semiclassical} regime
\beq \label{eq:schroeq}
\begin{aligned}
i\hbar u^\hbar_t +\frac{\hbar^2}{2} \Delta u^\hbar -Vu^\hbar & =  0, \qquad u^\hbar(t=0) =u^\hbar_0, \\ 
\|u_0^\hbar\|_{L^2(\mathbb{R}^d)} &=  1,  \qquad \hbar \ll 1
\end{aligned}
\ee
arises naturally in many problems of engineering and mathematical physics, see e.g. \cite{GMMP,RPK} and the references therein. 
A standard physical interpretation is that of the dynamics for a quantum particle, the behavior of which is expected to resemble classical mechanics as $\hbar\rightarrow 0$, hence the term ``semiclassical''.   

Semiclassical problems  appear in many applications; these include  long distance  propagation  for parabolic and hyperbolic wave equations \cite{GMMP,RPK}, or long-distance  paraxial propagation \cite{BEHJ,LLM,MS,LQO}. Certain mean-field limits of statistical mechanics  give rise to semiclassical limits \cite{golse} as well.  Molecular dynamics and modeling of chemical reactions is   another source of semiclassical limit problems \cite{T,saddlep}; it is there that singular saddle points arise as particularly important problems. Singular saddle points are also related to eigenvalue crossings, which develop even in  smooth systems. This has been an important motivation for their focused study \cite{FG,FGL}. 

Since  the direct solution of \eqref{eq:schroeq} becomes intractable for $\hbar \ll 1$, several asymptotic techniques have been developed for its approximation.
Semiclassical asymptotics can be said to be completely understood for problems with $V\in C^{1,1}(\mathbb{R}^d)$; difficulties arise for less regular potentials, which will henceforth be called ``non-smooth''. In particular, $W^{1,\infty}(\mathbb{R}^d)$ non-smooth potentials may arise as effective potentials in smooth systems \cite{FG,FGL,LT}, or from first principles modeling \cite{AFG,HTL,Miel,dHO}. 
For the relation of this type of non-smooth potentials and semiclassical limits of non-linear Schr\"odinger equations see, e.g., \cite{pinaud,golse}. Recent breakthroughs in the semiclassical limits of problem \eqref{eq:schroeq} over non-smooth potentials are \cite{AFFPG,FGL}. 

As has been highlighted in \cite{FGL}, the regularity of the potential is not as important as the overall behavior of the underlying classical flow, $\phi_t:(x,k)\mapsto(X(t),K(t)),$ where
\beq
\begin{aligned}
\label{eq:fflochar1}
\dot X(t)=2\pi K(t),& \qquad \dot K(t)=-\frac{1}{2\pi} \partial_x V (X(t)), \\
X(0)=x,&\qquad \hspace{0.0cm} K(0)=k.
\end{aligned}
\ee
It is well known that the flow $\phi_t$ defined by \eqref{eq:fflochar1} is well-defined for all $(x,k)\in \mathbb{R}^{2d}$ as long as $V\in C^{1,1}(\mathbb{R}^d)$, and that the problem \eqref{eq:fflochar1} has weak solutions (possibly many)  for all $(x,k)\in \mathbb{R}^{2d}$ as long as $V\in C^{1}(\mathbb{R}^d)$. This basic observation puts nicely in context why any $V \notin C^{1,1}(\mathbb{R}^d)$ is called non-smooth. 

Once we enter the regime of non-smooth potentials, the regularity of the flow is much more relevant than the smoothness of $V$.
For example, in \cite{FGL} there are results applicable, essentially, to any $V \in W^{1,\infty}(\mathbb{R}^d)$, \textit{under the assumption that the wavefunction does not fully interact with a singular saddle point}. One way to explain why singular saddle points are so different, is that the flow around one has an infinite Lyapunov exponent. (See also section \ref{sec:lou} for a more detailed discussion of singular saddle points.) Thus, out of the possible isolated non-smooth points, local maxima are the most challenging mathematically, since they give rise to singular saddle points. They are also particularly interesting physically,  since they model chemical reactions \cite{saddlep}. 

Finally, we note that despite the great progress of the last 30 years, the semiclassical limits for the non-smooth saddle point problem described in Remarque IV.3 of \cite{LP},  have not been computed before this work. We state a generalized version of this problem as Problem \ref{prblm:1}, in Section \ref{sbs:pssttl}, and solve it in Theorem \ref{thrm:selprin}.

\subsection{Statement of the Main Results: Semiclassical Estimates for $\mathbf{V\in C^{1,a}(\mathbb{R}^d)}$}

One of the main results of this paper is the computation of semiclassical limits for problems of the form \eqref{eq:schroeq} with potentials  $V \in C^{1,a}(\mathbb{R}^d)$, $a\in (0,1)$. 
For these results, a family of test functions (related to the widely used Banach algebra $\mathcal{A}$) will be used; namely

\begin{definition}[$\mathcal{B}_M$] We will denote by $\mathcal{B}_M$ the Banach space of functions $f:\mathbb{R}^{2d}\to \mathbb{R}$ defined in terms of the norm
\[
\vertiii{f}_M:=   \sum\limits_{m=0}^\infty M^{-m} \, \| \,|K|^m \widehat{f}(X,K)\|_{L^1(\mathbb{R}^{2d})} .
\]
The dual space, denoted by $\mathcal{B}_{-M}$ has norm
\[
\vertiii{f}_{-M}:=  \sup\limits_{\vertiii{\phi}_M=1} \langle f, \phi \rangle .
\]
\end{definition}

\begin{remark}  \upshape Some key observations for $\mathcal{B}_M$:
\begin{itemize}
	\item[(i)]  (Non-triviality of $\mathcal{B}_M$.) One easily sees the following: take a Schwarz test-function $\phi(x,k) \in \mathcal{S}(\mathbb{R}^{2d})$, which in addition has the property that its Fourier transform, $\widehat{\phi}(X,K)$, is supported on $|K|<L$. Then  $\phi \in \mathcal{B}_M$ for all $M>L$.
\item[(ii)]
(Relation to $\mathcal{A}$.) It is straightforward to observe that 
$$\|\phi(x,k)\|_{\mathcal{A}} = \| [\mathcal{F}_{k\to K}  \phi](x,K) \|_{L^1_KL^\infty_x} \leqslant \| \widehat{\phi}(X,K) \|_{L^1_KL^1_X} \leqslant \| \phi\|_{\mathcal{B}_M}, \qquad \forall M> 0.$$
\end{itemize}
\end{remark}

All the notations used here are  precisely defined in section \ref{sbs:not}. The definition of the algebra $\mathcal{A}$ is given in appendix \ref{sec:appA}. The main result is stated in the following theorem and its proof can be found in section \ref{sbs:prmth}.

\begin{theorem} \label{thrm:main}Consider the semiclassical IVP \eqref{eq:schroeq}, with
a potential $V \in C^{1,a}(\mathbb{R}^d)$ for some $a>0$, satisfying also $\int\limits_S \widehat{V}(S) |S| dS < \infty$. Denote by $W^\hbar(t)$ the Wigner transform (introduced in detail in section \ref{sub:WT}) of the wavefunction $u^\hbar$,
\[
W^\hbar(t) = \int\limits_y e^{-2\pi i k y} u^\hbar(x+\frac{\hbar y}{2},t) \overline{u}^\hbar(x-\frac{\hbar y}{2},t)dy,
\]
and by $\rho^\hbar(t)$ the solution of
\beq
\label{eqsert45}
\partial_t \rho^\hbar(t) + 2\pi k\cdot \partial_x \rho^\hbar(t) - \frac{1}{2\pi} \partial_x V \cdot \partial_k \rho^\hbar(t)=0, \,\,\,\,\,\,\,\,\,\,\,\,\,\,
\rho^\hbar(t=0)=\rho^\hbar_0.
\ee
Then there exists a constant $C>0$ so that
\beq \label{eq:dpexcluestform}
\vertiii{ W^\hbar(t) -\rho^\hbar(t) }_{-M} \leqslant C e^{C t} \Big( \vertiii{W^\hbar(0)-\rho^\hbar(0)}_{-M} + \hbar^a \Big)
\ee
In particular, if $\rho^\hbar_0 = W^\hbar[u_0^\hbar]$,
$
\vertiii{ W^\hbar(t) -\rho^\hbar(t) }_{-M} \leqslant C e^{C t}  \hbar^a .
$
\end{theorem}

\begin{remark} \upshape One should note that:
\begin{itemize}
	\item[(i)] In particular the assumptions allow for
$V$ with localized singularities of the form $C|x|^{1+a}$; 
see lemma \ref{lm:technV} for more details.
\item[(ii)] Using Theorem \ref{thrm:main}, a selection principle for the multivalued flow \eqref{eq:fflochar1} can be constructed, see Theorem \ref{thrm:selprin}. 
\item[(iii)] This result contains the semiclassical limit in the sense that, as long as $\vertiii{W^\hbar(0)-\rho^\hbar(0)}_{-M}=o(1)$,
\beq \label{eq:limitrewdq}
\lim\limits_{\hbar\to 0} \vertiii{W^\hbar(t)-\rho^\hbar(t)}_{-M}=0.
\ee
The selection of $\rho_0^\hbar$ so that the solution of classical problem \eqref{eqsert45} gives rise to a practical method is discussed in section \ref{secresopiijl}.
\end{itemize}
\end{remark}

Also, the estimates we develop  have a substantial impact on the treatment of nonlinear problems, where e.g. a priori bounds on $\|u(t)\|_{H^1}$ can be used through Sobolev embeddings to get  regularity  for the effective potential, $V=b|u|^p$.  An adaptation of  Lemma \ref{lm:coresemiclest} for the nonlinear Schr\"odinger equation can be found in \cite{atcMORE}.

\subsection{Statement of the Main Results: Numerical Investigation for ${\mathbf{ V(x)=-|x|}}$ }

Lemma \ref{lm:coresemiclest} illustrates very clearly why  the proof of Theorem \ref{thrm:main} cannot be extended  for $V \notin C^{1,a}$. Moreover, the numerical results discussed below indicate that there is one more assumption required to prove any version of Theorem \ref{thrm:main} when $a=0$. In that sense it seems that Theorem \ref{thrm:main} is sharp. 

On the other hand, the saddle point generated by $V(x)=-|x|$ is similar to that of $V(x)=-|x|^{1+a}$, $a\in (0,1)$ in many ways. If an estimate of the form \eqref{eq:dpexcluestform} was true, then a regularization similar to that described in Theorem \ref{thrm:selprin} would be possible. So a natural question arises: 
 {\em is it possible to outline the regime of validity of eq. \eqref{eq:limitrewdq}   for saddle points of the form $V(x)=-C|x|$?}

This question is answered positively, with the help of a numerical solver based on a posteriori error control. More specifically, for the numerical solution $\widetilde u^{\hbar}(t_n) \in L^2$ ($t^n$ being a discrete time level), it can be shown rigorously that we have an upper bound of the form 
 $$\|u^\hbar(t_n) - \widetilde u^\hbar(t_n)\|_{L^2} \leqslant E_n,$$  
 where $E_n$ is a computable quantity. The development of a  practical solver with a posteriori error control for the semiclassical Schr\"odinger equation with non-smooth potentials is a challenging task on its own; to the best of our knowledge the only available such solvers in the literature are \cite{KK,Dorf}. More details about the numerical method we use can be found in section \ref{sec:setup12233}.
Using this solver, it is possible to investigate reliably the behavior of $u^\hbar$, even though there is no a priori information  for the  behavior of the exact solution.  

\begin{remark} \upshape
Clearly one cannot investigate numerically the limit $\hbar\rightarrow 0$ by solving for particular small values of $\hbar$. However very often $\hbar \to 0$ is merely an approximation for concrete problems of the form \eqref{eq:schroeq} with $\hbar$  not  smaller than $10^{-4}$ \cite{T}. Thus, investigating the validity of our asymptotics for $\hbar \approx 10^{-2}$ to  $\hbar \approx 10^{-4}$ is interesting in itself -- in some cases more so than investigating the limit $\hbar\rightarrow 0$. Here we work for $\hbar \in [5\cdot 10^{-3}, 10^{-1}]$, having put the emphasis into ensuring stability and small error tolerance in the problems that we solve, rather than pushing computations for very small values of $\hbar$. In any case the qualitative behavior we observe seems to be quite robust for $\hbar\ll 1$; apparently it stabilizes  for $\hbar \approx 10^{-2}$. 
\end{remark}

The numerical results we obtain are described in some detail is section \ref{secresopiijl}.
Summarizing,  we note that the selection principle of Theorem \ref{thrm:selprin} appears to hold for a wavefunction interacting with the saddle point of $V(x)=-|x|$, under a \textit{non-interference condition}; for details see  definition \ref{defintrfffff}.  (The idea is that energy arriving to the saddle point in phase-space from many directions at the same time constitutes \textit{interference}.
) Singular wavepacket splitting cases can be successfully approximated; see Section \ref{sbs:noninr} and Appendix \ref{AppCwkb}. On the other hand, when { interference} takes place, we observe different behavior, with the semiclassical limit affected by \textit{quantum phase information}. Examples and quantitative aspects of this dependence on the phase are presented in section \ref{sec:interf11}.

These numerical results allow the formulation of a precise conjecture  for the validity of our semiclassical selection principle: Using the notations of Theorem \ref{thrm:main}, we propose that
\begin{conjecture*}
For potentials $V$ with localized singularities of the form $\pm C |x|$, 
\[
\lim\limits_{\hbar\to 0} \langle W^\hbar[u^\hbar(t)] - \rho^\hbar(t),\phi\rangle = 0 \qquad \forall \phi\in \mathcal{B}_{M},
\]
as long as there is no interference. 
\end{conjecture*} 

\begin{remark} \upshape
This can be seen as a refinement of the conditions derived by C. Fermanian-Kammerer, P. Gerard and C. Lasser. More specifically,  the assumption that ``the Wigner measure does not reach the set $S\setminus S^*$'', which appears in Theorem 2 of \cite{FGL}, can be refined to admit problems where the Wigner measure does interact with the singular saddle point (i.e. reaches $S\setminus S^*$) -- as long as there is no interference.
\end{remark}


\begin{remark}  \upshape
A  similar non-interference condition arises in \cite{Miller}. It is possible that the conjecture can be proved with methods similar to those used therein.
\end{remark}

The paper is organized as follows: in sections \ref{eqsewqas} and \ref{sec:lou}  we introduce preliminary material and some of the characteristics of the flow \eqref{eq:fflochar1} are presented. Section \ref{sec:selp} is devoted entirely in proving Theorem \ref{thrm:main} and stating the selection principle in Theorem \ref{thrm:selprin}. In section \ref{sec:setup12233} we describe the numerical method and some of its main characteristics. In Section  \ref{secresopiijl} we present  numerical results obtained in the case of non-smooth potentials and special attention is given to the interference and non-interference cases. Auxiliary and background material is presented in Appendices A, B and C.

\section{Phase-space methods for semiclassical asymptotics} \label{eqsewqas}

\subsection{The Wigner transform and Wigner measures}\label{sub:WT}

To study the semiclassical behavior of \eqref{eq:schroeq}, we use the Wigner transform. For a comprehensive introduction, as well as the state of the art for smooth potentials, one should consult the references \cite{LP,GMMP}. Here the aim is to present a brief but self-contained introduction.
For any $f\in L^2(\mathbb{R}^d)$,  its \textit{Wigner transform} (WT) is defined as
\beq
\label{eq:wt}
W^\hbar[f](x,k)=\int\limits_{y} e^{-2\pi i k y}f(x+\frac{\hbar y}2)\bar f(x-\frac{\hbar y}2)dy \quad \in \quad L^2(\mathbb{R}^{2d}).
\ee
This transform will be applied to the wavefunction $u^\hbar(t)$;
we will use the shorthand notations $W^\hbar(x,k,t)=W^\hbar[u^\hbar(t)](x,k)$, $W^\hbar_0=W^\hbar[u^\hbar_0]$ when there is no danger of confusion. In principle, the WT contains the same information as the original wavefunction, but unfolded in phase space, i.e. position-momentum space $\{(x,k)\} = \mathbb{R}^{2d}$. This physical information can be accessed through \textit{quadratic observables}: an operator valued observable $\mathbb{A}$ can be measured by
\cite{GMMP}
\beq\label{eq:tf}
A[u^\hbar](t)=\langle \mathbb{A} u^\hbar(t), u^\hbar(t) \rangle_x = \langle A_{\mathcal{W}}, W^\hbar[u^\hbar(t)]  \rangle_{x,k},
\ee
where $A_{\mathcal{W}}$ is the (semiclassically scaled) Weyl symbol of $\mathbb{A}$. More details on quadratic observables can be found in section \ref{sbs:compobs}.
Applying the WT to \eqref{eq:schroeq} we see that $W^\hbar(x,k,t)$  satisfies a well-posed equation in phase-space \cite{Mql}, namely
\begin{equation}
\label{eq1hsdvaol}
\begin{aligned}
\partial_t {{W}^\hbar}(x,k,t) +  2\pi k \cdot \partial_x {{W}^\hbar}(x,k,t) 
+ & i  \int{ e^{-2\pi i S y} \frac{V(x+\frac{\hbar}{2}y)-V(x-\frac{\hbar}{2}y)}{\hbar} dy  \,\, W^\hbar(x,k-S,t)dS}=0,\\
& {W}^\hbar(t=0)  =W^\hbar_0.
\end{aligned}
\end{equation}

The merit of the WT lies in its behavior as $\hbar \rightarrow 0$. The idea is that given a sequence of solutions of (\ref{eq:schroeq}), $\{u^{\hbar_n}(t)\}$ with ${\lim\limits_{n\rightarrow 0}\hbar_n=0}$, then (up to extraction of a subsequence) its  \textit{Wigner measure} (WM) $W^0(t)$ is defined as an appropriate weak-$*$ limit of
\beq
\label{eq:wm}
W^{\hbar}(t) \rightharpoonup W^0(t) \in \mathcal{M}^1_+(\mathbb{R}^{2d}),
\ee
where $\mathcal{M}^1_+(\mathbb{R}^{2d})$ are the probability measures on phase-space.
Moreover, the WM satisfies a Liouville equation,
\beq
\label{eq:wmeq}
\partial_t W^0(t) + 2\pi k\cdot \partial_x W^0(t) - \frac{1}{2\pi} \partial_x V \cdot \partial_k W^0(t)=0, \,\,\,\,\,\,\,\,\,\,\,\,\,\,
W^0(t=0)=W^0_0,
\ee
i.e. a formulation as in classical statistical mechanics. (See e.g. Th\'eor\`eme IV.1 of \cite{LP}, Section 7.1 of \cite{GMMP}, and Theorem \ref{thrm1} in Appendix \ref{sec:appA} of this paper.) For smooth potentials, problem \eqref{eq:wmeq}  can be efficiently solved, and its solution can be used to recover the macroscopic observables of the particle (and sometimes their probability densities; e.g., position and momentum densities). 
The Liouville equation \eqref{eq:wmeq} can be solved with the method of characteristics \cite{Evans}: consider the ODE for the characteristics (the classical trajectories),
\beq
\bac
\label{eq:char}
\dot X_{X_0,K_0}(t)=2\pi K_{X_0,K_0}(t), \,\,\,\,\,\,\,\,\,\,\, \dot K_{X_0,K_0}(t)=-\frac{1}{2\pi} \partial_xV(X_{X_0,K_0}(t)), \\
X_{X_0,K_0}(0)=X_0, \,\,\,\,\,\,\,\,\,\,\, K_{X_0,K_0}(0)=K_0.
\ea
\ee
Then a classical flow $\phi_t$ is defined in terms of
\beq
\phi_t(x,k)=(X_{x,k}(t),K_{x,k}(t)).
\ee
It is straightforward to see that the solution of \eqref{eq:wmeq} is given by
\[
W^0(t)=W^0_0 \circ \phi_{-t}.
\]
It is evident from \eqref{eq:char} why the regularity of $V\in C^{1,1}(\mathbb{R}^d)$ is a natural threshold for the validity of these types of results. For $V\in C^{1,1}(\mathbb{R}^d)$, the characteristics \eqref{eq:char} are well-defined for all $(X_0,K_0)\in\mathbb{R}^{2d}$, and thus the WM is unconditionally well defined at all times.  If $V\notin C^{1,1}(\mathbb{R}^d)$, then in general the Cauchy problem \eqref{eq:wmeq} is not well-posed over probability measures. 

In \cite{LP} it was further shown that for $V\in C^1(\mathbb{R}^d)$, and under appropriate additional technical assumptions, the WM does indeed satisfy \eqref{eq:wmeq}, which in general has multiple solutions. Thus,  the WM is one of the possible classical evolutions, but it is not known which one. In Remarque IV.3 of \cite{LP} a concrete example of $V\in C^1 \setminus C^{1,1}$, giving rise to a multivalued flow, is given -- namely the singular saddle point  $V(x)=-|x|^{1+a}$, $a\in(0,1)$. This ill-posedness is resolved in Theorem \ref{thrm:selprin}.
  
The class of potentials with conical singularities, $V\in W^{1,\infty}(\mathbb{R}^d)$,  e.g., $V(x)=-|x|\notin C^1$, arise as another natural threshold with respect to the regularity of flows. In particular, it is shown  that for potentials in $W^{1,\infty}(\mathbb{R}^d)$, the trajectories \eqref{eq:char} are well defined for \textit{almost all} initial data $(X_0,K_0)\in\mathbb{R}^{2d}$. On the level of the Liouville equation, this can be seen as well-posedness with initial data in $L^1\cap L^\infty$, \cite{AF,B}. 

\subsection{The smoothed Wigner transform}\label{sub:SWT}

It is well known that if $u^\hbar(x)$ exhibits oscillations at length-scales of $\hbar$, then $W^\hbar[u^\hbar](x,k)$ will exhibit oscillations at length-scales $\hbar$, and sometimes at smaller  scales as well. Therefore simply representing $W^\hbar$ numerically is prohibitively  expensive; and solving numerically \eqref{eq1hsdvaol} even more so. Very often a smoothed version of $W^\hbar$, $\widetilde{W}^\hbar = W^\hbar \ast G$, is used instead \cite{AMP,LP}; the motivation is that, if $A_{\mathcal{W}}$ is smooth enough,
\[
\langle \mathbb{A} u^\hbar(t), u^\hbar(t) \rangle = \langle A_{\mathcal{W}}, W^\hbar(t)  \rangle \approx  \langle A_{\mathcal{W}}, \widetilde{W}^\hbar(t)  \rangle.
\]
This can be made precise, i.e. in the limit the two transforms are equivalent \cite{LP},
\[
\lim\limits_{\hbar\rightarrow 0}   \langle \widetilde{W}^\hbar(t) -{W}^\hbar(t),\phi   \rangle =0  \,\,\,\,\,\,\,\,\, \forall \phi\in \mathcal{A}.
\]
(For the algebra of test functions in $\mathcal{A}$, see Appendix \ref{sec:appA}.)

In the remaining part of this paper,   we will denote $\widetilde{W}^\hbar(t)=\widetilde{W}^\hbar[u^\hbar(t)]$ the  \textit{smoothed Wigner transform} (SWT),  defined as
\begin{equation}
\label{eq:swtdef}
\widetilde{W}^\hbar(x,k)=\left({ \frac{2}{\hbar {\sigma_x \sigma_k}} }\right)^d \int\limits_{x,k} e^{-\frac{2\pi}{\hbar}\left[{ \frac{|x-x'|^2}{\sigma_x^2}+\frac{|k-k'|^2}{\sigma_k^2} }\right]}  W^\hbar (x',k')\,dx'\,dk'.
\end{equation}
Sometimes we will use the notation $\widetilde{W}^{\sigma_x,\sigma_k;\hbar}(x,k)$ when we want to denote explicitly the smoothing constants used.

For $\sigma_x \cdot\sigma_k \geqslant 1$ it can be shown that $\widetilde{W}^\hbar(x,k) \geqslant 0$ \cite{folland}. Often it is useful to use smaller values for the smoothing constants, $\sigma_x,\sigma_k<1$. In any case we will assume that the smoothing constants do not depend on $\hbar$, and are allowed to be in $\sigma_x,\sigma_k\in (0,1]$.
For more context on the SWT, including on the calibration of the smoothing parameters, see Appendix \ref{AppCoarse}.

\section{Flows with infinite Lyapunov exponents and loss of uniqueness} \label{sec:lou}

\begin{figure}[htb!]
\includegraphics[scale=0.4]{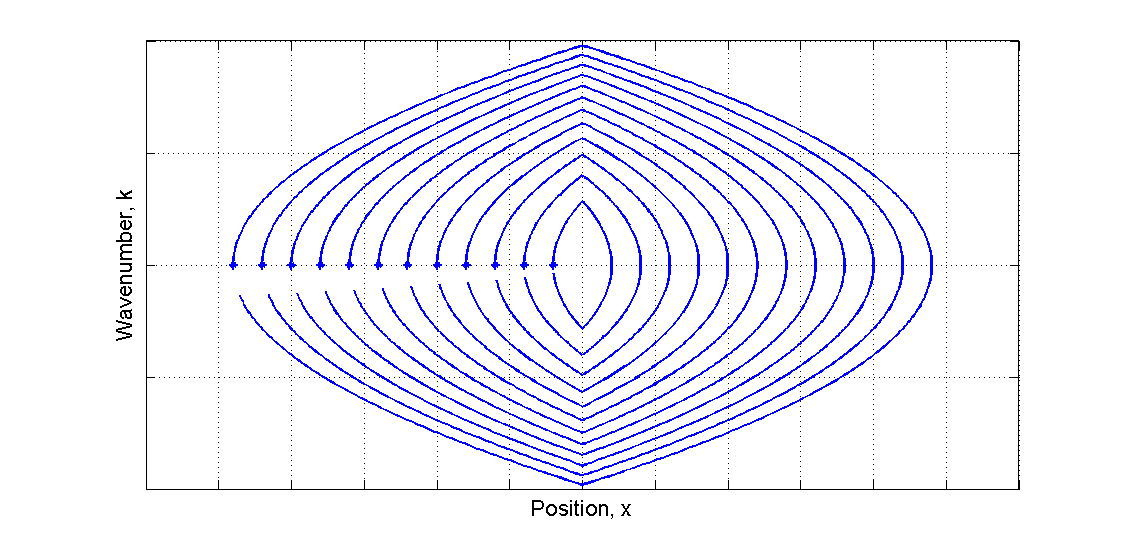}
\caption{For $V=|x|$, the characteristics in phase-space can be written explicitly (in each half space $\pm X >0$) as
$X(t)= -\frac{1}{2} sign(X)  t^2 + 2\pi K_0 t + X_0$,
$K(t)=-\frac{1}{2\pi} sign(X) t + K_0$.
 Here we see  plots of $(X(t),K(t))$ for various initial conditions $(X_0,0)$.   All trajectories eventually have corners, but all are uniquely defined. The picture is similar for $a \in (0,1)$.}
\label{Fig100}
\end{figure}

We investigate now, in some detail, the consequences of $V \notin C^{1,1}$. We consider the one-dimensional potentials
\beq \label{eq:vtheta}
V^\pm(x)=\pm |x|^{1+a}, \qquad a \in [0,1), \qquad x \in \mathbb{R}.
\ee
For $V^{+}=|x|^{1+a}$, the problem physically amounts to an oscillator. Although there is no strong solution of \eqref{eq:char} once the trajectory reaches $\{x=0\}$, by accepting weak solutions the flow is in fact well defined. Indeed, it is easy to check that the problem
\beq
\bac
\label{eq:char1}
\dot X(t)=2\pi K(t), \,\,\,\,\,\,\,\,\,\,\, \dot K(t)=-\frac{1}{2\pi} \partial_x V(X(t)), \\
X(0)=0, \,\,\,\,\,\,\,\,\,\,\, K(0)=K_0,
\ea
\ee
has a unique weak solution for all values of $K_0\in \mathbb{R}$. (See also
Figure \ref{Fig100}.)
So the impact of the conical singularity here is less smoothness of the trajectories,  but there is no loss of uniqueness. 
 
\begin{figure}[htb!]
\includegraphics[scale=0.45]{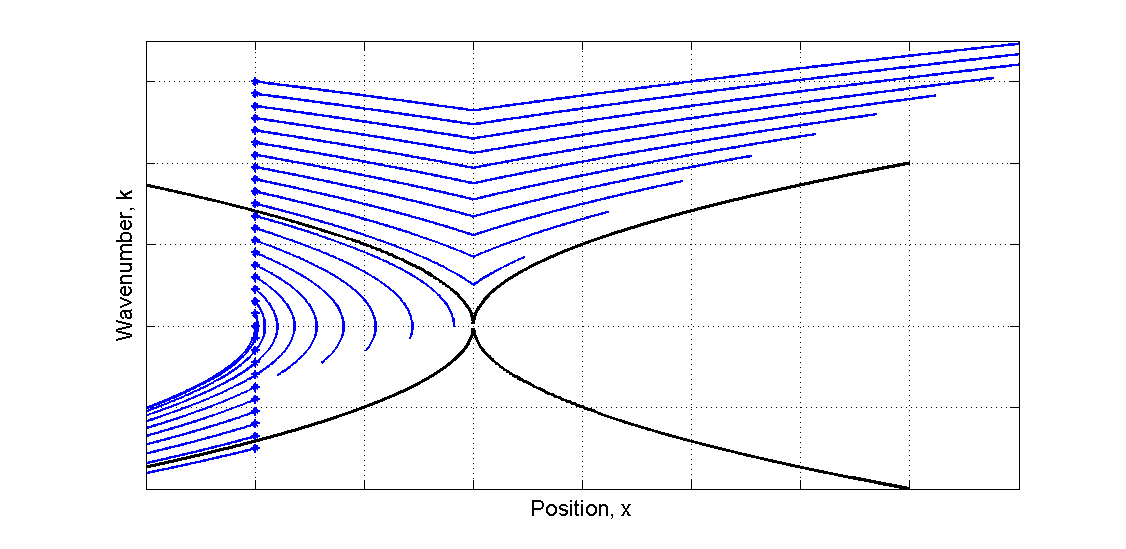}
\caption{For $V=-|x|$, the characteristics in phase-space can be written explicitly (in each half space $\pm X >0$) as
$X(t)= \frac{1}{2} sign(X)  t^2 + 2\pi K_0 t + X_0$,
$K(t)=\frac{1}{2\pi} sign(X) t + K_0$.
 Here we see  plots of $(X(t),K(t))$ for $t\in [0,T]$ and various initial conditions $(-1,K_0)$.   
  The separatrix $k=\pm\frac{1}{\pi}\sqrt{\frac{x}{2}}$ (shown in black) consists of two intersecting trajectories. Unlike regular saddle points, the two branches of the separatrix intersect in finite time over the fixed point $(0,0)$. The picture is similar for $a \in (0,1)$.}
\label{Fig101}
\end{figure}

For $V^{-}=-|x|^{1+a}$, the problem is a singular saddle point (Figure \ref{Fig101}). \textit{The trajectories that approach the fixed point $(x,k)=(0,0)$ now arrive in finite time, in contrast to what happens in regular saddle points.} Once they reach the fixed point, there is no unique continuation -- strong or weak. 
For example, let $a=0$ and consider the characteristic starting from $(X_0,K_0)=(1,-\frac{1}{\pi \sqrt{2}})$; then  each of the following 
\beq \label{eqwas101}
\bac
X(t)=\left\{ { 
\begin{array}{c} 
\frac{1}{2}t^2-\sqrt{2}t+1, \,\,\,\,\,\, t\leqslant \sqrt{2} \\
-\frac{1}{2}(t-\sqrt{2})^2, \,\,\,\,\,\, t >\sqrt{2}
\end{array}
}\right. \,\,\,\,\,\,\,\,\,\,\,\,\,\,\,\,\,
K(t)=\left\{ { 
\begin{array}{c} 
\frac{1}{2\pi}t-\frac{1}{\pi \sqrt{2}}, \,\,\,\,\,\, t\leqslant \sqrt{2} \\
-\frac{1}{2\pi}(t-\sqrt{2}), \,\,\,\,\,\, t >\sqrt{2}
\end{array}
}\right.
\ea
\ee
\beq \label{eqwas102}
\bac
\tilde{X}(t)=\left\{ { 
\begin{array}{c} 
\frac{1}{2}t^2-\sqrt{2}t+1, \,\,\,\,\,\, t\leqslant \sqrt{2} \\
\frac{1}{2}(t-\sqrt{2})^2, \,\,\,\,\,\, t >\sqrt{2}
\end{array}
}\right. \,\,\,\,\,\,\,\,\,\,\,\,\,\,\,\,\,
\tilde{K}(t)=\left\{ { 
\begin{array}{c} 
\frac{1}{2\pi}t-\frac{1}{\pi \sqrt{2}}, \,\,\,\,\,\, t\leqslant \sqrt{2} \\
\frac{1}{2\pi}(t-\sqrt{2}), \,\,\,\,\,\, t >\sqrt{2}
\end{array}
}\right.
\ea
\ee
\beq
\bac
\tilde{\tilde{X}}(t)=\left\{ { 
\begin{array}{c} 
\frac{1}{2}t^2-\sqrt{2}t+1, \,\,\,\,\,\, t\leqslant \sqrt{2} \\
0 , \,\,\,\,\,\, t >\sqrt{2}
\end{array}
}\right. \,\,\,\,\,\,\,\,\,\,\,\,\,\,\,\,\,
\tilde{\tilde{K}}(t)=\left\{ { 
\begin{array}{c} 
\frac{1}{2\pi}t-\frac{1}{\pi \sqrt{2}}, \,\,\,\,\,\, t\leqslant \sqrt{2} \\
0 , \,\,\,\,\,\, t >\sqrt{2}
\end{array}
}\right.
\ea
\ee 
are weak solutions of  \eqref{eq:char} past the interaction with the singularity. The respective explicit trajectories for $a\in (0,1)$ can be found in \cite{LP}. 

\textit{In other words, a classical particle with just enough momentum to reach this saddle point, can be scattered to the right, scattered to the left, or stay on the saddle point indefinitely -- or do combinations of the above.} There are genuinely different classical evolutions to choose from here. 
Moreover, if we take, e.g. a  particle starting at $(X_1(t),K_1(t))=\phi_t(-1,\frac{1}{\pi \sqrt{2}}+\delta)$, and $(X_2(t),K_2(t))=\phi_t(-1,\frac{1}{\pi \sqrt{2}}-\delta)$, we see that
\[
\forall t> \sqrt{2} \qquad
\lim\limits_{\delta \to 0} |(X_1(t),K_1(t)) - (X_2(t),K_2(t))| >0,
\]
therefore
\[
|(X_1(t),K_1(t)) - (X_2(t),K_2(t))| \leqslant e^{C_L t} 2\delta \qquad \mbox{ fails for every } C_L>0,
\]
see also Figure \ref{Fig101}.
\textit{In other words the flow, which is well-defined for almost all trajectories, has infinite Lyapunov exponent.}
 
With respect to more general potentials  with isolated conical singularities, i.e. 
\[
V(x)= V_0(x) + w(x)|g(x)|, \qquad V_0, w, g \mbox{ smooth},
\]
it was shown in \cite{FGL} that there is a similar behavior. Indeed,  away from the set $S=\{g(x)=0\}$, the potential is smooth, and therefore the characteristics are well-defined and smooth. As far as the semiclassical limit is concerned, if the WM is never supported on $S$, then the regular theory applies \cite{GMMP,FGL}. If however the WM arrives at some time on $S$, then the regular theory ceases to apply.
Moreover, the set $S$ should really be decomposed into the disjoint union $S=S_1 \cup S_0$, where
\beq
\label{eq:sets}
\bac
S_1= \{ g(x)=0 \mbox{ and } k\cdot \partial_x g(x) \neq 0 \}, \\
S_0=  \{ g(x)=0 \mbox{ and } k\cdot \partial_x g(x) = 0 \}.
\ea
\ee
If a characteristic  arrives at $S_1$,
its momentum will take it ``immediately'' out of $S$, and it will be continued uniquely -- with a corner, which essentially resembles the behavior observed for $V(x)=|x|$. However, if a characteristic  arrives at $S_0$, several classical evolutions are possible. 
The main result of \cite{FGL} is that, as long as the WM stays away from $S_0$, the uniquely defined flow indeed captures correctly its evolution. This motivates the following 

\begin{definition}  \label{definter} Consider a semiclassical family of problems \eqref{eq:schroeq}, with
\beq \label{eq:setconic}
V(x)= V_0(x) +  w(x)|g(x)|^{1+a}, \qquad V_0, w, g \mbox{ smooth}, \quad a \in [0,1).
\ee
We will say there is \textit{full interaction with the singularity of the flow}  if its WM reaches the set 
 \[
 S_0=  \{ g(x)=0 \mbox{ and } k\cdot \partial_x g(x) = 0 \}.
 \]
\end{definition}
So taking into account recent results, the study of ``full interaction with singular saddle points'' is in fact the natural generalization of the question put forth in in Remarque IV.3 of \cite{LP}.

\section{Proof of main Theorem} \label{sec:selp}
In this section we give the proof of the main result stated in Theorem \ref{thrm:main}. We begin by introducing notation, terminology and auxiliary technical results.  

\subsection{Notations} \label{sbs:not}

\begin{definition}[Fourier transform] \label{def:2} Given $f(x)$ its Fourier transform is defined as
\[
\widehat f (k):= \mathcal{F}_{x\to k}[f]= \int e^{-2\pi i x\cdot k} f(x) dx.
\]
For functions on phase-space $f(x,k)$, we will also use the Fourier transform in the 
second set of variables,
\[
\widehat f_2 (x,K):= \mathcal{F}_{k\to K}[f]= \mathcal{F}_2 f= \int e^{-2\pi i k\cdot K} f(x,k) dk.
\]
For the Fourier transform of functions of phase-space we will typically use the variables
\[
\widehat f (X,K):= \mathcal{F}_{x,k\to X,K}[f]= \mathcal{F} f= \int e^{-2\pi i [x\cdot X+ k\cdot K]} f(x,k) dxdk.
\]
\end{definition}

\begin{definition}[Schwarz test functions] We will denote by $\mathcal{S}(\mathbb{R}^d)$ the class of functions $\phi:\mathbb{R}^d\to\mathbb{C}$ for which
\[
\forall \quad \mbox{multi-indices} \quad a,b \quad  \exists \quad C_{a,b} \quad \mbox{so that} \quad |x^a \partial_x^b \phi(x)| \leqslant C_{a,b}.
\]	
As is well known, $\widehat{\phi} \in \mathcal{S}(\mathbb{R}^d) \Leftrightarrow \phi \in \mathcal{S}(\mathbb{R}^d)$.
\end{definition}

\begin{definition}
For $f\in C^{1,a}$ we define its norm by 
\[
\|f\|_{C^{1,a}} = \sup\limits_{x_1\neq x_2 } \frac{|\partial_x f(x_1) - \partial_x f(x_2)|}{|x_1-x_2|^a} + \sup\limits_{x} |f(x)| +\sup\limits_{x} |\partial_x f(x)| 
\]
\end{definition}

\begin{definition} For $a,b \in \mathbb{C}$, we will use the notation $a \preccurlyeq b$ with the understanding
\[
a \preccurlyeq b \quad \Leftrightarrow |a|=O(|b|).
\]
\end{definition}

\begin{definition}
Denote by $\mathcal{T}(t)$ the free-space propagator on phase-space,
\beq
\mathcal{T}(t) : f(x,k) \mapsto f(x-2\pi t k,k).
\ee

Given any function $f(x,k)$ on phase-space we denote for further reference
\begin{eqnarray}
T^{V}_\hbar f := &  \frac{2}{\hbar} Re \left[{ i\int{e^{2\pi i Sx}\widehat{V}(S)f(x,k-\frac{\hbar S}{2})dS} }\right] , \\
T^{V}_0 f := &  -\frac{1}{2\pi} \partial_x V \cdot \partial_k f , \\
\widehat T^{V}_\hbar f := &\mathcal{F} T^{V}_\hbar f= 2\int{ \widehat{V}(S) f (X-S,K) \frac{sin(\pi \hbar S\cdot K)}{\varepsilon} dS}, \\
\widehat{\mathcal{T}}(t)\widehat f (X,K) := & \mathcal{F} [\mathcal{T}(t)f] =\widehat{f}(X, K+2\pi X t),\\
\|f\|_{\mathcal{F}_j L^p}:= & \| \mathcal{F}_j f \|_{L^p}, \qquad \mbox{ for } j\in \{1,2,\emptyset \}. 
\end{eqnarray}
\end{definition}

\subsection{Auxiliary technical lemmas} \label{sbs:techndosc}

\begin{observation} \label{obs:109001} For any $f,g$ functions on phase-space, $t\in \mathbb{R}$ we have
\[
\langle \mathcal{T}(t)f,g\rangle = \langle f,\mathcal{T}(-t)g\rangle, \qquad \langle \widehat{\mathcal{T}}(t)f,g\rangle = \langle f,\widehat{\mathcal{T}}(-t)g\rangle,
\]
whenever the integrals exist. Moreover,
\[
\|\mathcal{T}(t) f\|_{\mathcal{F}L^p} = \| f\|_{\mathcal{F}L^p}, \qquad \|\widehat{\mathcal{T}}(t) f\|_{\mathcal{F}L^p} = \| f\|_{\mathcal{F}L^p}
\]
for all $p \in [1,\infty]$.
\end{observation}

\begin{observation} \label{obs:apriregW} If $W=W^\hbar[u](x,k)$ for some $u$ with $\|u\|_{L^2}=1$, then
\[
\widehat{W}_2 (x,K) = u(x-\frac{\hbar K}{2}) \overline {u}(x+\frac{\hbar K}{2}), \qquad \| \widehat{W}_2\|_{L^\infty_K L^1_x}=1
\]
In particular, it follows that
\[
\langle W , \phi \rangle \preccurlyeq \|\widehat{\phi}_2\|_{L^1_KL^\infty_x} \leqslant \| \phi\|_{\mathcal{F} L^1} \quad \Rightarrow \quad W \in \mathcal{B}_{-M}.
\]
\end{observation}

\medskip

\begin{lemma}[A specialized Liouville regularity estimate] \label{lm:splre}
Denote by $E(t)$ the propagator for the Liouville equation  \eqref{eqsert45},
 and assume that the potential $V$ satisfies
\[
\widehat{V}(S) |S| \in L^1.
\]
 Then there exists a constant $C_L>0$, depending on $V$ and $M$, so that
\[
\vertiii{\rho(t)}_M= \vertiii{E(t) \rho_0}_M \leqslant e^{C_L |t| } \vertiii{\rho_0}_M.
\]
\end{lemma}


\begin{remark} \upshape
 The point of this result is that smoothness in $k$ is in fact preserved by the flow. In that sense, this lemma can be seen as a counterpart of Proposition 1 of \cite{FGL}.
\end{remark}

Before we proceed to the proof , observe that
the assumption $\widehat{V}(S) |S| \in L^1$ is a little stronger than $V \in W^{1,\infty}$, but less strong than $V\in C^{1,1}$. 
For example if $a\in (0,1)$,  $V(x)=C|x|^{1+a} b(x)$ (where $b(x)$ is a smooth cutoff function of compact support, see, e.g., the statement of Problem \ref{prblm:1} in section \ref{sbs:pssttl}) then $\widehat{V}(S) |S| \in L^1$ while $V \notin  C^{1,1}$. Such potentials in particular are the ones that appear in the example of Remarque IV.3 of \cite{LP}. 

\medskip

\noindent {\bf Proof :} Without loss of generality, we  prove the result for $t>0$. 
We  work in the Fourier domain,
\[
\partial_t \widehat \rho -2\pi X\cdot \partial_K \widehat \rho + 2\pi \int\limits_S \widehat{V}(S) \widehat{\rho}(X-S,K) S \cdot K dS=0,
\]
for some $\rho_0 \in \mathcal{S}\cap \mathcal{B}_M$.
Now using $\mathcal{T}(t)$ and integrating in time, we can recast the Fourier-transformed Liouville equation  in mild form, namely
\[
\widehat{\rho}(t) = \widehat{\mathcal{T}}(t) \rho_0 + 2\pi \int\limits_{\tau=0} \widehat{\mathcal{T}}(t-\tau) \int\limits_S S\cdot K \widehat{V}(S) \widehat{\rho}(X-S,K,\tau) dS \,\, d\tau.
\]
By virtue of the Young inequality and observation \ref{obs:109001} it follows that
\[
\| \widehat {\rho}(t)\|_{L^1_{X,K}} \leqslant \|\widehat{\rho}_0 \|_{L^1} + \|\widehat{V}(S) \, |S|\,\|_{L^1_S} \int\limits_{\tau=0}^t \|\,|K| \, \widehat{\rho}(\tau) \|_{L^1_{X,K}};
\]
similarly, for any $m \in \mathbb{N}$,
\[
\bac
\| \,|K|^m \widehat {\rho}(t)\|_{L^1_{X,K}} \leqslant \|\, |K|^m\widehat{\rho}_0 \|_{L^1} + \|\widehat{V}(S) \, |S|\,\|_{L^1_S} \int\limits_{\tau=0}^t \|\,|K|^{m+1} \, \widehat{\rho}(\tau) \|_{L^1_{X,K}} \quad \Rightarrow \\

\Rightarrow \quad \sum\limits_{m} M^{-m} \| \,|K|^m \widehat {\rho}(t)\|_{L^1} \leqslant \sum\limits_m \left({  M^{-m}\|\, |K|^m\widehat{\rho}_0 \|_{L^1} + M\|\widehat{V}(S) \, |S|\,\|_{L^1} \int\limits_{\tau=0}^t M^{-m-1} \|\,|K|^{m+1} \, \widehat{\rho}(\tau) \|_{L^1} }\right) \quad \Rightarrow \\

\Rightarrow \quad \vertiii{\rho(t)} \leqslant \vertiii{\rho_0} + \|\widehat{V}(S) \, |S|\,\|_{L^1_S} M \int\limits_{\tau=0}^t \vertiii{\rho(\tau)} d\tau.
\ea
\]
The result follows by virtue of the Gronwall inequality, and by  density of $\mathcal{S}\cap \mathcal{B}_M$ in the space $\mathcal{B}_M$.
\qed

\begin{observation} \label{obs:Eadj} Denote by $E(t)$ the propagator of the Liouville equation (as in Lemma \ref{lm:splre} above), and $\widehat{E}(t)$ its Fourier transform, $\widehat{E}(t):\widehat{f} \mapsto \widehat{E(t)f}$. Then, for any $f,g$ functions on phase-space, $t\in \mathbb{R}$,
\[
\langle E(t)f,g\rangle = \langle f,E(-t)g\rangle, \qquad \langle \widehat{E}(t)f,g\rangle = \langle f,\widehat{E}(-t)g\rangle,
\]
whenever the integrals exist. In particular, $E(t)$ makes sense on functions belonging to $\mathcal{B}_{-M}$.
\end{observation}

\begin{lemma} \label{lm:coresemiclest}
Let $V\in C^{1,a}(\mathbb{R}^d,\mathbb{R})$, $W=W^\hbar[u]$ for some $\|u\|_{L^2}=1$, and $\phi$ be a test function on phase-space regular enough for the integrals below to exist. Then
\beq \label{eq:lm38a}
 \langle (T_\hbar^V - T^V_0) W, \phi \rangle_{x,k}  = \langle \int\limits_{s=-\frac{K}2}^{\frac{K}{2}} (\partial_x V(x+\hbar s)- \partial_x V(x))\cdot ds  \widehat{W}_2^\hbar, \widehat{\phi}_2 \rangle_{x,K} \preccurlyeq \hbar^a \| V\|_{C^{1,a}} \langle |\widehat{W}_2^\hbar|,  |K|^2 |\widehat{\phi}_2| \rangle_{x,K}.
\ee
Therefore, by straightforward application of observation \ref{obs:apriregW},
\beq \label{eq:lm38b}
 \langle (T_\hbar^V - T^V_0) W, \phi \rangle_{x,k} \preccurlyeq \hbar^a \| V\|_{C^{1,a}} M^2 \vertiii{ \phi}_M .
\ee
\end{lemma}

\noindent {\bf Proof:} Observe that
\[
\mathcal{F}_2 [T^{V}_\hbar W] = \frac{i}{\hbar}\int\limits_{k,S} e^{-2\pi i k\cdot K} \left[{
e^{2\pi i Sx}\widehat{V}(S)W(x,k-\frac{\hbar S}{2})-e^{-2\pi i Sx}\overline{\widehat{V}}(S)\overline{W}(x,k-\frac{\hbar S}{2})
}\right] dSdk.
\]
Since the Wigner transform is real valued, $\overline W = W$; since the potential is real valued $\overline{\widehat{V}}(S) = \widehat{V}(-S)$. Therefore,
\beq \label{eq:lsps123a}
\bac
\mathcal{F}_2 [T^{V}_\hbar W] = \frac{i}{\hbar}\int\limits_{k,S} e^{-2\pi i k\cdot K} \left[{
e^{2\pi i Sx}\widehat{V}(S)W(x,k-\frac{\hbar S}{2})-e^{-2\pi i Sx}{\widehat{V}}(-S) {W}(x,k-\frac{\hbar S}{2})
}\right] dSdk=\\
=\frac{i}{\hbar}\int\limits_{k,S} e^{-2\pi i [k\cdot K -S\cdot x]} \
\widehat{V}(S)W(x,k-\frac{\hbar S}{2})
 dSdk
-\frac{i}{\hbar}\int\limits_{k,S} e^{-2\pi i [k\cdot K +S\cdot x]} \
\widehat{V}(-S)W(x,k-\frac{\hbar S}{2})
 dSdk
=\\
=\frac{i}{\hbar}\int\limits_{k,S} e^{-2\pi i [k\cdot K -S\cdot x]} \
\widehat{V}(S)W(x,k-\frac{\hbar S}{2})
 dSdk
-\frac{i}{\hbar}\int\limits_{k,S} e^{-2\pi i [k\cdot K -S\cdot x]} \
\widehat{V}(S)W(x,k+\frac{\hbar S}{2})
 dSdk
=\\
=\frac{i}{\hbar}\int\limits_{k,S} e^{-2\pi i [k\cdot K -S\cdot x]} \
\widehat{V}(S) \left[{ W(x,k-\frac{\hbar S}{2}) -W(x,k+\frac{\hbar S}{2})}\right]
 dSdk
=\\
=i \widehat W_2(x,K) \int\limits_{S} e^{2\pi i S\cdot x} \
\widehat{V}(S)  \frac{ e^{-2\pi i \frac{\hbar K}{2}S} - e^{2\pi i \frac{\hbar K}{2}S} }{\hbar}
 dS
=i \widehat W_2(x,K) \frac{V(x-\frac{\hbar K}{2})-V(x+\frac{\hbar K}{2})}{\hbar}
\ea
\ee
On the other hand, it is trivial to check that
\beq \label{eq:lsps123}
\mathcal{F}_2 [T^{V}_0 W] = -i \partial_x V(x) \cdot K \mathcal{F}_2 W.
\ee
Thus combining equations \eqref{eq:lsps123}, \eqref{eq:lsps123a} it follows that
\[
\bac
\mathcal{F}_2 \left[{ (T_\hbar^V - T^V_0) W }\right]= i \widehat W_2(x,K) \frac{V(x-\frac{\hbar K}{2})-V(x+\frac{\hbar K}{2})}{\hbar} + i \partial_x V(x) \cdot K \mathcal{F}_2 W =\\
= i \widehat W_2(x,K) \int\limits_{s=-\frac{K}2}^{\frac{K}2} \left({ \partial_x V(x+\hbar s) - \partial_x V(x) }\right) \cdot ds \Rightarrow  \\
{ } \\

\Rightarrow \left|{ \langle (T_\hbar^V - T^V_0) W, \phi \rangle_{x,k} }\right| = \left|{ \langle \mathcal{F}_2(T_\hbar^V - T^V_0) W, \mathcal{F}_2\phi \rangle_{x,K} }\right| \leqslant \\

\leqslant \|V\|_{C^{1,a}} \hbar^a \langle \int\limits_{s=-\frac{K}2}^{\frac{K}2} |s|^a ds |\widehat W_2| , |\widehat {\phi}_2|\rangle_{x,K} \preccurlyeq \hbar^a \| V\|_{C^{1,a}} \langle |\widehat{W}_2^\hbar|, |K|^{a+1} |\widehat{\phi}_2| \rangle_{x,K}.
\ea
\]

The proof of eq. \eqref{eq:lm38a} is complete, since without loss of generality $a\leqslant 1$. For eq. \eqref{eq:lm38b} it suffices to observe that 
\[
\|K^2 \widehat{\phi}_2(x,K)\|_{L^1_KL^\infty_x} \leqslant \|K^2 \widehat{\phi}(X,K)\|_{L^1_{X,K}} \leqslant M^2 \vertiii{\phi}_M
\]
\qed

\begin{lemma} \label{lm:technV} Let 
\[
V(x)=C |x|^{1+a} b(x), \qquad a\in (0,1)
\]
where $b$ is a smooth cutoff function as in the statement of Problem \ref{prblm:1} (see Section \ref{sbs:pssttl}). It can be seen that
\[
\int\limits_{S} |\widehat{V}(S)|  \, |S| dS < \infty.
\]
\end{lemma}

\noindent {\bf Proof: } By observation,
\[
\| \widehat{V} \|_{L^\infty} \leqslant \|V\|_{L^1} < \infty.
\]
Moreover, observe that for $a\in(0,1)$
\[
\mathcal{F}_{x\to k} [|x|^{1+a}] =  C_{d,a} |k|^{-a-1-d}
\]
in the sense of distributions \cite{GelfSh}. It follows that
\[
|k|>1 \quad \Rightarrow \quad |\widehat{V}(k)| \leqslant C |k|^{-1-a-d}.
\]
The result follows by observing 
\[
\int \limits_S {|\widehat{V}(S)| \, |S| dS} \leqslant C \| \widehat{V} \|_{L^\infty}+ C \int\limits_{|S|>1} |k|^{-a-d} dk =
C \| \widehat{V} \|_{L^\infty}+ C \int\limits_{\rho=1}^{+\infty} \rho^{-a-d} \rho^{d-1} d\rho < \infty.
\]
\qed

We are now ready to give the proof of  the main theorem.  

\subsection{Proof of Theorem \ref{thrm:main}} \label{sbs:prmth}

Denote 
\[
h(t):= W^\hbar(t) - \rho^\hbar(t);
\]
by subtracting eq. \eqref{eqsert45} from \eqref{eq1hsdvaol} we find
\[
\bac
\partial_t h(t) +2\pi k \cdot \partial_x h + T_0^V h = (T_0^V-T^V_\hbar)W^\hbar \quad \Rightarrow \\
\Rightarrow \quad h(t) = E(t)(W^\hbar_0-\rho^\hbar_0) + \int\limits_{\tau=0}^t E(t-\tau) \big( (T_0^V-T^V_\hbar)W^\hbar(\tau)\big) d\tau \quad \Rightarrow\\
\Rightarrow \quad
\left|{ \langle h(t), \phi \rangle  }\right| \leqslant \left|{ \langle W^\hbar_0-\rho^\hbar_0, E(-t)\phi \rangle }\right|+ \int\limits_{\tau=0}^t \left|{\langle  \big( (T_0^V-T^V_\hbar)W^\hbar(\tau)\big), E(\tau-t) \phi \rangle }\right| d\tau ,
\ea
\]
where in the first step we used Duhamel's principle to take advantage of the propagator $E$ of the Liouville equation (see Lemma \ref{lm:splre}), and in the second step we used Observation \ref{obs:Eadj}.
At this point using Lemmas \ref{lm:splre}, \ref{lm:coresemiclest} and observation \ref{obs:apriregW} it follows that
\beq
\begin{aligned}
\left|{ \langle h(t), \phi \rangle  }\right| & \leqslant  \left|{ \langle W^\hbar_0-\rho^\hbar_0, E(-t)\phi \rangle }\right| + \hbar^a \| V\|_{C^{1,a}} M^2 \int\limits_{\tau=0}^t  \vertiii{E(\tau-t)\phi}_M d\tau  \leqslant \\
& \leqslant e^{C_L t} \vertiii{W^\hbar_0-\rho^\hbar_0}_{-M} \vertiii{\phi}_M + \hbar^a \| V\|_{C^{1,a}} M^2 \vertiii{\phi}_M \int\limits_{\tau=0}^t  e^{C_L \tau } d\tau.
\end{aligned}
\ee

By a simple estimate of the $d\tau$ integral, it follows that
\beq \label{eq:finspelestw}
\left|{ \langle h(t), \phi \rangle  }\right|  \leqslant C e^{C_L t} \Big( \vertiii{W^\hbar_0-\rho^\hbar_0}_{-M} + \hbar^a \Big) \vertiii{\phi}_M
\ee
for some constant $C$ independent of $\hbar$. \qed

\begin{observation} \label{obs:lstobsggt}
 One strategy is to choose $W_0^\hbar = \rho_0^\hbar$, in which case the first term on the rhs of eq. \eqref{eq:finspelestw} simply drops out. As we already discussed,  in many cases it is desirable that $\rho^\hbar_0$ is a smoothed version of $W^\hbar_0$, so that the interference terms are suppressed (see section \ref{sub:SWT}, Appendix \ref{AppCoarse}). The requirement
$
\vertiii{W^\hbar_0 - \rho^\hbar_0}_{-M} =o(1)
$
prescribes a particular family of smoothing strategies that can be used and still be covered by Theorem \ref{thrm:main}. For an  explicit example of smoothing so that
\beq \label{eq:vertiiconstr}
\vertiii{W^\hbar_0 - \widetilde{W}^\hbar_0}_{-M} =o(1)
\ee
see Lemma \ref{lm:Msmstr}.
\end{observation}

\subsection{Passage to the limit $\hbar \to 0$} \label{sbs:pssttl} With Theorem \ref{thrm:main} at hand, we can now proceed to resolve the following

\begin{problem} \label{prblm:1} Let $b(x) \in \mathcal{S}(\mathbb{R}^d)$ be a cutoff function with $b(x)=1$ for $|x|<2L$, $b(x)=0$ for $|x|>4L$. Let
\[
V(x)=-|x|^{1+a}b(x),  \qquad a \in [0,1),
\]
and  $\{ u_0^\hbar\}_{\hbar \in (0,1)} $ be a family of initial data for \eqref{eq:schroeq} so that
\beq \label{eq:iddwltua}
W^0_0 = w-* \lim\limits_{\hbar \to 0} W^\hbar[u^\hbar_0] \mbox{ exists, } \quad supp \, W^0_0 \cap \{ \frac{1}{2}(2\pi k)^2-|x|^{1+a} b(x) = 0 \} \neq \emptyset.
\ee

Compute $W^0(t)=w-* \lim\limits_{\hbar \to 0} W^\hbar[u^\hbar(t)]$, or show that it is not well-defined.
\end{problem}

\begin{remark} \upshape Some clarifications are in order:

\begin{enumerate}
\item[(i)] The cutoff $b$ does not play any substantial role, and is included only for technical reasons. With out loss of generality we will assume that $L$ is large enough so that it doesn't affect our computations.

\item[(ii)] The $w-*$ limit is taken with respect to the test functions $\mathcal{B}_M$, e.g, $\langle W^0_0, \phi \rangle = \lim\limits_{\hbar \to 0} \langle W^\hbar[u^\hbar_0], \phi \rangle$.

\item[(iii)] It follows from eq. \eqref{eq:iddwltua} that problem \eqref{eq:wmeq} for the evolution of the WM in time has multiple weak solutions. If $a \in (0,1)$, it is known by \cite{LP} that $W^0(t)$ is one of these. In that case, we need to compute the selection principle, i.e. a practical criterion to select the correct one. If $a=0$, it is not  known rigorously whether  $W^0(t)$ is related to some appropriate weak solution of \eqref{eq:wmeq}.

\item[(iv)] In view of Theorem \ref{thrm:main},
\[
\lim\limits_{\hbar \to 0} \langle W^\hbar[u^\hbar(t)], \phi \rangle = \lim\limits_{\hbar \to 0} \langle \rho^\hbar(t), \phi \rangle.
\]
So the question is simplified to (the ``purely classical'') computation of the concentration limit $\hbar \to 0$ of \eqref{eqsert45}. To fix ideas, we will work with $\rho^\hbar_0 = W^\hbar_0$. The only assumption of Theorem \ref{thrm:main} which is not obviously satisfied, is that $\widehat{V}(S) |S| \in L^1$. For that, we refer to Lemma \ref{lm:technV}.
\end{enumerate}
\end{remark}

\begin{theorem}[Selection principle] \label{thrm:selprin} Assume we are in the setting of Problem \ref{prblm:1}, and in addition $a\in (0,1)$. Denote
\beq \label{eq:prtpps}
\begin{aligned}
S^+&=\{(x,k)\in\mathbb{R}^2 \,\, | \,\, H(x,k)=\frac{1}{2}(2\pi k)^2-|x|^{1+a}b(x) > 0\}, \\
S^-&=\{(x,k)\in\mathbb{R}^2 \,\, | \,\, H(x,k)=\frac{1}{2}(2\pi k)^2-|x|^{1+a}b(x) < 0\}, \\
S&=\{(x,k)\in\mathbb{R}^2 \,\, | \,\, H(x,k)=\frac{1}{2}(2\pi k)^2-|x|^{1+a}b(x) = 0\}.
\end{aligned}
\ee
\begin{itemize}
\item[(i)] Denote by $\chi_\Omega$ the indicator function for the domain $\Omega \subseteq \mathbb{R}^{2d}$. Then
 $\lim\limits_{\hbar \to 0} W^\hbar[u^\hbar(t)]$ is well defined if and only if 
$w-* \lim\limits_{\hbar \to 0} W^\hbar_0 \chi_{S^+} $  exists. 
\item[(ii)]
Moreover, if (i) holds, $w-* \lim\limits_{\hbar \to 0} W^\hbar_0 \chi_{S^-}$  also exists, and
\beq \label{eq:concsel}
w-* \lim\limits_{\hbar \to 0} W^\hbar[u^\hbar(t)] =  w-* \lim\limits_{\hbar \to 0} W^\hbar_0 \chi_{S^-} \circ \phi_{-t} +
w-* \lim\limits_{\hbar \to 0} W^\hbar_0 \chi_{S^+} \circ \phi_{-t}.
\ee
\end{itemize}
\end{theorem}

\begin{remark} \upshape Some technical clarifications:
\begin{enumerate}
\item[(i)] It is explained in detail below  how eq. \eqref{eq:concsel} makes sense. For a simple explicit case, see also example \ref{ex:1selccc}.
\item[(ii)] In \cite{LP} it is discussed in some detail how $w-* \lim\limits_{\hbar \to 0} W^\hbar_0 \chi_{S^+} $ may fail to exist.
\item[(iii)] We will make some standard additional regularity assumptions on the initial data, namely
\[
\forall \hbar, \quad \max\limits_{|a|, |b| \leqslant d+1} \| x^a \partial_x^b u^\hbar_0(x) \|_{L^2} < \infty.
\]
This is sufficient to assure that $\rho^\hbar_0=W^\hbar_0 \in L^1(\mathbb{R}^{2d}) \cap L^\infty(\mathbb{R}^{2d})$; see Theorem \ref{thrm:wigregureg}. Although this is not strictly necessary, it simplifies considerably some technical points in the proof below.
\end{enumerate}
\end{remark}

\medskip

\noindent {\bf Proof of Theorem \ref{thrm:selprin}: } Phase space was partitioned into the disjoint union $\mathbb{R}^2 = S^+ \cup S \cup S^-$ in eq. \eqref{eq:prtpps}, accordingly
\[
\bac
\rho^\hbar_0 = \rho^\hbar_+ + \rho^\hbar_- + \rho^\hbar_z, \qquad \mbox{ where } \quad
\rho^\hbar_\pm= \rho^\hbar_0 \chi_{S^\pm}, \quad \rho^\hbar_z= \rho^\hbar_0 \chi_{S}.
\ea
\]
Since $\rho^\hbar_0 \in L^1 \cap L^\infty$ and $S$ is of measure zero, it follows automatically that $\rho^\hbar_z=0$. 
Therefore
$\rho^\hbar_0 = \rho^\hbar_+ + \rho^\hbar_-,$ and it suffices to solve each of the problems 
\begin{eqnarray}
\label{eqsert45a}
\partial_t \rho_+^\hbar(t) + 2\pi k\cdot \partial_x \rho_+^\hbar(t) - \frac{1}{2\pi} \partial_x V \cdot \partial_k \rho_+^\hbar(t)=0, \,\,\,\,\,\,\,\,\,\,\,\,\,\,
\rho^\hbar(t=0)=\rho^\hbar_+, \\
\label{eqsert45b}
\partial_t \rho_-^\hbar(t) + 2\pi k\cdot \partial_x \rho_-^\hbar(t) - \frac{1}{2\pi} \partial_x V \cdot \partial_k \rho_-^\hbar(t)=0, \,\,\,\,\,\,\,\,\,\,\,\,\,\,
\rho^\hbar(t=0)=\rho^\hbar_-,
\end{eqnarray}
separately.
The point, of course, is that by construction,
$\rho^\hbar_\pm \circ \phi_{-t}$ stays supported inside $S^\pm$,  $\forall t \in \mathbb{R}, \,\, \hbar \in (0,1)$. 
The restriction of the flow $\phi_t$ on each of the sets $S^+$, $S^-$, will be denoted by $\phi^+_t$, $\phi^-_t$ respectively.

\begin{claim} The flow $\phi^\pm_t$  is well-defined and continuous, i.e.
\[
\phi^\pm_t \in C(S^\pm,S^\pm) \qquad \forall t>0.
\]
\end{claim}

\noindent {\bf Proof of the claim: } It follows in exactly the same way as Proposition 1 of \cite{FGL}.  
\qed

Therefore each of $\phi^\pm_t$ can be extended to the closure of its domain. By abuse of notation (but without real danger of confusion) we will denote this extension as $\phi^\pm_t$
\[
\phi^\pm_t \in C(\overline{S^\pm},\overline{S^\pm}) \qquad \forall t>0.
\]

So in solving each of \eqref{eqsert45a}, \eqref{eqsert45b} we will work exclusively on the respective domains $\overline{S^\pm}$. Thus, for $f \in \mathcal{S}(\overline{S^\pm})$
\beq \label{eq:clouselppr}
 \langle \rho^\hbar_\pm \circ \phi^\pm_t, f \rangle  =  \langle \rho^\hbar_\pm  , f \circ \phi^\pm_t \rangle 
\ee

To conclude we observe that since
$\rho^\hbar = \rho^\hbar_+ + \rho^\hbar_-$, and we know that $w-* \lim\limits_{\hbar} \rho^\hbar$ exists, then necessarily $w-* \lim\limits_{\hbar} \rho^\hbar_-$ exists if and only if $w-* \lim\limits_{\hbar} \rho^\hbar_+$ exists. In that  case, and observing that
  $f \circ \phi^\pm_t$ stays continuous for all times,
\[
\lim\limits_{\hbar \to 0}  \langle \rho^\hbar_\pm(t), f \rangle  =  \langle \lim\limits_{\hbar \to 0} \rho^\hbar_\pm  , f \circ \phi^\pm_t \rangle 
\]

For the same reason, if $w-* \lim\limits_{\hbar} \rho^\hbar_+$ doesn't exist, we cannot pass to the $\hbar \to 0$ limit of eq. \eqref{eq:clouselppr}.

The proof of Theorem \ref{thrm:selprin} is  now complete.

\qed

 It is clear that if Theorem \ref{thrm:main} was valid for $a=0$, then Theorem \ref{thrm:selprin} would follow, with the same proof. This is the motivation behind the numerical investigation of its validity for $a=0$ which follows in the next sections.

\medskip
Now, let us look at a concrete example:
\begin{example} \label{ex:1selccc} Assume we are in the setting of Theorem \ref{thrm:selprin}, we take $d=1$ and let 
\[
u_0^\hbar = \hbar^{-\frac{1}4} e^{-\frac{\pi}{2} \left({ \frac{x-x_0}{\sqrt{\hbar}} }\right)^2 + i m(\hbar) \frac{ \sqrt{|x_0|^{1+a}} (x-x_0) }{\hbar}}
\]
for some $x_0<0$.
If $\lim\limits_{\hbar \to 0} m(\hbar) = 1$, then the WM is
\[
W^0_0=w-* \lim\limits_{\hbar\to 0} W^\hbar[u^\hbar_0] = \delta(x-x_0,k-\frac{\sqrt{2 |x_0|^{1+a}}}{2\pi})
\]
is supported on the separatrix $S$. 
Since $W^\hbar_0$ is a Gaussian in phase space with an effective support of $O(\hbar^{\frac{1}2})$, it follows that
if e.g., $m(h) = 1 + \hbar^{\frac{1}6} sin(\frac{1}\hbar)$, then $w-* \lim\limits_{\hbar} \rho^\hbar_\pm$ doesn't exist, since the mass of $W^\hbar_0$ oscillates between $S^+$ and $S^-$. If, on the other hand, $m(h) = 1 + \hbar^{6} sin(\frac{1}\hbar)$, then the oscillations would be negligible in the limit, and $w-* \lim\limits_{\hbar} \rho^\hbar_\pm$ exists.
\end{example}

Related examples can also be found in \cite{AP,AP2}.

\section{The numerical method} \label{sec:setup12233}

\subsection{Solving the semiclassical Schr\"odinger equation with conical singularities}\label{sec:setup}
The numerical solution of \eqref{eq:schroeq} is complicated from the theoretical as well as from the practical point of view. The main difficulty is that the solution of \eqref{eq:schroeq} oscillates with wavelength $\mathcal{O}(\hbar)$ thus standard numerical methods require very fine meshes (space and time)  to resolve adequately this high oscillatory behavior. Further the solution might exhibit caustics, making its numerical approximation even more difficult. Finally the relatively low smoothness of the potential $V$ means that several tools widely used in the numerical analysis and simulation of such problems are now not available. 

{ Popular methods for the numerical solution of \eqref{eq:schroeq} are time-splitting spectral methods and Crank-Nicolson finite element / finite difference methods. The standard Crank-Nicolson finite element / finite difference methods} suffer from a very restrictive dispersive relation, cf. \cite{JMS}, connecting the space and time mesh sizes with the parameter $\hbar$ thus requiring considerable computational resources  in order to produce accurate solutions for $\hbar\ll 1$.  In an attempt to relax this restrictive dispersive relation Bao, Jin and Markowich in \cite{BJM} proposed  time-splitting spectral methods for the numerical solution of \eqref{eq:schroeq}. This is widely considered to be the preferred approach for semiclassical problems; however it requires $V \in C^2$ at least for any kind of rigorous convergence result.  

A different approach to overcome this difficulty is based on adaptivity. 
Adaptive methods are widely used  in recent years  to construct accurate numerical approximations to  a broad class of problems { with substantially reduced computational cost by creating appropriately nonuniform meshes in space and time.}
There are several ways to { propose an adaptive strategy. One such approach  is based on rigorous a posteriori error control. The idea is to estimate the error in some natural norm by}
\begin{equation}
\label{GApost}
\|u-U\|\le \mathcal{E}(U) 
\end{equation}
where $\mathcal{E}(U)$ a computable quantity  depending on the approximate solution $U$ and the data of the problem. { A crucial property that the estimator $\mathcal{E}(U)$ must satisfy, is to converge with the same order as the numerical method. It is then said that $\mathcal{E}(U)$ decreases with optimal order with respect to the mesh discretization parameters.}
%
%
The existing literature on adaptive methods based on a posteriori error bounds for the numerical approximation of \eqref{eq:schroeq} is very limited. 
Very recently the authors  presented in \cite{KK}, an adaptive algorithm for the numerical approximation of \eqref{eq:schroeq},  based on a posteriori error estimates of optimal order. The proposed adaptive method
%
proved to be competitive with the best available methods in the literature not only for the approximation of the solution of \eqref{eq:schroeq} but as well as for its observables, c.f. \cite{KK}.  

Here we want to investigate the behavior of a quantum problem, for which we don't have even any qualitative a priori information. (E.g. the percentage of mass scattered in different directions after the interaction with the singularity.) Hence a posteriori error control is particularly useful, as it provides a rigorous, quantitative grasp on the quantum interaction -- making meaningful the subsequent comparison to the classical asymptotics.

\subsection{The CNFE method}\label{sec:syn}
{ In \cite{KK} the authors consider the initial-and-boundary value problem
\begin{equation}
\label{semiclassical} \left \{
\begin{aligned}
&i\hbar u_t^{\hbar}+\frac{\hbar^2}2\varDelta u^{\hbar}-Vu^{\hbar} = f    &&\quad\mbox{in ${\varOmega}\times (0,T]$,}& \\
&u^{\hbar}=0  &&\quad\mbox{on $\partial \varOmega\times [0,T]$,}&  \\
&u^{\hbar}(t=0)=u^{\hbar}_0 &&\quad\mbox{in ${\varOmega}$},&
\end{aligned}
\right.
\end{equation}
where $\varOmega\subset\mathbb{R}^d$ is a bounded domain and  $f\in L^\infty\left([0,T];L^2(\varOmega)\right)$ is a forcing term. They discretize \eqref{semiclassical} by a Crank-Nicolson finite element (CNFE) scheme and prove a posteriori error estimates of optimal order. One of the main features of the considered finite element spaces is that they are allowed to change in time.}
%
The optimal order a posteriori error bounds are derived in the $L^{\infty}_t L^2_x$ norm and the analysis includes time-dependent potentials. Furthermore the derived a posteriori estimates are valid for $L^{\infty}_t L^{\infty}_x$-type potentials as well, in contrast to the existing results in the literature which require smooth $C^1_t C^2_x$-type potentials. 

{The analysis in \cite{KK} is based on the reconstruction technique, proposed by Akrivis, Makridakis \& Nochetto, for the heat equation, cf. \cite{AMN, MN1}. In \cite{KK} the authors, following this technique, introduce a novel time-space reconstruction for the CNFE scheme, appropriate for the Schr\"odinger equation \eqref{semiclassical}. A posteriori estimates for \eqref{semiclassical} and the CNFE method were also proven by D\"ofler in \cite{Dorf}, but the estimator was not of optimal order in time.}

{ The main results of \cite{KK} can be summarized as follows: The approximations $U^n(x)$ of $u^{\hbar}(x,t_n)$, $0\le n\le N,$  are computed for a non-uniform time grid $0=:t_0<t_1<\cdots<t_N=:T$ of $[0,T]$. For each $n$, $U^n$ belongs to a finite element space (which depends on $n$) consisting of piecewise polynomials of degree $r$. By $U(x,t)$ we denote the piecewise linear interpolant between the nodal values $U^n$. 
More specifically, for $t\in[t_{n-1},t_n],$ $U(x,t):=\displaystyle\frac{t-t_{n-1}}{t_n-t_{n-1}}U^n(x)+\displaystyle\frac{t_n-t}{t_n-t_{n-1}}U^{n-1}(x)$. Then 
\begin{equation}
\label{AEE}
 \|(u^{\hbar}-U)(t)\|_{L^2(\varOmega)}\le \Ec_N^0 + \Ec_N^\text{S} + \Ec_N^\text{T}, \quad \forall t\in [0,T],
\end{equation} 
where $\Ec_N^0, \ \Ec_N^\text{S} ,\  \Ec_N^\text{T}$  are all computable quantities. More precisely, $\Ec_N^0$  accounts for the initial error, while $\Ec_N^\text{S}, \ \Ec_N^\text{T}$  are the space and time estimators respectively. These estimators are used to refine appropriately the time and space mesh sizes, thus creating an adaptive algorithm. The algorithm is said to converge up to a preset tolerance \textsc{Tol} if, after appropriate refinements, we obtain an approximate solution $U$ of $u$ with 
$$\Ec_N^0 + \Ec_N^\text{S} + \Ec_N^\text{T}<\textsc{Tol}.$$
In particular, in view of \eqref{AEE}, we will then have that
\begin{equation}
\label{eqapost}
 \|(u^{\hbar}-U)(t)\|_{L^2(\varOmega)}\le\textsc{Tol}, \quad \forall t\in[0,T].
 \end{equation}}

%
%
%
 The adaptive algorithm of \cite{KK} provides efficient error control for the solution and its observables for small values of Planck's constant $\hbar$, and in particular reduces substantially the computational cost as compared to uniform meshes.   
It is very difficult to obtain such results via standard techniques and without adaptivity, especially when non-smooth potentials are considered. In addition, it is to be emphasized that as long as the adaptive algorithm converges, we can guarantee rigorously, based on the a posteriori error analysis, that the \emph{total $L^\infty_t L^2_x$ error remains below a given tolerance, \textsc{Tol}}.
For more details, see \cite{KK}. 

\subsection{Validation of the CNFE scheme}\label{numexp}
{ We consider the one-dimensional spatial case of \eqref{semiclassical}, $\varOmega=(a,b)$, and we proceed to a series of numerical experiments which (a) validate the method and the estimators in \eqref{AEE} in terms of accuracy and (b) highlight the advantages of adaptivity. We consider the numerical solution of \eqref{semiclassical}, obtained by the CNFE scheme, with initial condition
\begin{equation*}
u^{\hbar}_0(x) = a_0(x)\text{e}^{i\frac{S_0(x)}{\hbar}},
\end{equation*}
where $a_0$ may or may not depend on $\hbar$. For the spatial discretization we use finite element spaces consisting of B-splines of degree $r,\, r\in\mathbb{N}$. The theoretical order of convergence for the CNFE scheme is $2$ in time and $r+1$ in space; thus the expected order of convergence of the estimator $\mathcal{E}_N^{\text{S}}$ is $r+1$ and of   $\mathcal{E}_N^{\text{T}}$ is $2$.

Next, our purpose  is to verify numerically the aforementioned order of convergence for the estimators, for smooth and non-smooth potentials $V$. To this end, let $\ell\in\mathbb{N}$ count the different realizations (runs) of the experiments. We consider uniform partitions in both time and space, and let $M(\ell)+1$ and $N(\ell)+1$ denote the number of nodes in space (of $[a,b]$) and in time (of $[0,T]$), respectively. Then $\varDelta x(\ell):={ \displaystyle \frac{{b-a}}{{M(\ell)}} }$ and $\varDelta t(\ell):={ \displaystyle\frac{ T}{{N(\ell)}} }$ denote the space and time discretization parameters (of the $\ell^{\text{th}}$ realization), respectively. The \emph{experimental order of convergence (EOC)} s computed for the space estimator $\mathcal{E}_N^{\Ss}$ as follows:
\begin{equation}
\label{EOCS}
\text{EOC}_{\Ss}:=\frac{\log\Big(\Ec_N^{\Ss}(\ell)/ \Ec_N^{\Ss}(\ell+1)\Big)}{\log\Big(M(\ell+1)/M(\ell)\Big)},
\end{equation}
where $\Ec_N^{\Ss}(\ell)$ and $\Ec_N^{\Ss}(\ell+1)$ denote the values of the space estimators in two consecutive implementations with mesh sizes $\varDelta x(\ell)$ and $\varDelta x(\ell+1)$, respectively. Similarly, for the time estimator $\mathcal{E}_N^{\Tt}$ the EOC is computed as 
\begin{equation}
\label{EOCT}
 \text{EOC}_{\Tt}:=\frac{\log\Big(\Ec_N^{\Tt}(\ell)/ \Ec_N^{\Tt}(\ell+1)\Big)}{\log\Big(N(\ell+1)/N(\ell)\Big)}. 
 \end{equation}

}
First, let us look at a a smooth double well potential problem with initial data
 \begin{align}
& V(x)=(x^2-0.25)^2, \quad a_0(x)= \mathrm{e}^{-\frac{25}{2}x^2}, \quad S_0(x)= -\frac{1}{5}\ln\Big(\mathrm{e}^{5(x-0.5)}+\mathrm{e}^{-5(x-0.5)}\Big), \quad \text{ with } \hbar=0.25. \label{DWP} 
\end{align}
  The computational domain  is $[a,b]\times[0,T]=[-2,2]\times[0,1]$. 
\begin{table}[htb!]
\begin{subtable}{0.5\linewidth}
\centering
\caption{Space Estimator}
\begin{tabular}{|c||c|c|}\hline
$M$ &  $ \Ec_N^{\Ss}$ & $\text{EOC}_{\Ss}$  \\ \hline\hline
$35$	   &   $7.4125$e$-01$ &--              \\
$50$  &	$1.6791$e$-01$ &$4.1633$   \\
$70$  &	$4.1761$e$-02$ &$4.1354$   \\
$100$  &	$9.7450$e$-03$ &$4.0799$  \\
$145$ &	$2.1714$e$-03$ &$4.0407$  \\ 
$200$  &	$5.9598$e$-04$ &$4.0205$ \\\hline
\end{tabular}
\end{subtable}%
\begin{subtable}{0.5\linewidth}
\centering
\caption{Time Estimator}
\begin{tabular}{|c||c|c|}\hline
$N$ &  $ \Ec_N^{\Tt}$ & $\text{EOC}_{\Tt}$  \\ \hline\hline
$80$	   &   $1.7266$e$-02$ &--              \\
$160$  &	$3.9316$e$-03$ &$2.1347$   \\
$320$  &	$9.6275$e$-04$ &$2.0299$   \\
$640$  &	$2.3943$e$-04$ &$2.0076$  \\
$1280$ &	$5.9784$e$-05$ &$2.0018$  \\ 
$2560$  &	$1.4942$e$-05$ &$2.0003$ \\\hline
\end{tabular}%
\end{subtable}%
\caption{$\text{EOC}_{\Ss}$ and  $\text{EOC}_{\Tt}$ for Double Well potential} \label{DWtbl}
\end{table}%
For the double well potential \eqref{DWP} we use cubic B-splines for the spatial discretization. The results are shown in Table \ref{DWtbl}. The predicted theoretical order of convergence is observed for both the space and time estimators. 

Now let us look at a problem with a non-smooth potential, namely
 \begin{align}
& V(x) = 10 |x|, \quad a_0(x)=\hbar^{-\frac14}\mathrm{e}^{-\frac{\pi}{2\hbar}(x-x_0)^2}, \quad S_0(x)=25\sqrt{1.5}(x-x_0), \quad \text{ with } \hbar=0.5. \label{NSP}
\end{align}
We use quartic B-spline for the space discretization and  $[a,b]\times[0,T]=[-4,4]\times[0,0.1]$ is the computational domain.
The numerical results are shown in Table \ref{NStbl} demonstrating the correct order of convergence for the estimators. It is worth noting that in this case the wavepacket passes over the non-smooth point $x=0$ during the simulation time.
\begin{table}[htb!]
\begin{subtable}{0.5\linewidth}
\centering
\caption{Space Estimator}
\begin{tabular}{|c||c|c|}\hline
$M$ &  $ \Ec_N^{\Ss}$ & $\text{EOC}_{\Ss}$  \\ \hline\hline
$ 800$   & $1.4268$e$-01$ &--              \\
$1000$  &	$4.5514$e$-02$ &$5.1204$   \\
$1200$  &	$1.8094$e$-02$ &$5.0594$   \\
$1600$  &	$4.3186$e$-03$ &$4.9799$  \\
$2000$ &	$1.4135$e$-03$ &$5.0051$  \\
$3200$  &	$1.3481$e$-04$ &$4.9998$  \\ \hline
\end{tabular}
\end{subtable}%
\begin{subtable}{0.5\linewidth}
\centering
\caption{Time Estimator}
\begin{tabular}{|c||c|c|}\hline
$\Delta t\times10^6$ &  $ \Ec_N^{\Tt}$ & $\text{EOC}_{\Tt}$  \\ \hline\hline
$10$	   &   $1.5731$e$-03$ &--              \\
$5.724$  &	$5.1538$e$-04$ &$2.0001$   \\
$3.629$  &	$2.0715$e$-04$ &$2.0001$   \\
$1.768$  &	$4.9168$e$-05$ &$1.9999$  \\
$1.012$ &	$1.6109$e$-05$ &$2.0000$  \\ 
$0.312$ &	$1.5312$e$-06$ &$1.9999$  \\  \hline
\end{tabular}
\end{subtable}
\caption{$\text{EOC}_{\Ss}$ and  $\text{EOC}_{\Tt}$ for non-smooth potential} \label{NStbl}
\end{table}

\begin{figure}[htb!]
%
{\hspace{-1.5cm}
\includegraphics[scale=0.38]{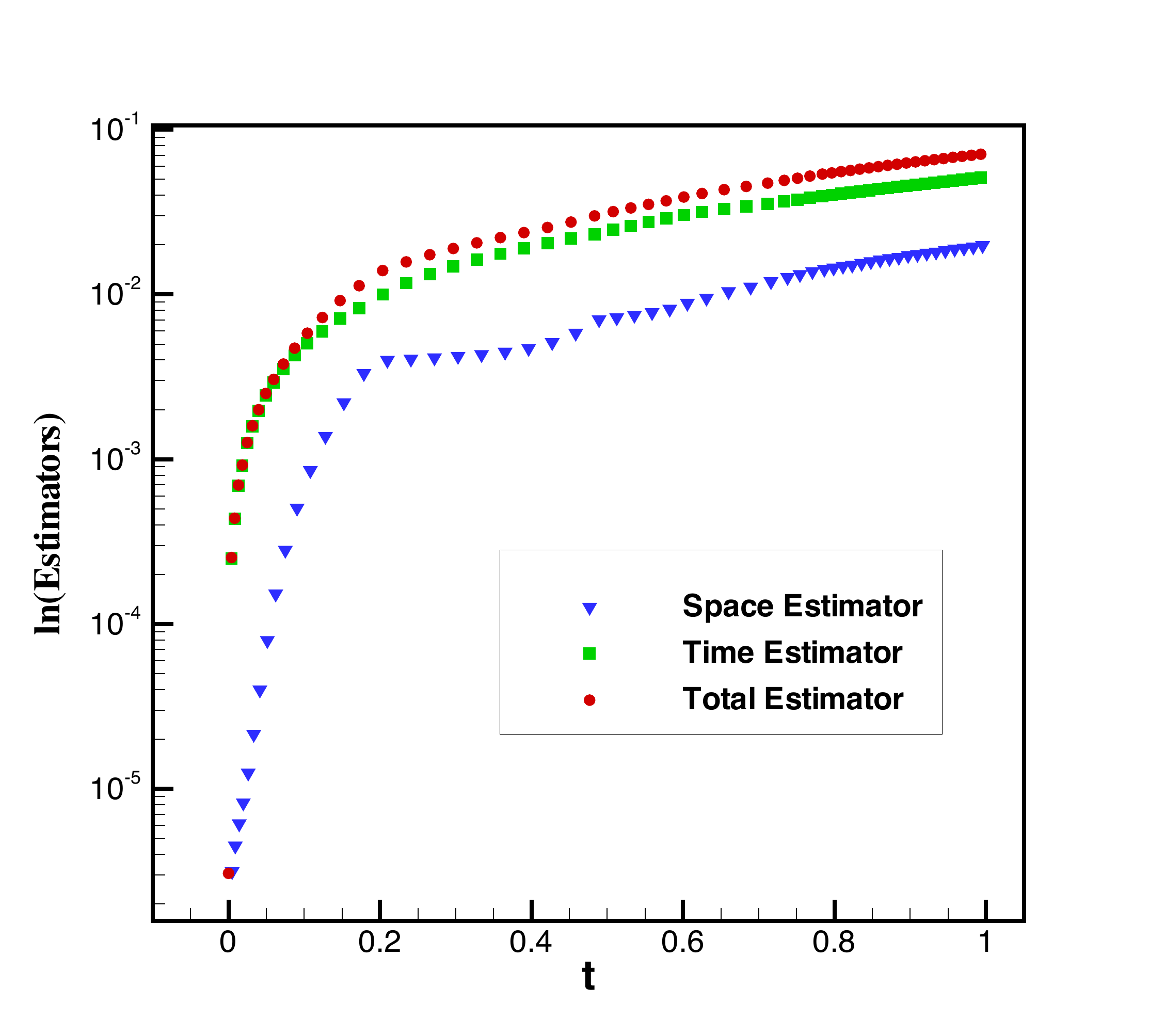}}
 {\hspace{-0.7cm}
\includegraphics[scale=0.34]{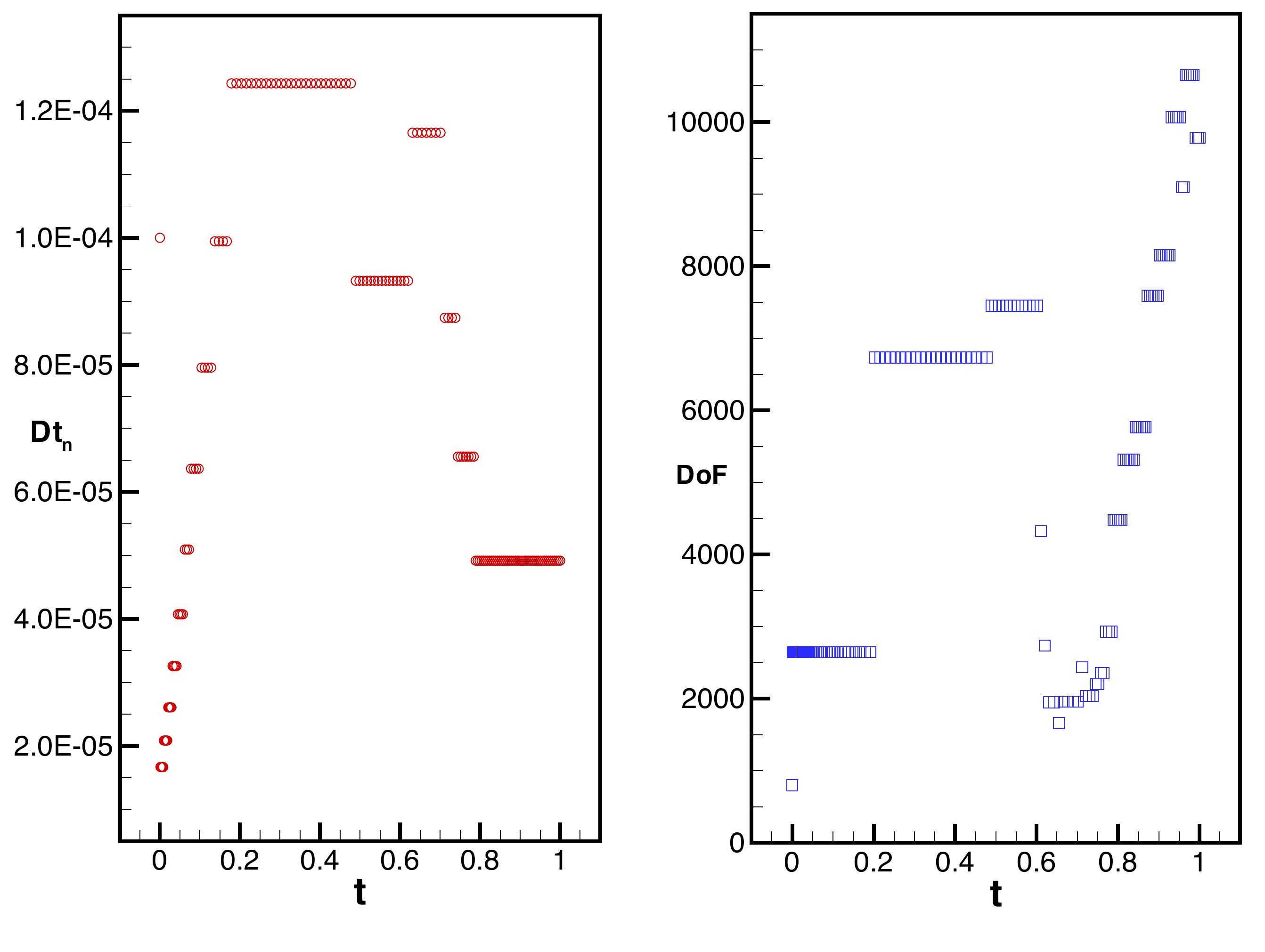} 
%
%
%
\caption{Evolution of estimators  in logarithmic scale (left) and variation of the time-steps $\varDelta t_n$ and the degrees of freedom (DoF) versus $t$ (right) during  adaptivity for  $V(x,t)=\displaystyle\frac {x^2}2\cdot\frac {1}{t+0.05}$. \label{TDPf}}
}
\end{figure}

Finally, to observe the benefits of adaptivity, we consider a   time dependent potential, namely
 \beq
 V(x,t)=\frac{ x^2 }{2(t+0.05)}, \quad a_0(x)=\mathrm{e}^{-\lambda^2(x-0.5)^2}, \quad S_0(x)=5(x^2-x) \quad \text{ with } \hbar=1. \label{TDP}
\ee
The computational domain is  $[a,b]\times[0,T]=[-1,2]\times[0,1]$ and we discretize space using cubic B-splines.
In Figure ~\ref{TDPf},  we plot the evolution of the estimators in logarithmic scale and the variation in time of the time-steps $\varDelta t_n:=t_n-t_{n-1}$ and of the degrees of freedom. This is a characteristic example where intensive adaptivity is observed, in both time and space.

\subsection{Approximation of quadratic observables } \label{sbs:compobs}
%
%

As always when discretizing problems in free space, we have to make sure the computational domain used is large enough so that (for the initial data $u_0^\hbar$ and timescale $T$ in question) the solutions of problems \eqref{semiclassical}, \eqref{eq:schroeq} are close to each other. This follows from standard localization arguments, and it is easy to check it in practice (by measuring how much mass reaches the endpoints) and poses no particular difficulty here.
Hence eq.  \eqref{eqapost} can be interpreted as an approximation between the numerical solution $U$ and the exact solution of the free space problem \eqref{eq:schroeq}. Here we
 discuss systematically how this bound can be used for the approximation of quadratic observables of the wavefunction $u^\hbar(t)$.

The quadratic observable with symbol $A_{\mathcal{W}}(x,k)$ is measured for a state $u^\hbar$ through
\[
\bac
A[u^\hbar](t)=
\langle{ \,\,  W^\hbar[u^\hbar] \,\, , \,\, A_{\mathcal{W}} \,\, }\rangle  
= \int {e^{-2\pi i K(X+Y) } A_{\mathcal{W}}\left({ \frac{X+Y}2,\hbar K}\right)   dK } \,\,  u^\hbar(X) dX \,\, \overline{u^\hbar}(Y) dY.
\ea
\]
We will be concerned with two special types of observables, namely  \textit{observables of position}, for $A_{\mathcal{W}}=A_{\mathcal{W}}(x)$
\beq \label{eqbomdebr1}
\bac
A[u^\hbar](t)=\langle{ \,\,  W^\hbar[u^\hbar] \,\, , \,\, A \,\, }\rangle  
= \int { A_{\mathcal{W}}(x)u^\hbar(x,t) \overline{u^\hbar(x,t)} dx},
\ea
\ee 
and \textit{separable observables}, $A_{\mathcal{W}}=A_1(x)A_2(k)$,
\beq
\bac \label{eqbomdebr2}
A[u^\hbar](t)=
\langle{ \,\,  W^\hbar[u^\hbar] \,\, , \,\, A \,\, }\rangle  
= \hbar^{-1} \int {e^{-2\pi i K \frac{(X+Y)}{\hbar} } A_2\left({  K}\right)   dK } \quad A_1 \left({ \frac{X+Y}2 }\right)\quad  u^\hbar(X) dX \quad \overline{u^\hbar}(Y) dY.
\ea
\ee

These observables are essentially controlled by the $L^\infty_tL^2_x$ norm of the wavefunction; this is made more precise in the following

\begin{lemma}[Approximation of observables] If $\|u^\hbar-U\|_{L^2} \leqslant \textsc{Tol}$ as in \eqref{eqapost}, then for every observable of position 
\beq \label{eqbomdebr123}
|A[u^\hbar](t) -  A[U](t)| \leqslant  \textsc{Tol}(  \|U\|_{L^2} + \|u^\hbar\|_{L^2}  )\| A_{\mathcal{W}}\|_{L^\infty} = O( \textsc{Tol}),
\ee
while for any separable observable
\beq \label{eqbomdebr12}
|A[u^\hbar](t) -  A[U](t)| \leqslant \hbar^{-\frac{1}2}   \textsc{Tol}(  \|U\|_{L^2} + \|u^\hbar\|_{L^2}  )\| A_{\mathcal{W}}\|_{L^2} = O( \hbar^{-\frac{1}2}  \textsc{Tol}).
\ee
\end{lemma}

The proof follows by inspection of equations \eqref{eqbomdebr1}, \eqref{eqbomdebr2}. 

\vspace{0.25 cm}
\begin{remark} \upshape The estimate \eqref{eqbomdebr12} is far from sharp; in fact for regular, localized observables the $\hbar^{-\frac{1}2}$ is very pessimistic. Still, carrying out rigorously a sharper microlocal estimate for a non-smooth problem is outside the scope of this work. We note that, even in an imperfect way, it is seen rigorously that the $L^2$ approximation of the wavefunction does indeed control the observables.
\end{remark}


\subsection{Particles for the Liouville equation} \label{sbs:particles}

To approximate numerically the solution of \eqref{eqsert45},  we use a particle method; decompose the initial condition
\[
\rho^\hbar_0 \approx \sum\limits_{j=1}^N M_j\delta(x-X_j,k-K_j);
\]
then the center of each particle moves along its respective trajectory, in accordance to \eqref{eq:char}. (See also the caption of Figure \ref{Fig101} for an explicit form of the trajectories.) Thus 
\[
\rho^\hbar(t) \approx P^\hbar(t)= \sum\limits_{j=1}^N M_j\delta(x-X_j(t),k-K_j(t)).
\]
The advantage in this case is that we know explicitly the trajectories, and therefore 
\[
\langle \rho^\hbar(t) - P^\hbar(t),\phi \rangle=
\langle \rho^\hbar_0 - P^\hbar(0),\phi \rangle.
\]
This makes it easy to generate approximations of observables of $\rho^\hbar(t)$, i.e.  
$$\int{\rho^\hbar(x,k,t)A_{\mathcal{W}}(x,k)dxdk} \approx \sum\limits_j M_jA_{\mathcal{W}}(X_j(t),K_j(t))$$ 
with predetermined accuracy.

\section{Numerical results} \label{secresopiijl}

In this section we present a series of one-dimensional numerical experiments, investigating whether an appropriate version of Theorem \ref{thrm:main} can be seen to hold for $a=0$, i.e. if
\beq \label{eq:chtesrt}
\langle W^\hbar[u^\hbar(t)] - \rho^\hbar(t),\phi\rangle = o(1)
\ee
holds over saddle points of the form $V(x)=-C|x|$.
More specifically, we work with the non-smooth potential of type \eqref{eq:setconic}, namely we take 
\begin{equation}
\label{VAbsx}
V(x) = 1 + (1+\tanh(4(x+2.5)))(1+\tanh(-4(x-2.5)))\frac{(-|x|+4)}{8}.
\end{equation}
\begin{figure}
\includegraphics[scale=0.25]{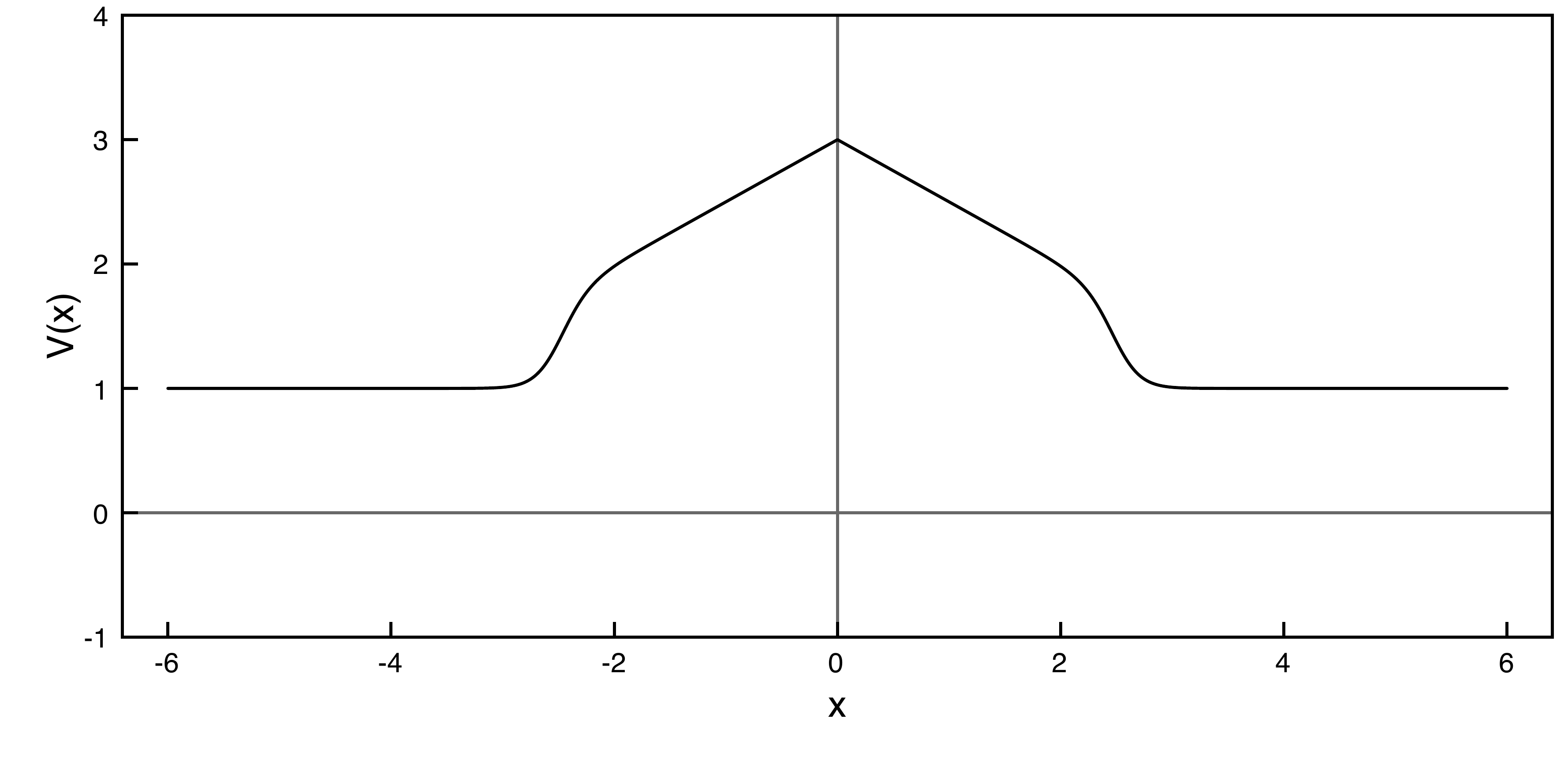}
\caption{The non-smooth potential V of eq. \eqref{VAbsx}.}\label{MAbsX}
\end{figure}
This potential 
incorporates the non-smoothness at $x=0$ with a smooth transition to a constant value away from it. Note that in a neighborhood of $x=0$, $V$ is exponentially close to $-\frac{|x|}{2}+3$, see Figure \ref{MAbsX}.

We compute 
the numerical solution of problem \eqref{semiclassical}, which is known to approximate well problem \eqref{eq:schroeq} as long as the effective support of the solution doesn't reach the boundary of the computational domain. It will be referred to as the ``exact wavefunction'', and denoted by $u^\hbar$ in the sequel. (``Exact'' in the sense that the full, quantum dynamics are used.)
The wavefunction $u^\hbar$ is computed with a prescribed error tolerance of
  $\textsc{Tol}\approx 0.01$ (more specifically $\textsc{Tol}\in [0.005,0.02]$).

  We also compute the numerical solution to \eqref{eqsert45} with initial data $\rho^\hbar_0 = \widetilde{W}^\hbar_0$,  i.e. a smoothing of $W^\hbar_0$. We write for future reference
\beq
\label{eqsert4567}
\partial_t \rho^\hbar(t) + 2\pi k\cdot \partial_x \rho^\hbar(t) - \frac{1}{2\pi} \partial_x V \cdot \partial_k \rho^\hbar(t)=0, \,\,\,\,\,\,\,\,\,\,\,\,\,\,
\rho^\hbar(t=0)=\widetilde{W}^\hbar_0.
\ee  
  This will be referred to as the ``classical SWT'', and denoted by $\rho^\hbar(t)$ in the sequel. As was discussed in section \ref{sbs:particles}, (almost all) the trajectories can be computed explicitly.
The initial data $u_0^\hbar$ are chosen so that there is full interaction with the singularity of the flow, in the sense of definition \ref{definter}.

So we have two reliable computations; one for the full quantum dynamics of the problem, and one for a semiclassical model inspired by Theorem \ref{thrm:main}. We proceed to measure a number of observables against $u^\hbar(t)$, $\rho^\hbar(t)$ before, during, and after interaction with the saddle point. This is equivalent to checking whether eq. \eqref{eq:chtesrt}
holds for a number of test functions (the Weyl symbols of the observables).

In the process of setting up the numerical experiments and interpreting the results, a clear dichotomy arises between problems with and without \textit{interference}. A brief, precise definition can be given as follows:

\begin{definition} \label{defintrfffff} Given $u_0^\hbar$, $V$ consider the problem 
\beq \label{eq:ioghh}
\partial_t f^\hbar(t) + 2\pi k\cdot \partial_x f^\hbar(t) - \frac{1}{2\pi} \partial_x V \cdot \partial_k f^\hbar(t)=0, \,\,\,\,\,\,\,\,\,\,\,\,\,\,
f^\hbar(t=0)=W^\hbar[u^\hbar_0].
\ee
We say that interference is observed on a point $(x_*,k_*)$ of phase-space if (non-negligible for $\hbar=o(1)$) amounts of mass of $f$ arrive to $(x_*,k_*)$ at the same time from different directions.
\end{definition}

Clearly, interference is only possible where two trajectories intersect in finite time. For our potential $V$ as in eq. \eqref{VAbsx}, this is only the point $(0,0)$. If one wavepacket approaches $(0,0)$ from one side, there is no interference going on. Interference would be taking place if two wavepackets arrive on $(0,0)$, one from the right and one from the left, at the same time.

\subsection{Non-interference problems} \label{sbs:noninr}

For values of $\hbar$ ranging from $5\cdot 10^{-1}$ to $5\cdot10^{-3}$, we simulate the evolution in time of wavepackets of the form
\begin{equation}
\label{u0nonint}
\bac
 u_0(x) = a_{0}(x) \mathrm{e}^{i m \frac{S_{0}(x)}{\hbar}}, \quad
 a_{0}(x)= \hbar^{-\frac14}\mathrm{e}^{-\frac{\pi}{2} (\frac{x-x_0}{\sqrt{\hbar}})^2 }, \quad
 S_{0}(x) = \sqrt{|x_0|}(x-x_0),
\ea
\end{equation}
for $x_0=-1.5$, $m \in [0.8165,1.4289]$.

When $m=1$, the SWT $\widetilde{W}^\hbar[u^\hbar]$ of this problem is centered on $(-1.5,\sqrt{\frac{3}2})$; this point reaches zero in $t=\sqrt{6}$, and roughly half the mass of the quantum particle -- should \eqref{eq:chtesrt} hold -- is expected to pass to $\{ x>0 \}$, while the other half should reach close to $x=0$ and then be  reflected back to $\{ x<0 \}$. By perturbing the value of $m$ in the initial data, the amount of mass expected to cross over to $\{ x>0 \}$ changes (from no mass crossing over, to all the mass crossing over in the extreme cases). In all of these case studies the interaction with the singularity starts around $t=1.3$ and is over around $t=2.45$. Thus e.g. before the interaction with $x=0$ the classical and quantum solutions should agree very well, which provides one more opportunity to validate and check our computations.
We  look at $\rho^\hbar(t)$ and $W^\hbar[u^\hbar(t)]$ in phase space, and we measure the observables with symbols
\beq \label{eqobsset}
A_{\alpha,\beta,j}(x,k) =  x^\alpha \,\, k^\beta  \,\, \chi_{[0,4]}((-1)^j x) \chi_{[-1,1]}(k), 
\quad \quad \mbox{ for } \alpha,\beta \in \mathbb{N}_0, \,\,\, \alpha + \beta \leqslant 2, \,\,\,\, j \in \{1,2\}.
\ee
The precise measurement of these observables corresponds to
\begin{equation*}
\begin{aligned}
&\langle{ \,\,  W^\hbar \,\, , \,\, x^\alpha k^\beta   \,\, \chi_{[0,4]}((-1)^j x) \chi_{[-1,1]}(k) \,\, }\rangle =\\ { } \\
&=\int {e^{-2\pi i K(X+Y)} \left[{ \chi_{[0,4]}((-1)^j\frac{X+Y}{2}) \chi_{[-1,1]}(\hbar K) (\frac{X+Y}2)^\alpha (\hbar K)^\beta }\right]\,  dK } \quad  u^\hbar(X)\, dX \quad \overline{u^\hbar}(Y)\, dY.
\end{aligned}
\end{equation*}
For $\beta=0$ these are observables of position only, so by the estimate \eqref{eqbomdebr123} we have a very good approximation. 
For $\beta>0$ we do not attempt to saturate the estimate \eqref{eqbomdebr12}, since it is quite clear from the numerical results that it is not necessary. Our findings are fully consistent for both types of observables, as we will see below.

\begin{figure}[htb!]
\includegraphics[width=148mm,height=74mm]{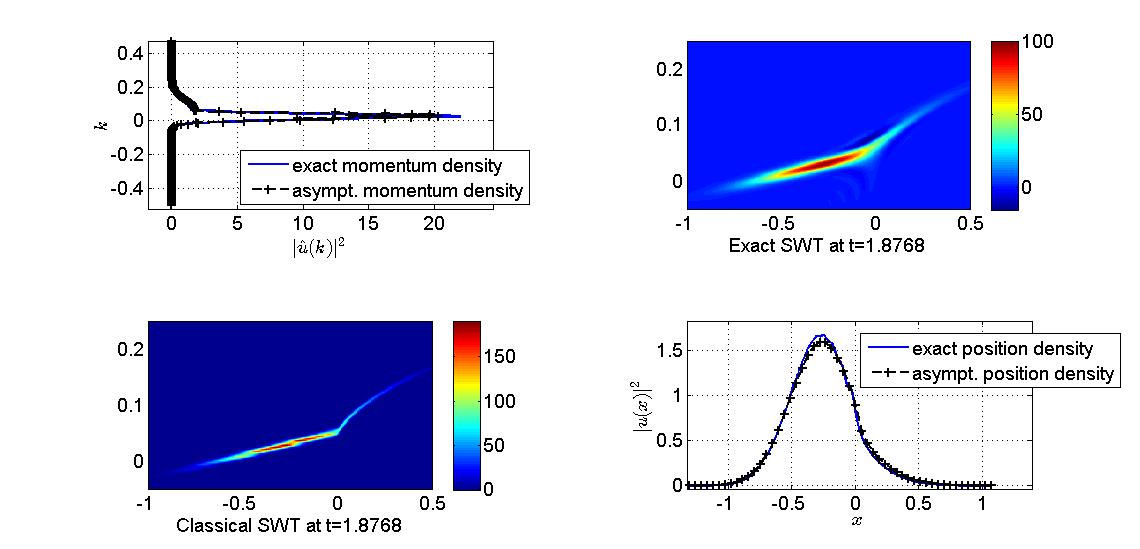} 
\caption{ Numerical result for $\hbar=10^{-2}$, $m=0.9186$.
Top right: Exact SWT. Top left: momentum density ($dx$ integral of SWT). Bottom right: position density ($dk$ integral of the SWT). Bottom left: $\rho^\hbar(t)$. (Note that in the SWT plots the wavenumber is scaled with $\frac{1}{2\pi}$.)
}
\label{figA1}
\end{figure}

\begin{figure}[htb!]
\begin{tabular}{l r}
\includegraphics[width=68mm,height=60mm]{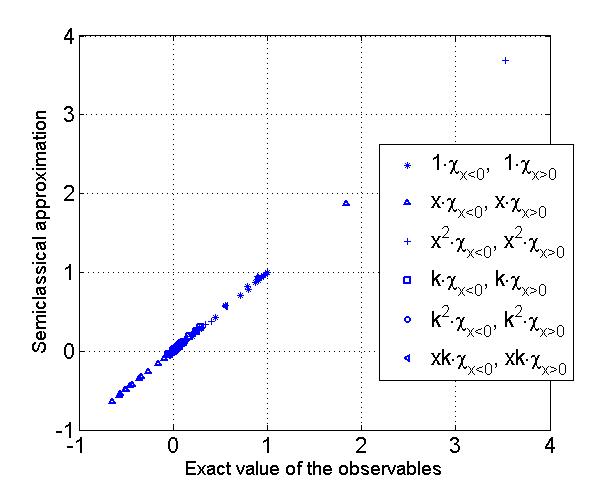} &
\includegraphics[width=68mm,height=60mm]{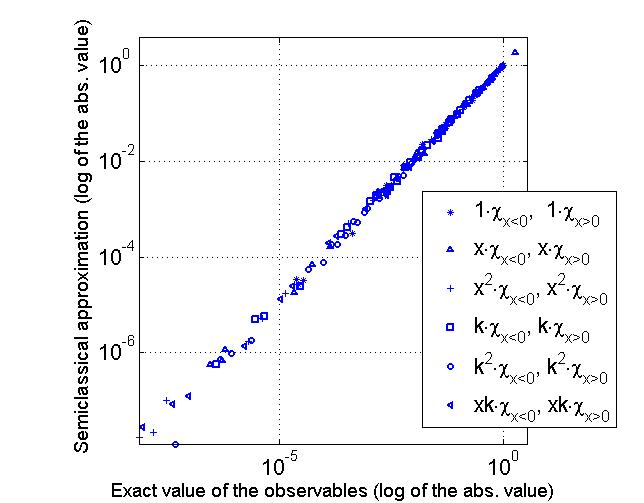} 
\end{tabular}
\caption{ Observable measurements, for the observables in \eqref{eqobsset} and $\hbar=10^{-2}$, $m\in\{ 0.8165, 0.8777,0.9186,1.0206,1.4289\}$, at times $t\in [1.1788, 2.3577]$. The $x-$coordinate of each point is the measurement on the numerical solution $U(t)$, $\mathcal{A}_{quant}=\langle A(x,k), W^\hbar [U(t)]\rangle$, and the $y-$coordinate is the corresponding classical measurement $\mathcal{A}_{cl}=\langle A(x,k), P^\hbar (t)\rangle$. Note that the times used here roughly span the interaction time, in which any discrepancy between the classical and quantum dynamics could occur. More specifically the interaction starts around $t=1$ and is over by $t=2.4$ for all the problems. 
}
\label{figA2}
\end{figure}

The agreement we find between the quantum dynamics and the proposed semiclassical asymptotics  is striking already from relatively large values of $\hbar$. This is not entirely unexpected, as away from $x=0$ the Liouville equation \eqref{eq:wmeq} is in fact identical with the full quantum dynamics \eqref{eq1hsdvaol}. The finding is that ``nothing non-classical happens'' on $x=0$ either, as can be clearly seen in Figures \ref{figA1}, \ref{figA2} and \ref{figA3}.

\begin{figure}[htb!]
\begin{tabular}{l r}
\includegraphics[width=68mm,height=60mm]{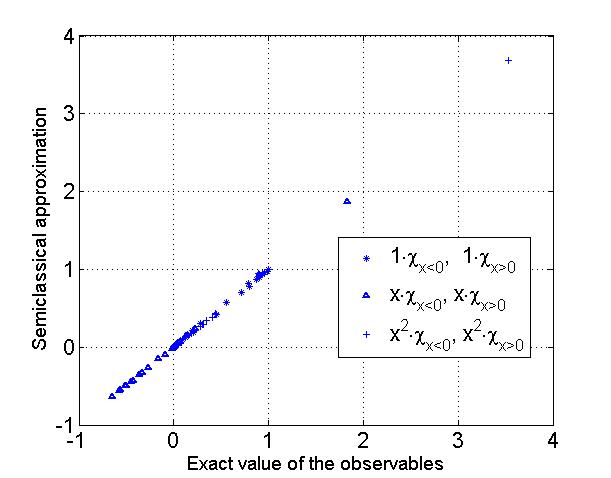} &
\includegraphics[width=68mm,height=60mm]{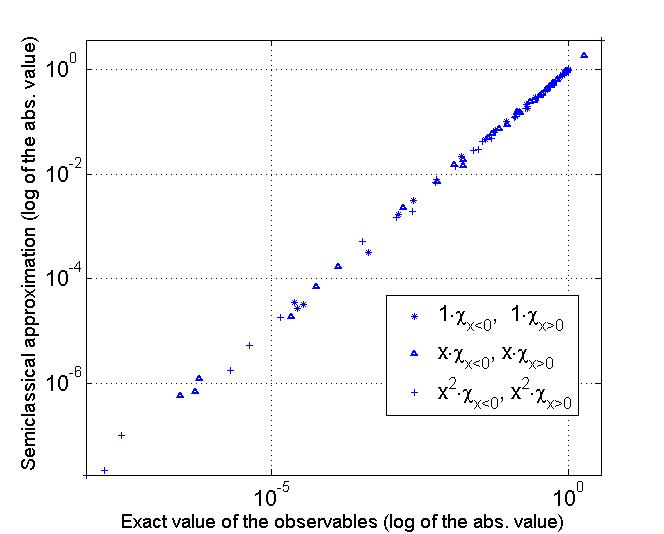} 
\end{tabular}
\caption{ Measurements of observables of position only. The qualitative behavior is consistent with the larger dataset. These benefit from better accuracy, by virtue of eq. \eqref{eqbomdebr123}. It is clear that qualitatively the picture doesn't change when we include observables depending on momentum as well (i.e. as in Figure \ref{figA2}). This is not surprising, since the simple estimate of eq. \eqref{eqbomdebr12} is known to be pessimistic.
}
\label{figA2323}
\end{figure}

\begin{figure}[htb!]
\includegraphics[width=64mm,height=50mm]{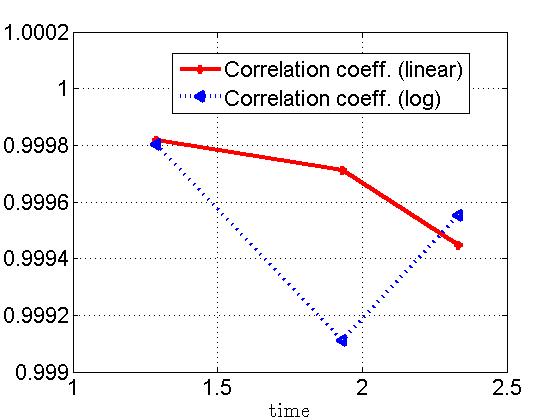} 
\caption{ Given two vectors of measurements $\vec x$ and $\vec y$, where we expect $x_i \approx y_i$, a standard way to measure how well they line up is through the correlation coefficient $\rho_{x,y}=\frac{\langle \vec x, \vec y \rangle}{\|\vec x\| \, \|\vec y\|}$, with $\rho_{x,y}=1$ if the two vectors are exactly aligned.
Here we plot the correlation coefficients for groupings of measurements that correspond to early, high and late interaction stages. We use both the linear and log scaling (as in Figure \ref{figA2}). The agreement is striking (and most probably numerical errors are comparable to any quantum-classical discrepancies). 
}
\label{figA3}
\end{figure}

Qualitatively this behavior also appeared in investigating problems with different envelopes and other values of $\hbar$. This creates a compelling sense that in non-interference problems, eq. \eqref{eq:chtesrt} is valid. In Appendix \ref{AppCwkb} an even more singular example can be seen to be correctly captured by the regularized semiclassical asymptotics.

\subsection{Collision of two wave packets}  \label{sec:interf11}

Since this is an one-dimensional problem, the only way to have interference is by one wavepacket arriving to $x=0$ from the left, and one from the right at the same time. So we consider the collision of two wave packets, symmetrically located around $x=0$, traveling with same velocities and opposite directions. The two wavepackets have a phase difference of an angle $2\pi\theta, \ 0\le\theta\le 1$. They meet over the corner point of the potential, they interact and continue to travel in opposite directions until they are completely separated.  

 The initial datum, $\rho^\hbar_0=\widetilde{W}^\hbar_0[u_0^\hbar]$ is symmetric around $(0,0)$, up to exponentially small terms, \textit{for all $\theta$}. To see that, we compute
 \begin{equation} \label{eqwtbdtwwvf}
\begin{aligned}
W^\hbar[u_0^\hbar](x,k)  =& \frac{2}{\hbar} e^{ -\pi \frac{(x-x_0)^2}{\hbar} - 
 4\pi \frac{(k-\frac{\sqrt{|x_0|} }{2\pi}  )^2}{\hbar} }+ 
 \frac{2}{\hbar} e^{ -\pi \frac{(x+x_0)^2}{\hbar} - 
 4\pi \frac{(k+\frac{\sqrt{|x_0|} }{2\pi}  )^2}{\hbar} }\\
& +2 Re \left({
\frac{2}{\hbar}   e^{ -\frac{\pi}{\hbar}x^2 - \frac{4\pi}{\hbar}k^2  } e^{- 2\pi i \theta - \frac{2i\sqrt{|x_0|}}{\hbar}x + 4\pi i  x_0 k}
 }\right);
 \end{aligned}
 \end{equation}
 For all practical purposes the third term is suppressed by the smoothing (since it is highly oscillatory), and with it all trace of $\theta$ in $\widetilde{W}^\hbar[u^\hbar_0]$.
  The flow is also symmetric around $(0,0)$, hence one quickly observes that $\rho^\hbar(t)$ in this problem predicts a distribution of mass symmetric around zero.

A crucial observable we study closely for this problem, is
 the amount of {mass} located to each side of $x=0$ after the crossing is completed, $\int\limits_{\pm x>0} |u^\hbar_0(x,t_*)|dx$ for $t_*$ sufficiently large. An eventual mass imbalance means that there are interactions going on not included in the classical dynamics of \eqref{eqsert45}.

\begin{figure}[htbp]
\begin{center}
\includegraphics[width=8cm, height=8cm]{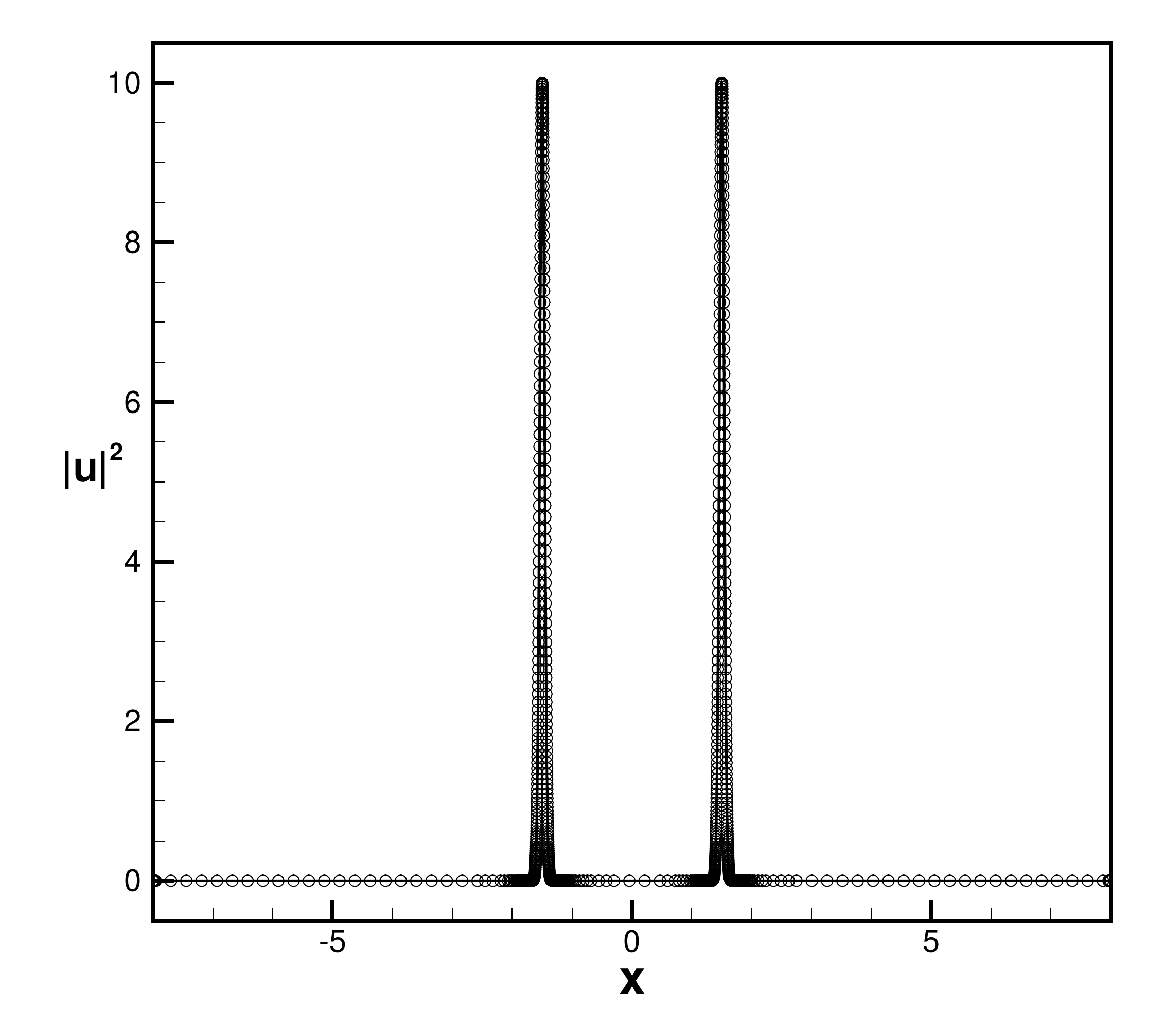}
\includegraphics[width=8cm, height=8cm]{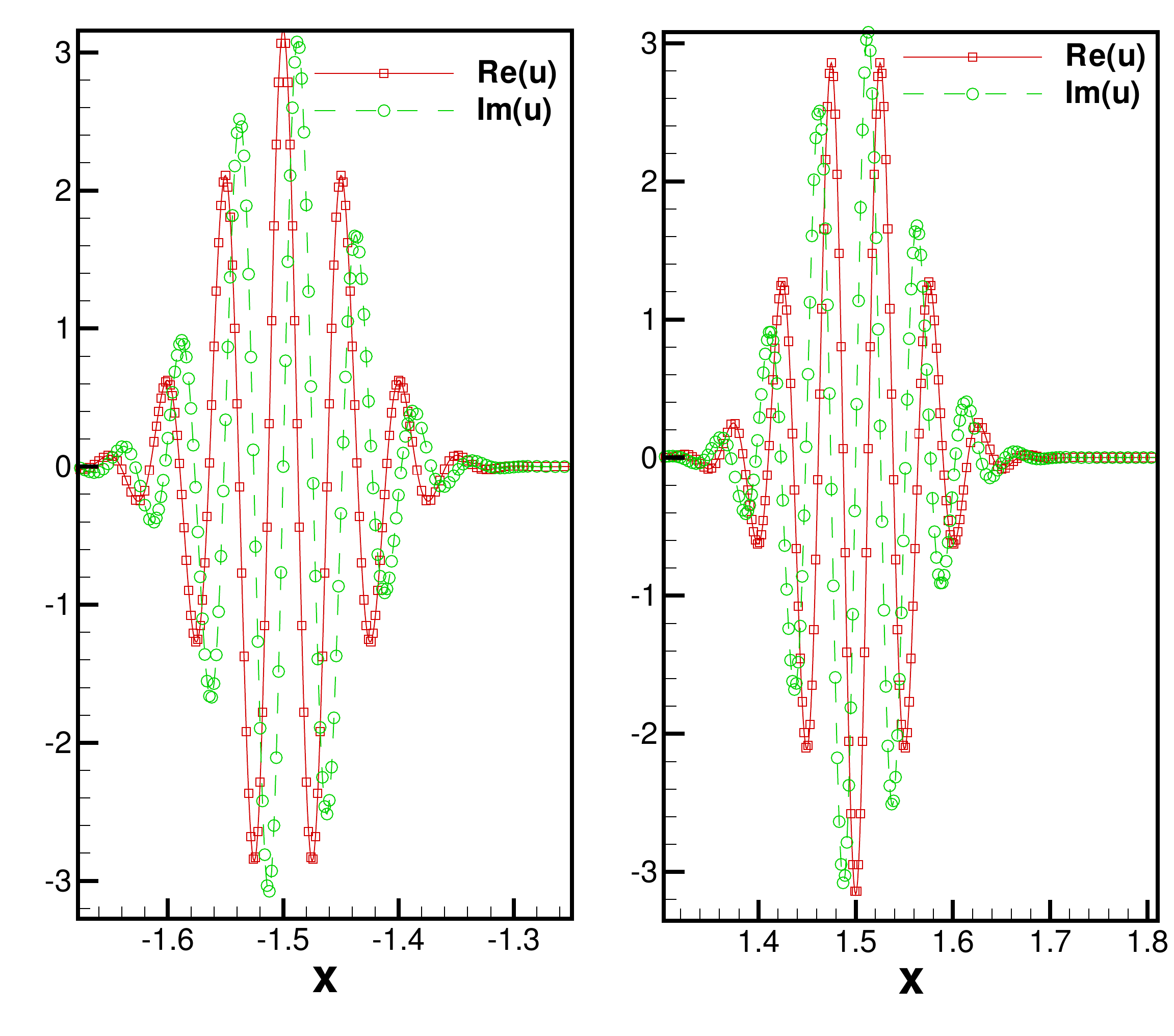} 
\caption{ Visualization of $u_0$ (as defined in \eqref{u0coll}, for $\hbar=1e-2$).
Left: Position density. Right: Real and imaginary parts for the two components of $u_0$ (right) }
\label{CollU0}
\end{center}
\end{figure}

\begin{figure}[htb!]
\includegraphics[width=148mm,height=74mm]{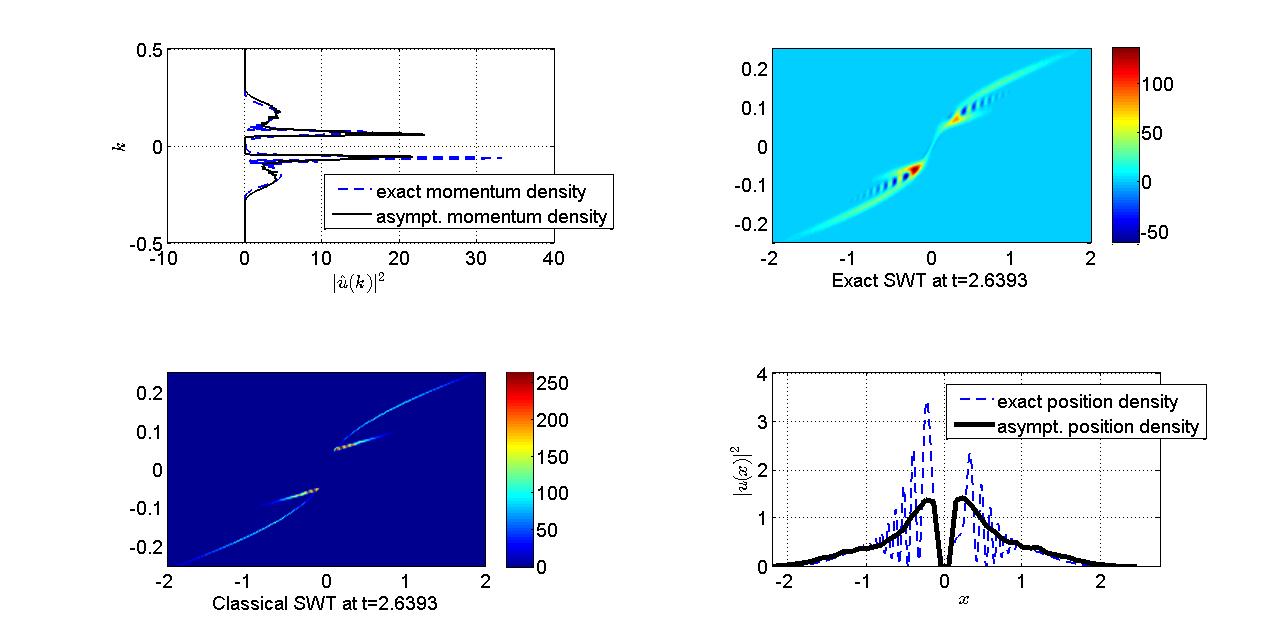} 
\caption{ Numerical result for $\hbar=5\cdot 10^{-3}$, $\theta=\frac{18}{24}$.
Top right: Exact SWT. Top left: momentum density ($dx$ integral of SWT). Bottom right: position density ($dk$ integral of the SWT). Bottom left: $\rho^\hbar(t)$. (Note that in the SWT plots the wavenumber is scaled with $\frac{1}{2\pi}$.)
}
\label{figA0001}
\end{figure}

\begin{figure}[htbp]
\begin{center}
\includegraphics[scale=0.27]{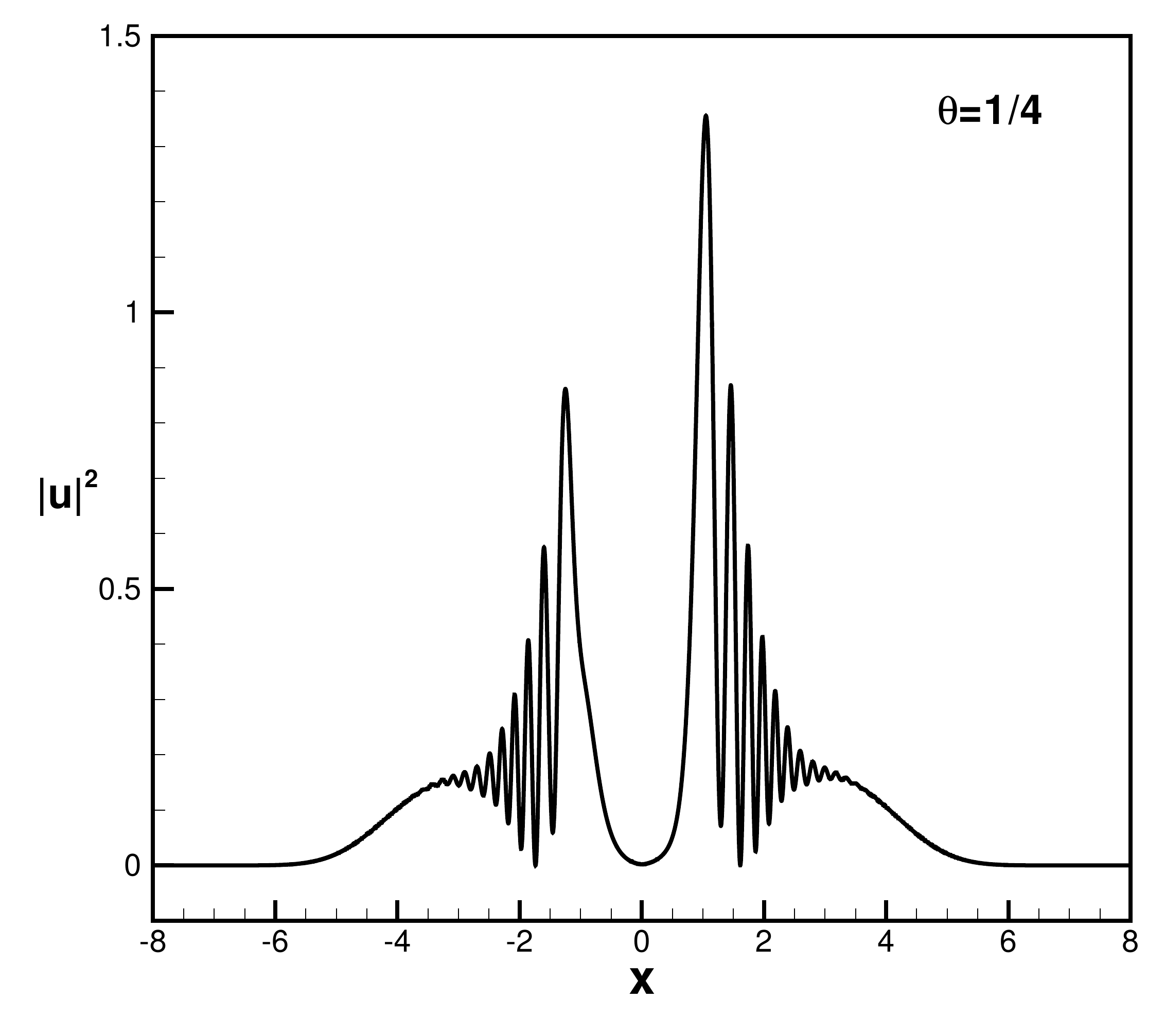}
\includegraphics[scale=0.27]{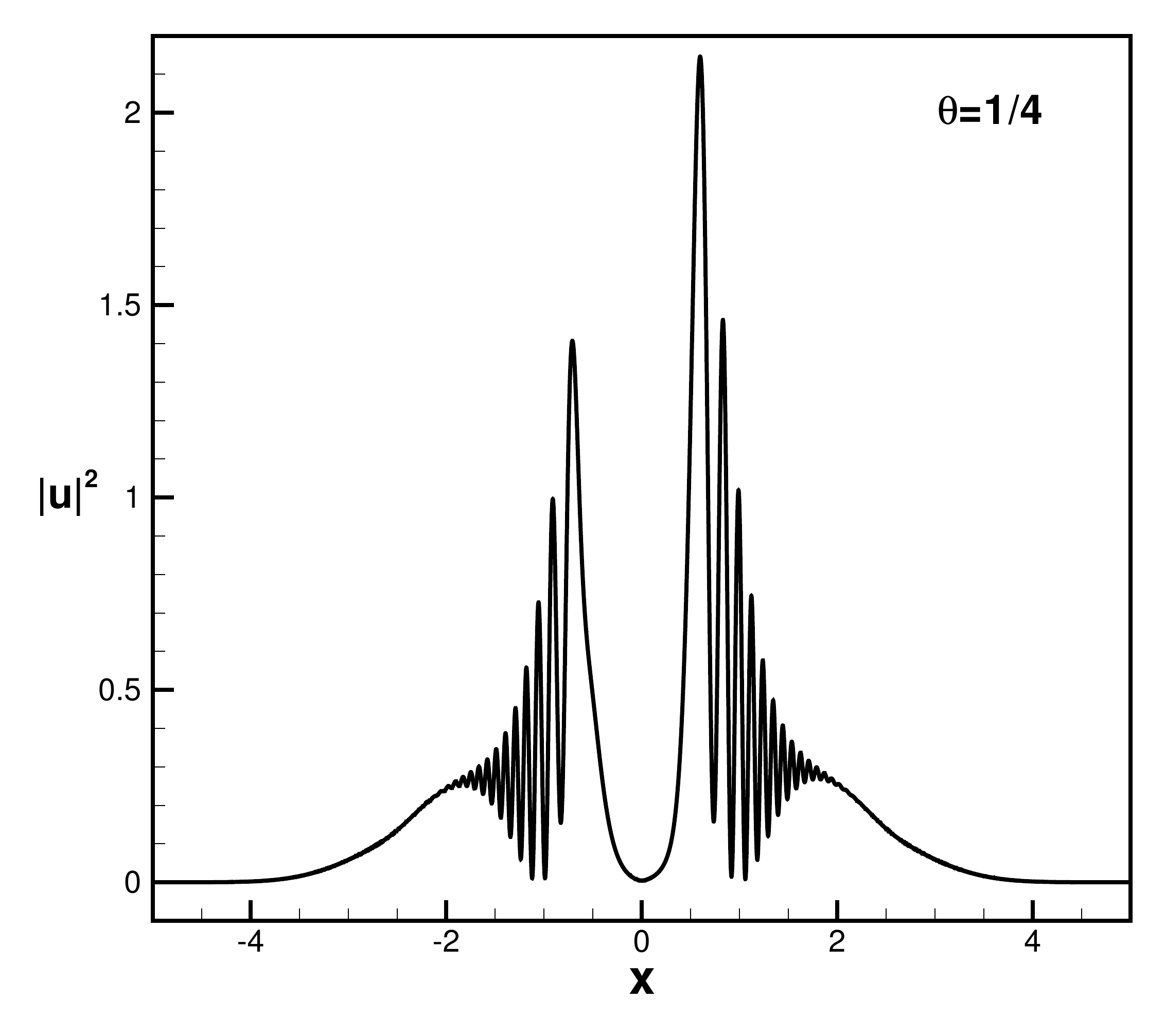}
\includegraphics[scale=0.27]{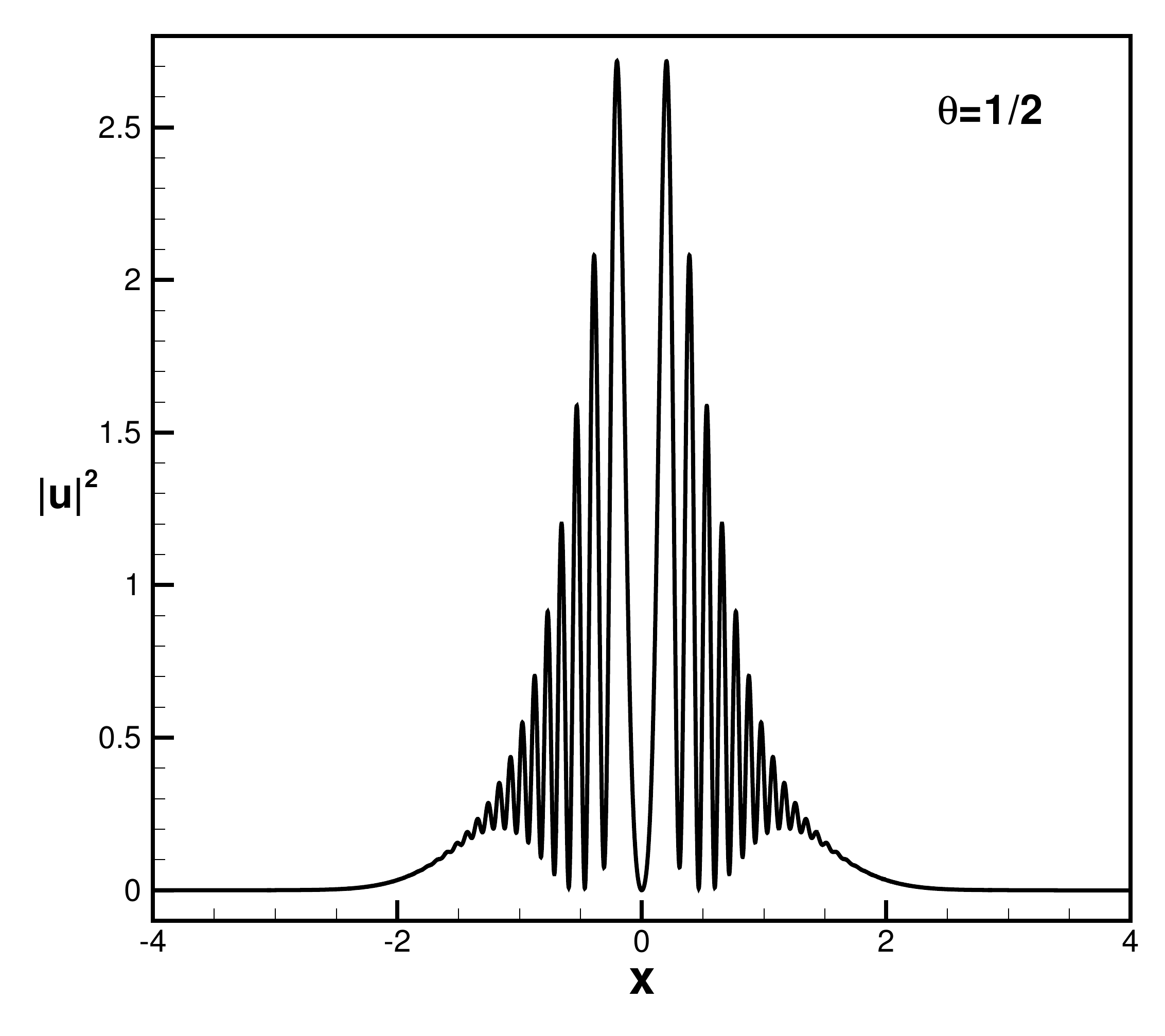}
\includegraphics[scale=0.27]{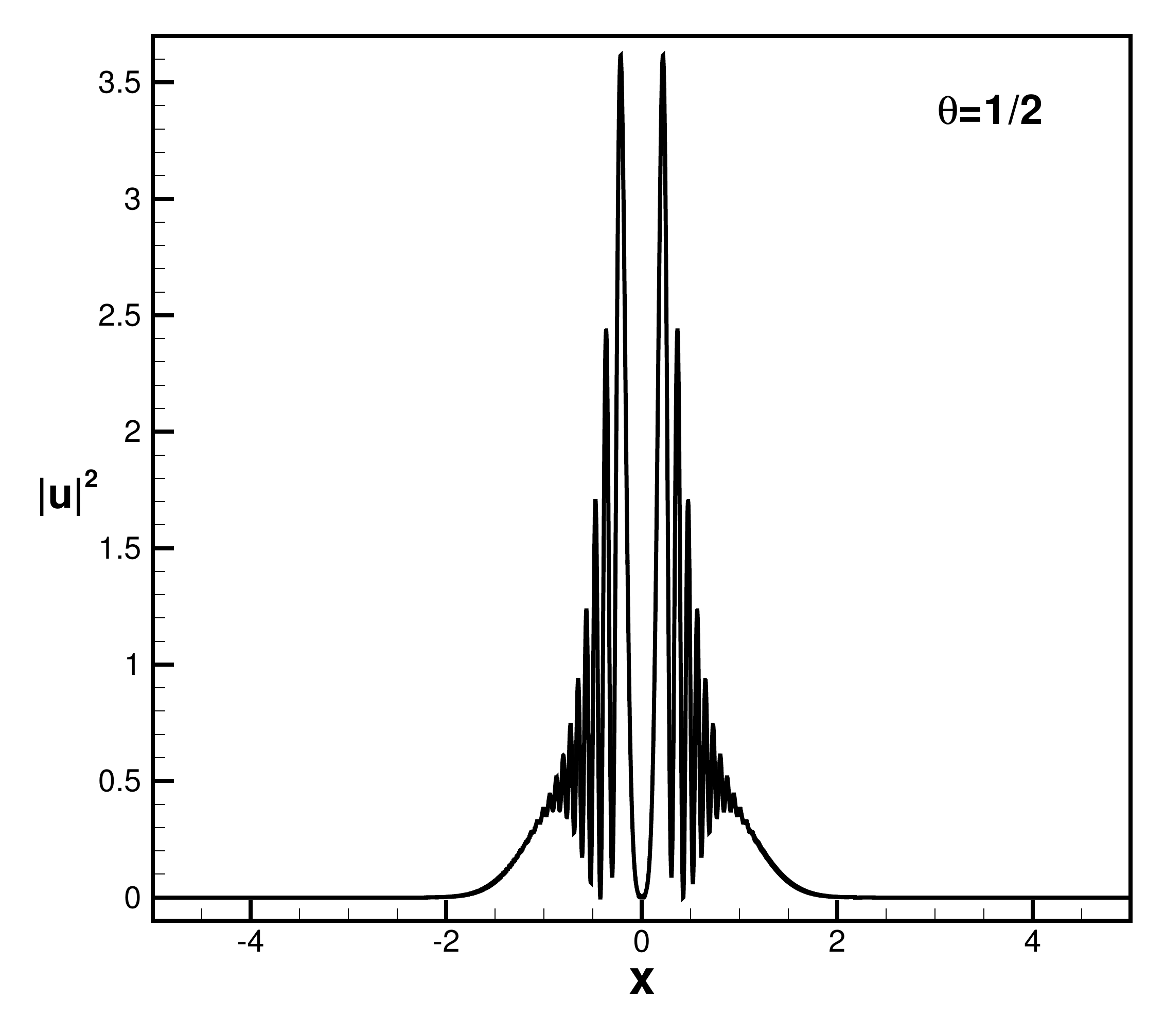}
\includegraphics[scale=0.27]{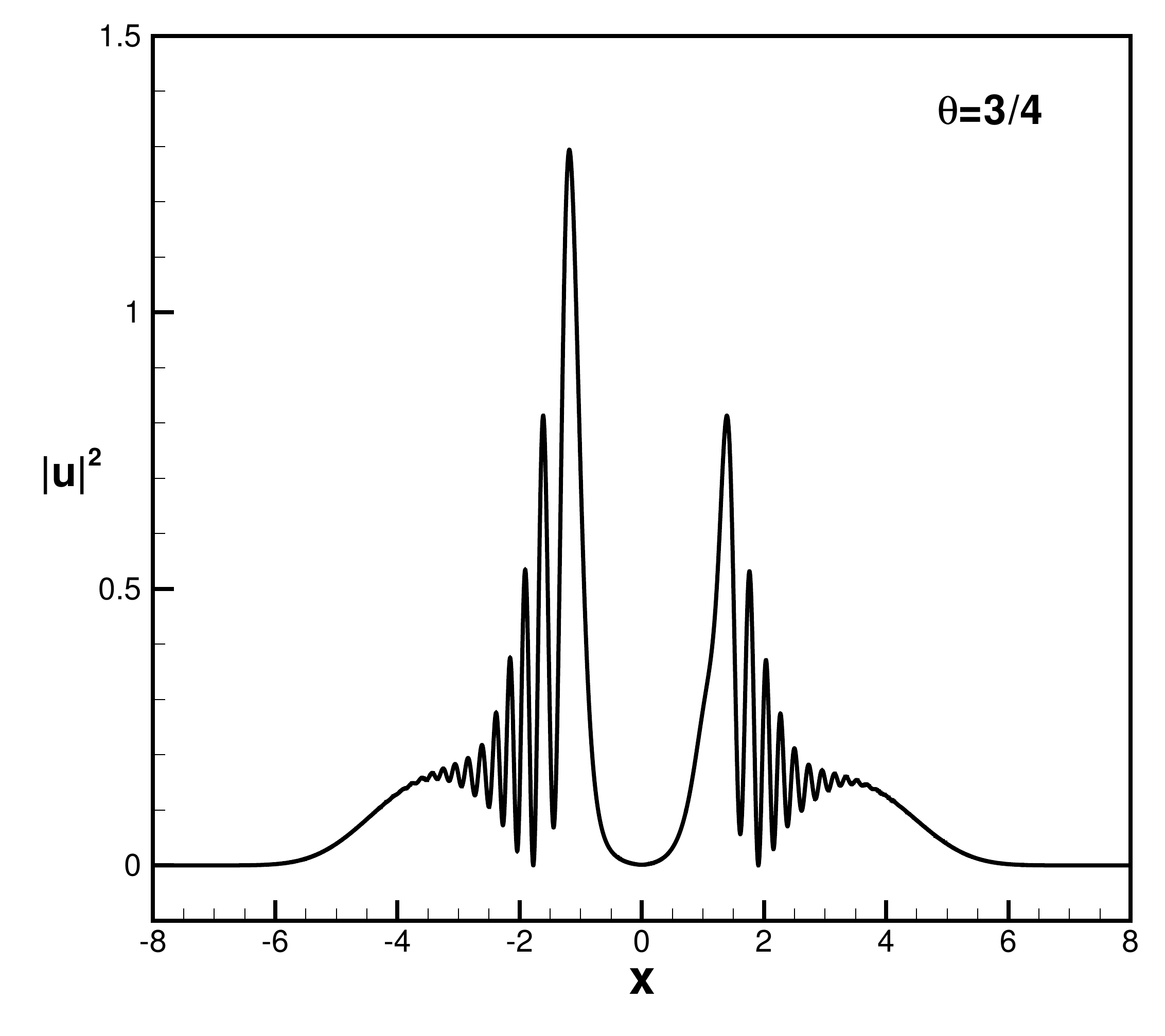}
\includegraphics[scale=0.27]{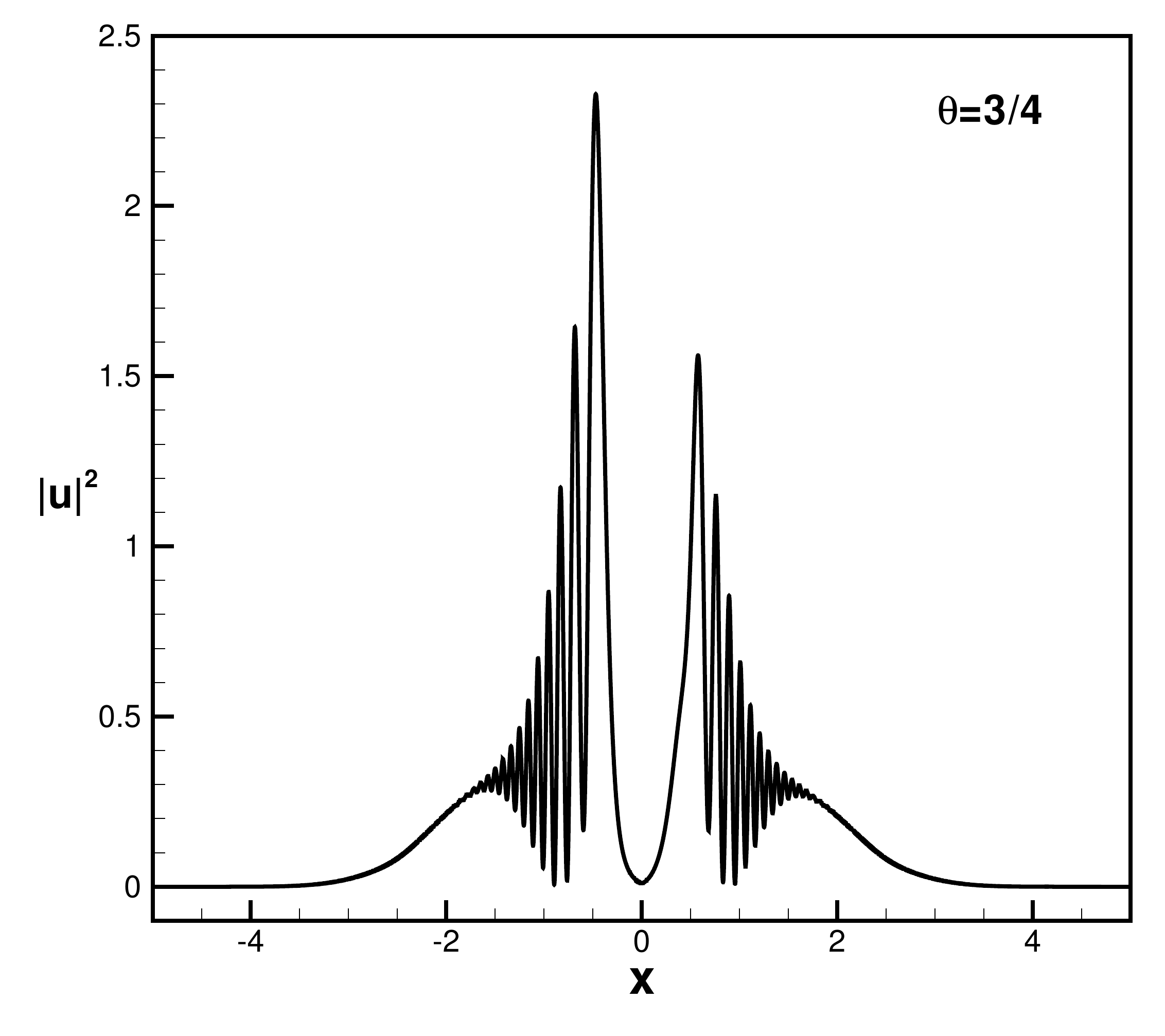}
\caption{Mass distribution for the problem with initial data \eqref{u0coll} after the interaction for various values of $\theta$: $\hbar=10^{-2}$(left) and $\hbar=5\cdot 10^{-3}$(right)}
\label{Coll2}
\end{center}
\end{figure}

The computational domain is taken sufficiently large to avoid possible interactions with the boundary and we discretize in space using quintic B-splines. The initial condition is of the form 
\begin{equation}
\begin{aligned}
\label{u0coll}
& u_0(x) = a_{0,1}(x) \mathrm{e}^{i \frac{S_{0,1}(x)}{\hbar}}  + a_{0,2}(x) \mathrm{e}^{i \frac{S_{0,2}(x)}{\hbar}}  \mathrm{e}^{i 2\pi \theta} ,\quad 0\le\theta\le 1, \\
& a_{0,1}(x)= \hbar^{-\frac14}\mathrm{e}^{-\frac{\pi}{2} (\frac{x-x_0}{\sqrt{\hbar}})^2 },  \quad a_{0,2}=\hbar^{-\frac14}\mathrm{e}^{-\frac{\pi}{2} (\frac{x+x_0}{\sqrt{\hbar}})^2 }, \\
& S_{0,1}(x) = \sqrt{|x_0|}(x-x_0), \quad S_{0,2}(x) =- \sqrt{|x_0|}(x+x_0). 
\end{aligned}
\end{equation}
  In Figure \ref{CollU0} the graphs of $u_0(x)$, $|u_0(x)|^2$ are shown for $\hbar=10^{-2}$. The wave packets are located initially at $x_0=-\frac32$ and $-x_0=\frac32$ respectively. The initial step of the adaptive algorithm resolves correctly the profile of $u_0$  producing an initial mesh, depicted also in Figure \ref{CollU0},  of around $1000$ points with an initial error bound approximately $ 10^{-9}$.  In what follows the total error \eqref{AEE} is kept under $10^{-2}$. 
In Figure  \ref{Coll2} the mass distribution is shown for three values of the parameter $\theta=\frac14, \frac12, \frac34$ and two values of  Planck's constant $10^{-2}, \ 5\cdot10^{-3}$. The snapshots correspond to a time where the two wave packets have interacted with each other over the corner of the potential and continue to move away from it.  
In Figure \ref{figA0001} we see the numerical approximations for $\widetilde{W}^\hbar[u^\hbar]$ and $\rho^\hbar$ at a time after the interaction. The classical approximation is completely symmetric, while the quantum result is not. This is a case of a ``microscopic'' (i.e. invisible in the WM of the problem) feature, the phase-difference $\theta$, playing a non-negligible ``macroscopic'' role.

%

 This non-symmetry of the mass distribution depends on the phase separation of the two wave packets and on the value of $\hbar$.  Since mass is conserved -- analytically  as well as numerically -- the excess mass in one side is compensated by less mass on the other side of the corner. We measure this by the excess mass percentage ($\text{EMP}$) after the interaction,
\[
EMP=\|U(t_*) \chi_{x>0}\|^2_{L^2} -  \|U(t_*) \chi_{x<0}\|^2_{L^2} \in [-1,1].
\]
For time $t_* > 2.5$ so that the interaction is complete, and the two waves travel away from $x=0$ in opposite directions.
For $\theta=\frac14$ the wave packets have a $\frac{\pi}{2}$ phase difference and more mass, is located to the right of the corner, $\text{EMP}\approx 5\%$. In a completely analogous way for a phase separation of $\frac{3\pi}{2}$ ($\theta=\frac34$) the exactly same amount of excess mass is shifted to the left of the corner.  However for $\theta=0, \frac12, 1$ the mass is distributed equally around the corner, $\text{EMP}=0$. 

\begin{figure}[htbp]
\begin{center}
\includegraphics[scale=0.5]{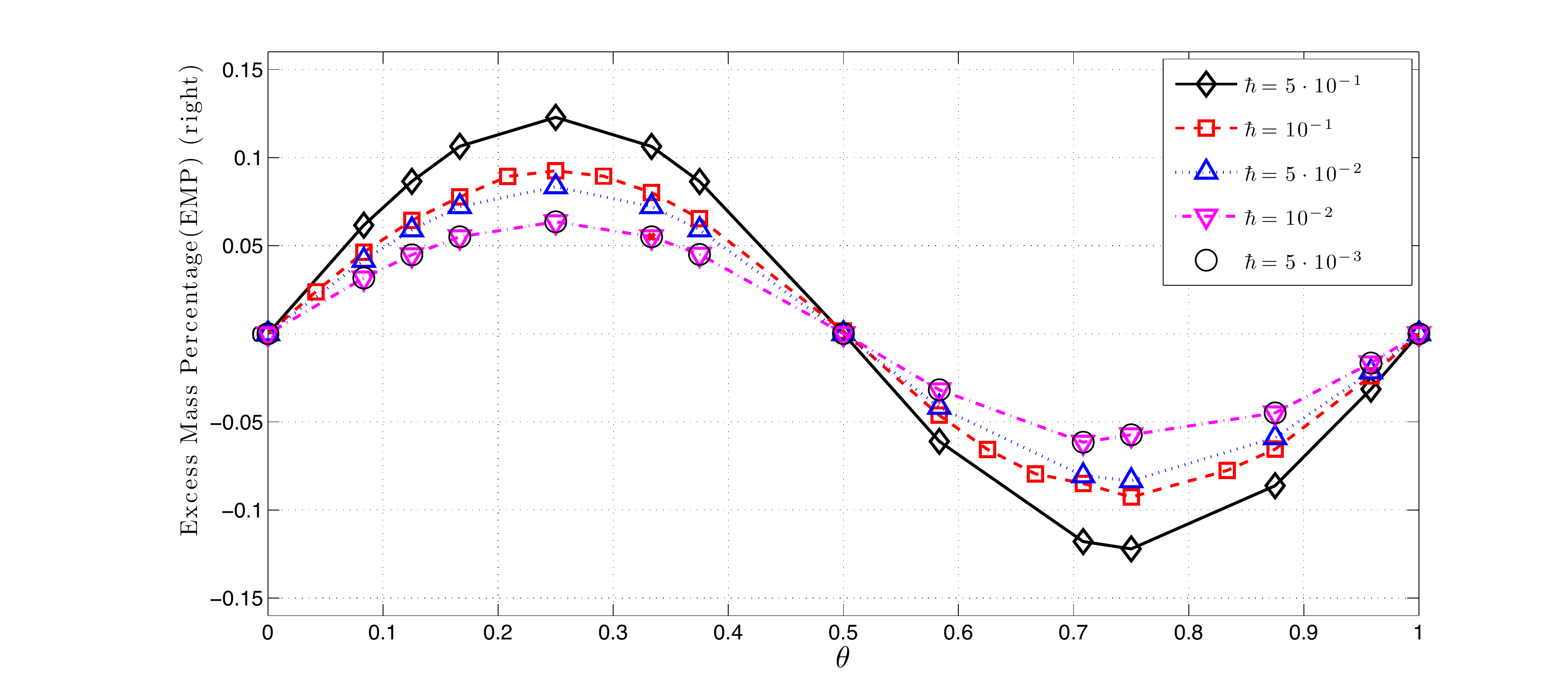}
\caption{EMP distribution for various values of  $\theta$ and $\hbar$.}
\label{massdistr}
\end{center}
\end{figure}

The dependance of the mass imbalance on  Planck's constant $\hbar$ is not easily visible from Figure \ref{Coll2}. To clarify the situation we run several numerical experiments for a variety of values of $\hbar$ and $\theta$ :
\begin{align*}
& \theta = (2, 3, 4, 6, 8, 9, 12, 14, 17, 18, 21, 23) / 24 ,  \\
& \hbar = 5\cdot 10^{-1}, 10^{-1}, 5\cdot 10^{-2}, 10^{-2}, 5\cdot 10^{-3} . 
\end{align*}
The results are summarized in Figure \ref{massdistr} where the variation of the EMP to the right of the corner $x=0$ is shown. The dependence on the value of $\hbar$ is evident. The location of maximum and minimum values of $\text{EMP}$ depend solely on the value of $\theta$  and occurs for $\theta=\frac14$ and $\theta=\frac34$ respectively but does not depend on $\hbar$. The value of this maximum and minimum depend on the value of $\hbar$, and seem to stabilize for $\hbar$ small enough. For $\hbar=\mathcal{O}(1)$,  $\text{EMP}\approx 12.5\%$, and it reaches an apparent limiting value of $\text{EMP}\approx 5.5\%$  for $\hbar=5\cdot 10^{-3}$.  We also notice  for $\theta=0, \frac12, 1$ and for any value of $\hbar$ the mass is distributed evenly around the corner, $\text{EMP}=0$. The behavior encoded in Figure \ref{massdistr} seems to persist even if change the envelopes $a_{0,1}$, $a_{0,2}$; i.e. there seems to be a quantum scattering operator that depends only on the phase difference of the interfering waves.

\begin{remark} \upshape While the discrepancy between $\rho^\hbar(t)$, $W^\hbar(t)$ after the interference effects is clear, one could think that this is only due to the smoothing of the initial data. In other words, is maybe $f^\hbar$, defined as in eq. \eqref{eq:ioghh}, close to $W^\hbar(t)$?
 First of all, working with the full detail of $W^\hbar_0$ is not a practical asymptotic method, as is clearly seen by plotting it for small $\hbar$ (see eq. \eqref{eqwtbdtwwvf} for the explicit form). Still, we did investigate its behavior; $f^\hbar(t)$ gives rise to comparable (sometimes larger) mass imbalances as $\rho^\hbar$. In other words,
\[
\lim\limits_{\hbar\to 0}\langle \rho^\hbar(t_*)- W^\hbar(t_*), \phi\rangle \neq 0, \qquad \qquad \lim\limits_{\hbar\to 0}\langle f^\hbar(t_*)- W^\hbar(t_*), \phi\rangle \neq 0,
\]
at least for the observables corresponding to ``mass on the left/right''. This reinforces the conclusion that the interference effect is genuinely quantum, and thus cannot be captured by the solution of a Liouville equation on its own. This is also consistent with the modeling of \cite{coherent}.
\end{remark}

\appendix
\section{Background on the Schr\"odinger equation and the Wigner transform} \label{sec:appA}
The Schr\"odinger equation (\ref{eq:schroeq}) is well-posed on $L^2(\mathbb{R}^d)$ for real potentials $V$ in Kato's class, i.e. if $V=V_1+V_2$ and
\beq
V_1 \in L^\infty, \,\,\,\,\,\,\, V_2 \in L^p, \,\, p> \max \{2,\frac{d}2 \}
\ee
or
\beq
V_1 \in L^\infty, \,\,\,\, V_2 \in L^2_{loc}, \,\, \,\, \mbox{ and } \exists C>O \mbox{ such that } V_2>-C(1+|x|^2).
\ee
Practically all physically interesting cases are covered by these conditions -- unlike the situation in classical mechanics. 
The Wigner transform (WT),
\beq
W^\hbar:L^2(\mathbb{R}^d)\times L^2(\mathbb{R}^d) \rightarrow L^2(\mathbb{R}^{2d}): f,g \mapsto W^\hbar[f,g]=\int e^{-2\pi i k y} f(x+\frac{\hbar y}2) \overline{g}(x-\frac{\hbar y}2)\,dy,
\ee
 seen as a bilinear mapping is essentially unitary in $L^2$, in the sense that
\beq \label{eq:wl2ret}
\|W^\hbar[f,g]\|_{L^2(\mathbb{R}^{2d})}=\hbar^{-\frac{d}2} \|f\|_{L^2(\mathbb{R}^{d})}\|g\|_{L^2(\mathbb{R}^{d})}.
\ee
This allows the construction of an $L^2$ propagator for the Wigner equation out of the Schr\"odinger propagator \cite{Mql}.
We would like to interpret the WT as a phase-space probability density in the sense of classical statistical mechanics; it has e.g. the correct marginals as position and momentum density
\beq
\int{W^\hbar[f](x,k)\,dk} = |f(x)|^2, \,\,\,\,\, \int{W^\hbar[f](x,k)\,dx} = |\widehat f(x)|^2.
\ee
However this picture cannot be taken too literally, since the WT has negative values in general \cite{Cohen,Hlaw}.
In fact, it has been realized that when smoothed with an appropriately large kernel, the WT becomes non-negative. Skipping over some details, this can be seen as an equivalent reformulation of the Heisenberg uncertainty principle: one can get a valid (i.e. a priori non-negative) probability that a particle occupies a region in phase-space only if that region is large enough. This leads to the definition of the Husimi transform,
\beq
H^\hbar[f](x,k)=\left({ \frac{2}{\hbar} }\right)^d e^{-\frac{2\pi}{\hbar}\left[{ |x|^2+|k|^2 }\right]} \ast W^\hbar [f] \geqslant 0 \,\,\, \forall f \in L^2.
\ee
The Husimi transform is used to prove the positivity of the WM, since, as can be readily checked, $W^\hbar[u^\hbar]$ and $H^\hbar[u^\hbar]$ are close in weak sense as $\hbar\shortrightarrow 0$ \cite{LP}.

The particular topology used for weak-$*$ convergence $W^\hbar[u^\hbar], H^\hbar[u^\hbar] \rightharpoonup W^0$ is built on the algebra of test functions $\mathcal{A}$, introduced in \cite{LP} and defined as
\beq \label{eqA}
\mathcal{A} = \{ \phi \in C(\mathbb{R}^{2d}) \,|\,\, \int{ \mathop{sup}\limits_{x} |\mathcal{F}_{k \shortrightarrow K}[\phi(x,k)]|dK } < \infty \}.
\ee
The main result for Wigner measures in smooth problems, precisely stated, is the following

\begin{theorem}[Wigner Measures for the linear Schr\"odinger equation \cite{LP,GMMP}] \label{thrm1} Let the real valued potential $V$ be in Kato's class, and assume there exists a $C>0$ such that $V(x)\geqslant -C(1+|x|^2)$. Assume moreover that the family of initial data $\{u_0^{\hbar_n}\}$, for a sequence $\lim\limits_{n\rightarrow \infty} \hbar_n=0$, has the following properties
\begin{itemize}
\item
\textbf{($\hbar$-oscillation)}
If $F_\phi(R)$ is defined as
\[
F_\phi(R)= \limsup\limits_{n\rightarrow \infty} \int\limits_{|k|\geqslant \frac{R}{\hbar_n}} { |\widehat{ \phi u^{\hbar_n}_0 }|^2\, dk},
\]
then, for all continuous, compactly supported $\phi$
\[
\lim\limits_{R\rightarrow \infty} F_\phi(R)=0.
\]
\item
\textbf{(compactness)} If $G(R)$ is defined by
\[
G(R)= \limsup\limits_{n\rightarrow \infty} \int\limits_{|x|\geqslant R} { | u^{\hbar_n}_0 |^2\, dx},
\]
then
\[
\lim\limits_{R\rightarrow \infty} G(R)=0.
\]
\end{itemize}

\vspace{0.5cm}
Then, for a semiclassical family of problems of the form (\ref{eq:schroeq}), and any timescale $T>0$, the following hold:
\begin{itemize}
\item
There exists a subsequence of the initial data, $u^{\hbar_{m_n}}_0$, so that their Wigner transform converges in $\mathcal{A}'$ weak-$*$ sense to a probability measure,
\[
\bac
\forall \phi \in \mathcal{A} \,\,\, \lim\limits_{n\rightarrow \infty} \langle{ W^{\hbar_{m_n}}_0  - W_0^0, \phi}\rangle=0, \,\,\,\,\,\,\,\,\,\,\,\,\,\,
W^0_0 \in \mathcal{M}^1_+(\mathbb{R}^{2d})
\ea
\]
\item
For $t\in[0,T]$,
define $W^0(t)$ as the propagation of the initial Wigner measure $W^0_0$ under the Liouville equation \eqref{eq:wmeq}. Then
\[
W^\hbar[u^{\hbar_n}(t)]=W^{\hbar_n}(t) \rightharpoonup W^0(t)
\]
in $\mathcal{A}'$ weak-$*$ sense. 
\end{itemize}
\end{theorem}

\medskip

Finally, one should note that if $f,g$ are ``nice enough'', then their WT $W^\hbar[f,g]$ will also be ``nice'':

\begin{definition}{$\Sigma^m$} \label{def:SIGMA} We will say that $f\in L^2(\mathbb{R}^d)$ belongs to  $\Sigma^m$ if
\[
\| f\|_{\Sigma^m}=\max\limits_{|a|, |b|\leqslant m} \| x^a \partial_x^b f \|_{L^2} < +\infty.
\]
\end{definition}

\noindent {\bf Remark:} It is clear that $f \in \Sigma^m \Rightarrow \widehat{f} \in \Sigma^m$.

\begin{theorem} \label{thrm:wigregureg}If $f,g \in \Sigma^m(\mathbb{R}^{d})$ (see Definition \ref{def:SIGMA} above), then
\beq \label{eq:nwjhg}
W^\hbar[f,g] \in \Sigma^m(\mathbb{R}^{2d}).
\ee
Moreover,
\beq \label{eq:nwjhg2}
\mbox{ if } m>2d, \mbox{ then } W^\hbar[f,g] \in L^1\cap L^\infty.
\ee
\end{theorem}

\noindent {\bf Proof:} Observe that if
\[
R:F(x,y) \mapsto F(x+Cy,x-Cy), 
\]
the Wigner transform can be seen just as a composition of
\[
W^\hbar[f,g] = \mathcal{F}_{y\to k} R_{\frac{\hbar}2} f(x)\overline{g}(y).
\]

So the strategy for the proof of \eqref{eq:nwjhg} is clear; show that $f(x)\overline{g}(y) \in \Sigma^m(\mathbb{R}^{2d})$, and then show that each of the operators $R_C$, $\mathcal{F}_{y\to k}$ are bounded on $\Sigma^m(\mathbb{R}^{2d})$.

So, if
$f,g \in \Sigma^m(\mathbb{R}^{d}),$
it readily follows that
\[
\bac
\| f(x)  \overline{g}(y) \|_{\Sigma^m(\mathbb{R}^{2d})} =
\max\limits_{|a_1|+|a_2|, |b_1|+|b_2| \leqslant m} \|  x^{a_1}y^{a_2} \partial_x^{b_1} \partial_y^{b_2} f(x)\overline{g}(y) \|_{L^2} \leqslant\\

\leqslant \max\limits_{|a_1|, |b_1| \leqslant m} \|  x^{a_1} \partial_x^{b_1}  f(x) \|_{L^2} \quad
\max\limits_{|a_2|, |b_2| \leqslant m} \|  y^{a_2}  \partial_y^{b_2} \overline{g}(y) \|_{L^2} \leqslant \| f\|_{\Sigma^m(\mathbb{R}^{d})} \| g\|_{\Sigma^m(\mathbb{R}^{d})}.
\ea
\]

Now assume $F(x,y) \in \Sigma^m(\mathbb{R}^{2d})$;
\[
\bac
\|R_C F\|_{\Sigma^m(\mathbb{R}^{2d})} \leqslant  \max\limits_{|a_1|+|a_2|, |b_1|+|b_2| \leqslant m} \|  x^{a_1}y^{a_2} \partial_x^{b_1} \partial_y^{b_2} F(x+Cy,x-Cy) \|_{L^2}=\\

= \max\limits_{|a_1|+|a_2|, |b_1|+|b_2| \leqslant m} \|  \left({ \frac{x+Cy \,\, + \,\, (x-Cy)}{2} }\right)^{a_1} \left({ \frac{x+Cy  \,\,- \,\, (x-Cy)}{C} }\right)^{a_2}
\partial_x^{b_1} \partial_y^{b_2} F(x+Cy,x-Cy) \|_{L^2} \leqslant \\

\leqslant C' \|F\|_{\Sigma^m(\mathbb{R}^{2d})}.
\ea
\]

Finally
\[
\bac
\|\mathcal{F}_{y\to k} F\|_{\Sigma^m(\mathbb{R}^{2d})} \leqslant  \max\limits_{|a_1|+|a_2|, |b_1|+|b_2| \leqslant m} \|  x^{a_1}k^{a_2} \partial_x^{b_1} \partial_k^{b_2}  \widehat{F}_2(x,k) \|_{L^2}=\\

\leqslant  \max\limits_{|a_1|+|a_2|, |b_1|+|b_2| \leqslant m}  (2\pi)^{b_2-a_2}  \| 
x^{a_1} \partial_x^{b_1} \partial_y^{a_2} \big(y^{b_2} F(x,y) \big) \|_{L^2}
 \leqslant C'' \| F\|_{\Sigma^m(\mathbb{R}^{2d})}.
\ea
\]

Now to prove \eqref{eq:nwjhg2}; by virtue of the Sobolev embedding Theorem \cite{Evans}
\[
\| F \|_{L^\infty(\mathbb{R}^{2d})} \leqslant \| F\|_{H^{d+1}(\mathbb{R}^{2d})} \leqslant \| F\|_{\Sigma^{d+1}(\mathbb{R}^{2d})};
\]
moreover
\[
\bac
\| F \|_{L^1(\mathbb{R}^{2d})} = \int\limits_{x,y} |F(x,y)|dxdy =  \int\limits_{x,y} |F(x,y)|\frac{(1+|x|^2+|y|^2)^r}{(1+|x|^2+|y|^2)^r} dxdy \leqslant\\

\leqslant \| (1+|x|^2+|y|^2)^r F(x,y) \|_{L^2(\mathbb{R}^{2d})}  \| \frac{1}{(1+|x|^2+|y|^2)^r} \|_{L^2(\mathbb{R}^{2d})} \leqslant 

\| F\|_{\Sigma^{4r}(\mathbb{R}^{2d})} \sqrt{ \int\limits_{\rho=0}^\infty \frac{\rho^{2d-1} d\rho}{(1+\rho^2)^{2r}} }
\ea
\]
which is finite for $r> \frac{d}2$.
\qed

\section{The Smoothed Wigner Transform}
\label{AppCoarse}

As was mentioned, sometimes flexibility in the calibration of the smoothing is required. Several approaches for the smoothing of the Wigner transform have been studied \cite{Cohen,Hlaw}, and there exist trade offs for the different choices and scalings of smoothing kernels. We use a Gaussian smoothing in what we call the Smoothed Wigner transform (SWT). This has the advantage that it leads to entire analytic functions of known order and type, thus making available a great toolbox of results for their asymptotic study \cite{AMP}.

The SWT was introduced in 
 \eqref{eq:swtdef}. Observe that
\begin{equation}
\begin{aligned}
\label{eq:swtcompo}
\widetilde{W}^\hbar[u](x,k)&=\left({ \frac{2}{\hbar {\sigma_x \sigma_k}} }\right)^d \int\limits_{x,k} e^{-\frac{2\pi}{\hbar}\left[{ \frac{|x-x'|^2}{\sigma_x^2}+\frac{|k-k'|^2}{\sigma_k^2} }\right]}  W^\hbar (x',k')\,dx'\,dk'=\\
  &  =\left({ \frac{\sqrt{2}}{ \sqrt{\hbar} {\sigma_x }}   }\right)^d \int\limits_y{
e^{-2\pi i k y - \frac{\hbar \pi}{2}\sigma_k^2 y^2}
  \int\limits_{x'} { e^{-\frac{2\pi}{\hbar} \frac{|x-x'|^2}{\sigma_x^2}  } u(x'+\frac{y\hbar}2) \overline{u}(x'-\frac{y\hbar}2)\, dx'\,dy} },
\end{aligned}
\end{equation}
therefore only $d$ convolutions are needed (i.e. in $x$), as the smoothing in $k$ can be performed as part of the FFT.

\begin{figure} 
\includegraphics[height=8cm,width=5cm]{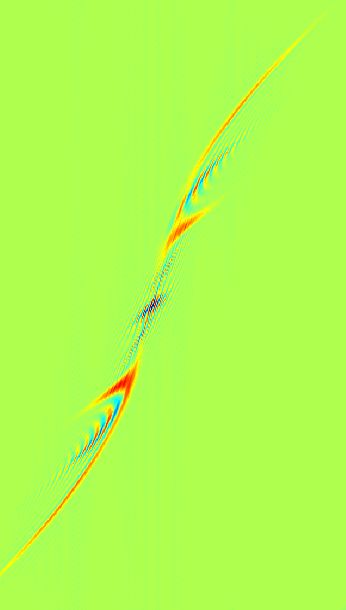}
\includegraphics[height=8cm,width=5cm]{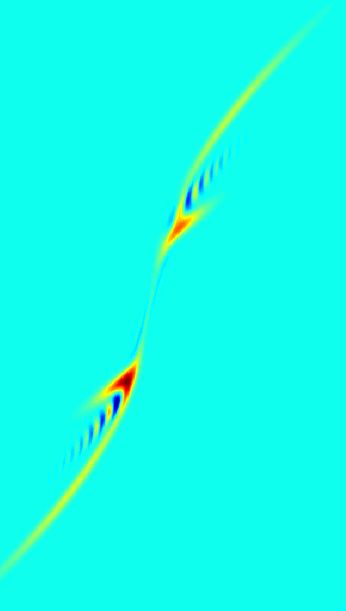}
\includegraphics[height=8cm,width=5cm]{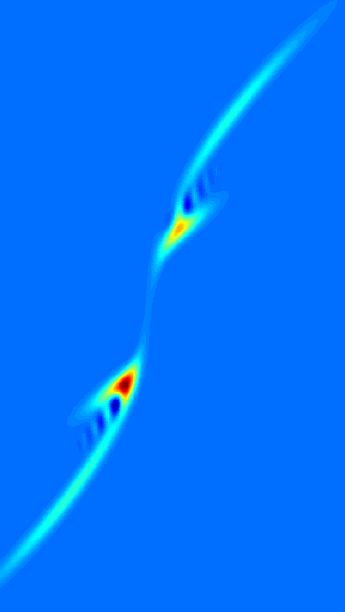}
\caption{Smoothing of the Wigner transform. Left; the WT -- dominant features correspond to oscillations that vanish in the limit. Middle; fine smoothing -- the most spurious oscillations are gone, there is good resolution and some very negative values. Right; coarser smoothing, more appropriate for computational use -- there are still non-negligible negative values, but the dominant features of the density are clearly positive while definition has not been overly smeared.}
\label{FigSWT}
\end{figure}

To implement this transform numerically, we will use the FFT. First of all recall that
\begin{lemma} \label{psffft}
For any function $f \in \mathcal{S}(\mathbb{R})$,
\[
 \sum\limits_{j \in \mathbb{Z}} f(jh) e^{-2\pi i k nh} = \frac{1}{h} \sum\limits_{j \in \mathbb{Z}} \widehat{f}(k + \frac{j}{h}).
\]
\end{lemma}
This is a direct corollary of the Poisson summation formula, and the starting point of any use of the FFT to approximate the Fourier transform of a continuous function. (The requirement $f \in \mathcal{S}(\mathbb{R})$ can be relaxed; the details along this direction are outside the scope of this work.)

 Similarly, we can create an appropriate version of the Poisson summation formula for the evaluation of the $dy$ integral in \eqref{eq:swtcompo} as an FFT:
\begin{lemma} If $f \in \mathcal{S}(\mathbb{R})$, denote by
$$
S_{a,b}(X,y)=e^{-\frac{\pi }{2} \hbar\sigma_k^2y^2  }  \int\limits_{x'} {e^{-\frac{2\pi}{\hbar \sigma_x^2}(X-x')^2} f(x'+\hbar b y) \overline{f}(x'-\hbar b y)\,  dx'}.
$$
Then
$$
\sum\limits_{j \in \mathbb{Z}} S_{a,b}(X,j) e^{-2\pi i 2Ka j} = \frac{\sigma_k^2}{b \sqrt{2}} \sum\limits_{j \in \mathbb{Z}} \widetilde{W}[u](X,\frac{aK+j}{b}).
$$
\end{lemma}

Observe that the integral $\int\limits_{x'} {e^{-\frac{2\pi}{\hbar \sigma_x^2}(X-x')^2} f(x'+\hbar b y) \overline{f}(x'-\hbar b y)  \,dx'}$ only needs to be computed in a small interval in $x'$ for each $X$ because of the Gaussian localization.
%
%

\bigskip

It is handy to note that smoothing in the $k$ direction introduces small changes in the $\vertiii{\cdot}_{-M}$ norm (see also Observation \ref{obs:lstobsggt}):

\begin{lemma} \label{lm:Msmstr} Let $\|u_0^\hbar\|_{L^2}=1$, $W^\hbar_0 = W^\hbar[u^\hbar_0]$, and $\rho^\hbar = \widetilde{W}^{0,\sigma_k;\hbar}_0$. Then
\[
\vertiii{ W^\hbar_0-\rho^\hbar }_{-M} \leqslant  \hbar \frac{\pi}{2} \sigma_k^2 M^2.
\]
\end{lemma}

\noindent {\bf Proof: }
\[
\bac
|\langle \rho^\hbar - W^\hbar_0 , \phi \rangle| = |\langle \mathcal{F}_2 W_0^\hbar, (1-e^{-\frac{\hbar\pi}{2} \sigma_k^2 K^2}) \mathcal{F}_2 \phi \rangle| \leqslant \frac{\hbar\pi}{2} \sigma_k^2   \||K^2| \widehat{\phi}_2 \|_{L^1_KL^\infty_x} \leqslant \frac{\hbar\pi}{2} \sigma_k^2 M^2 \vertiii{\phi}_M.
\ea
\]

We used Observation \ref{obs:apriregW}, and the fact that
\[
\|K^2 \widehat{\phi}_2(x,K)\|_{L^1_KL^\infty_x} \leqslant \|K^2 \widehat{\phi}(X,K)\|_{L^1_{X,K}} \leqslant M^2 \vertiii{\phi}_M.
\]
\qed

\noindent {\bf Remark: } Observe that in practice 
\[
\lim\limits_{\hbar \to 0} \langle W^\hbar_0 - \widetilde{W}^\hbar_0 , \phi \rangle =0
\]
(which is significantly weaker than \eqref{eq:vertiiconstr} or lemma \ref{lm:Msmstr} above) is often sufficient in practice. In particular, the use of some smoothing in the $x$ direction as well, does not seem to hurt the quality of approximation in our numerical examples. In any case, we have an explicit, uniform upper bound for the effect of smoothing in the $k$ variables only.

\section{Slicing in two of a WKB wavefunction} \label{AppCwkb}

\begin{figure}[htb!]
\begin{tabular}{l r}
\includegraphics[width=68mm,height=60mm]{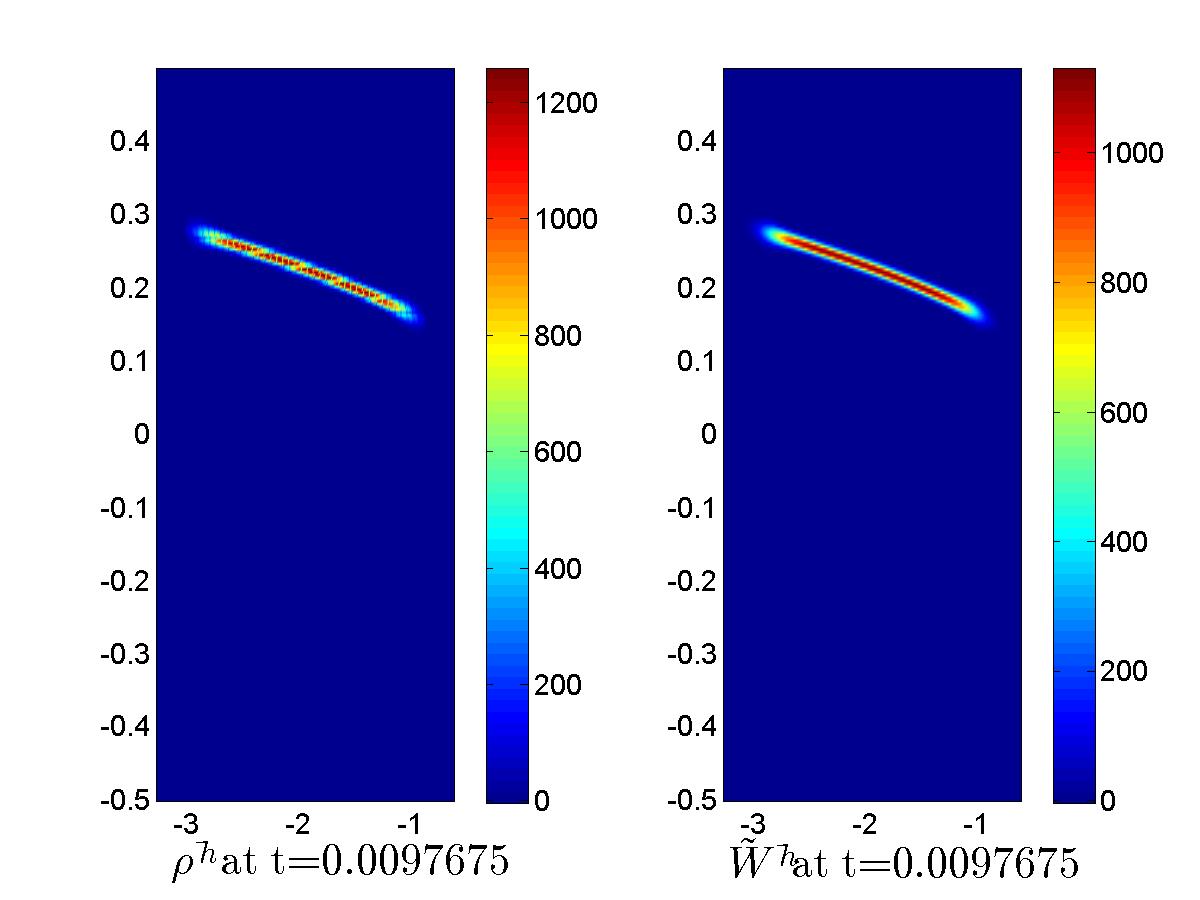} &
\includegraphics[width=68mm,height=60mm]{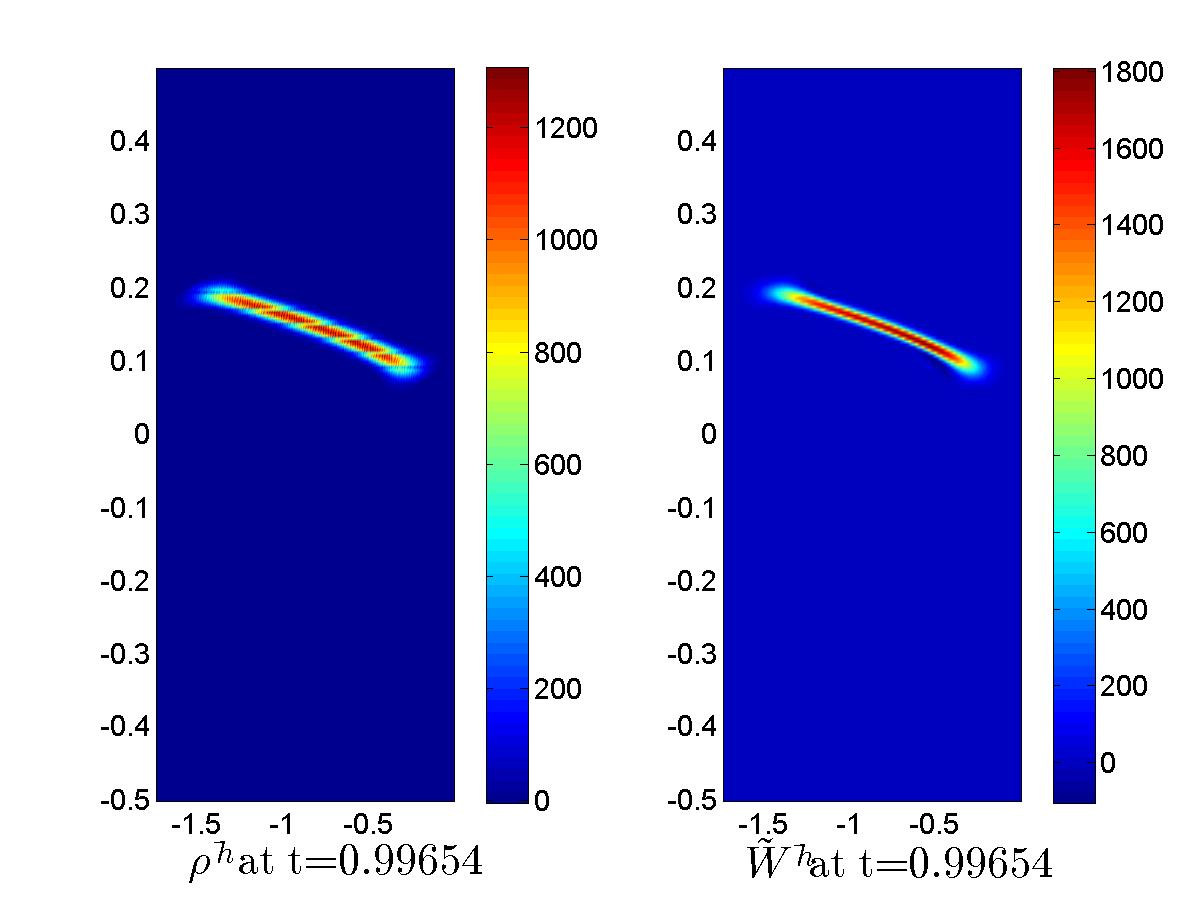} \\
\includegraphics[width=68mm,height=60mm]{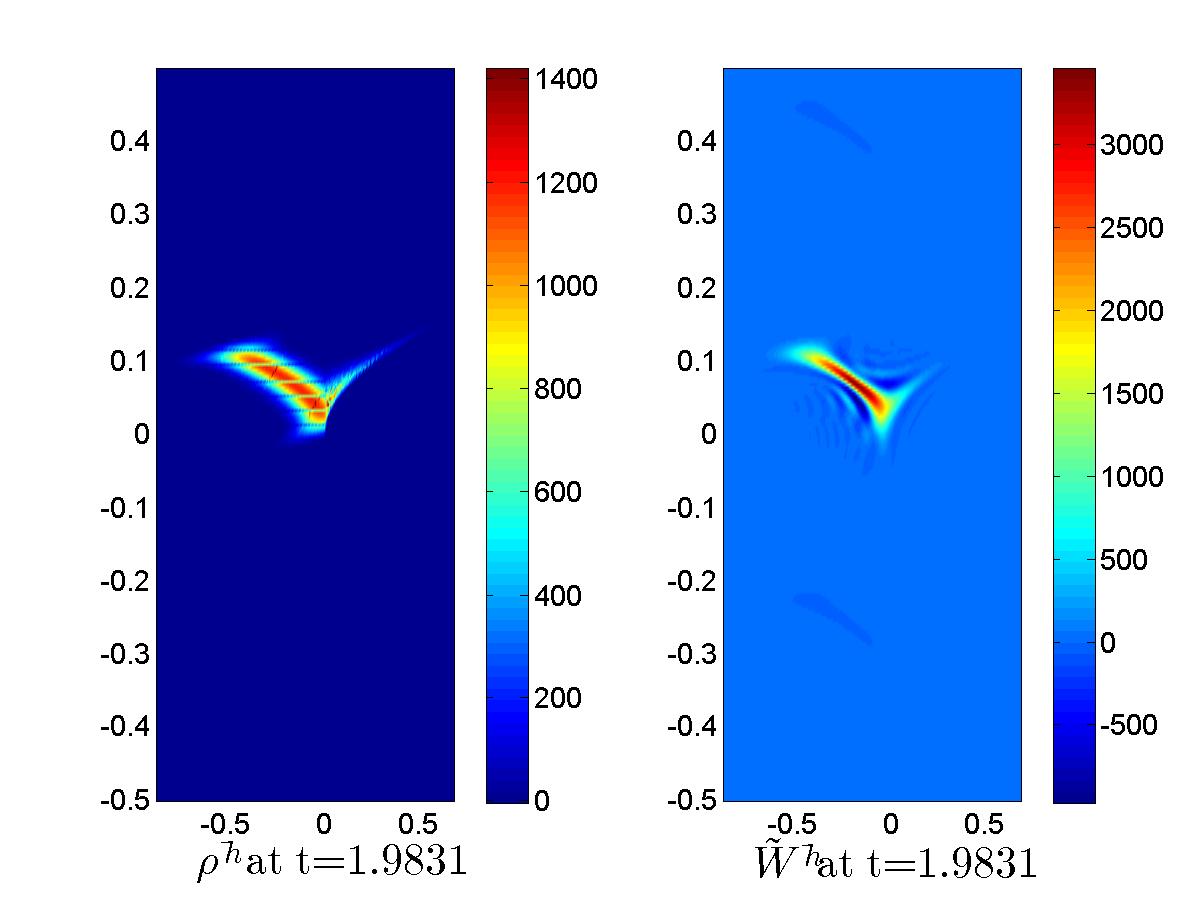} &
\includegraphics[width=68mm,height=60mm]{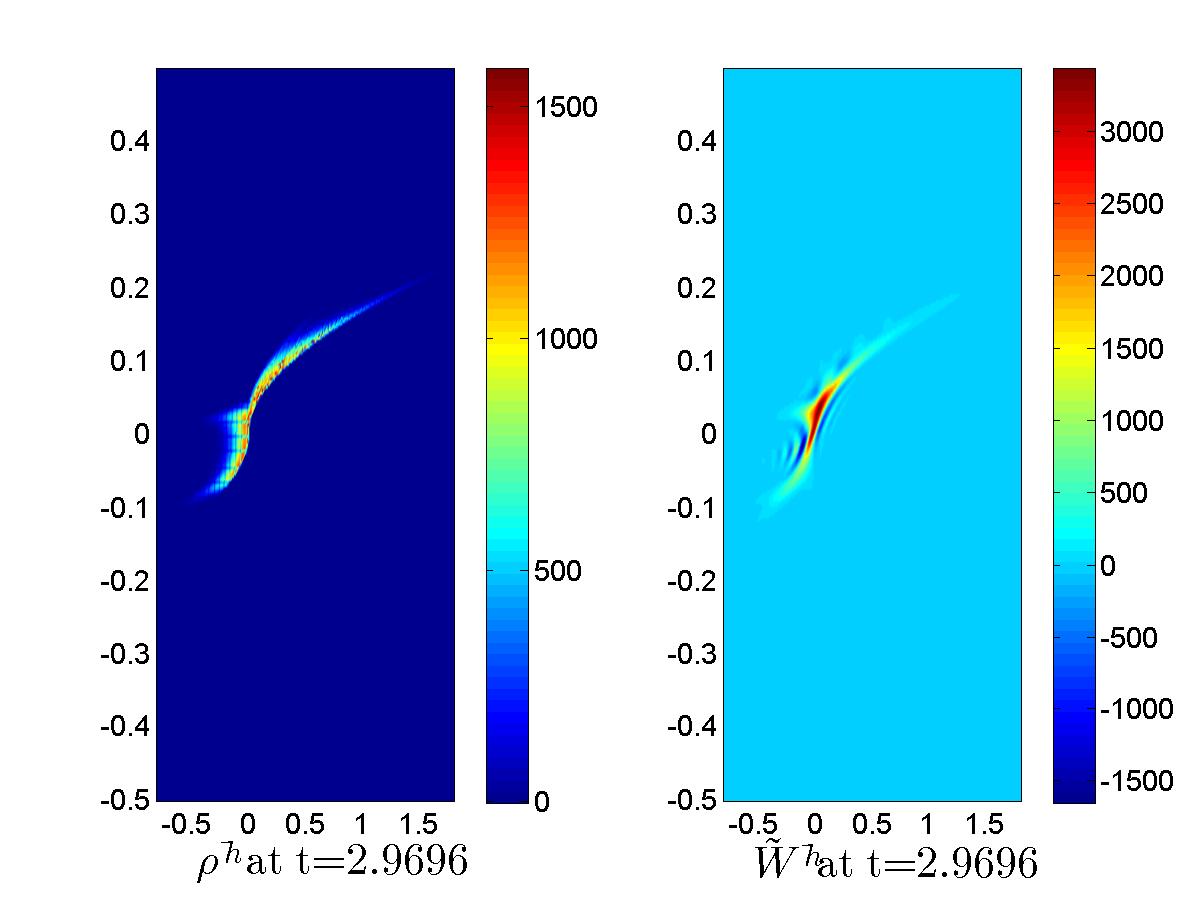} \\
\includegraphics[width=68mm,height=60mm]{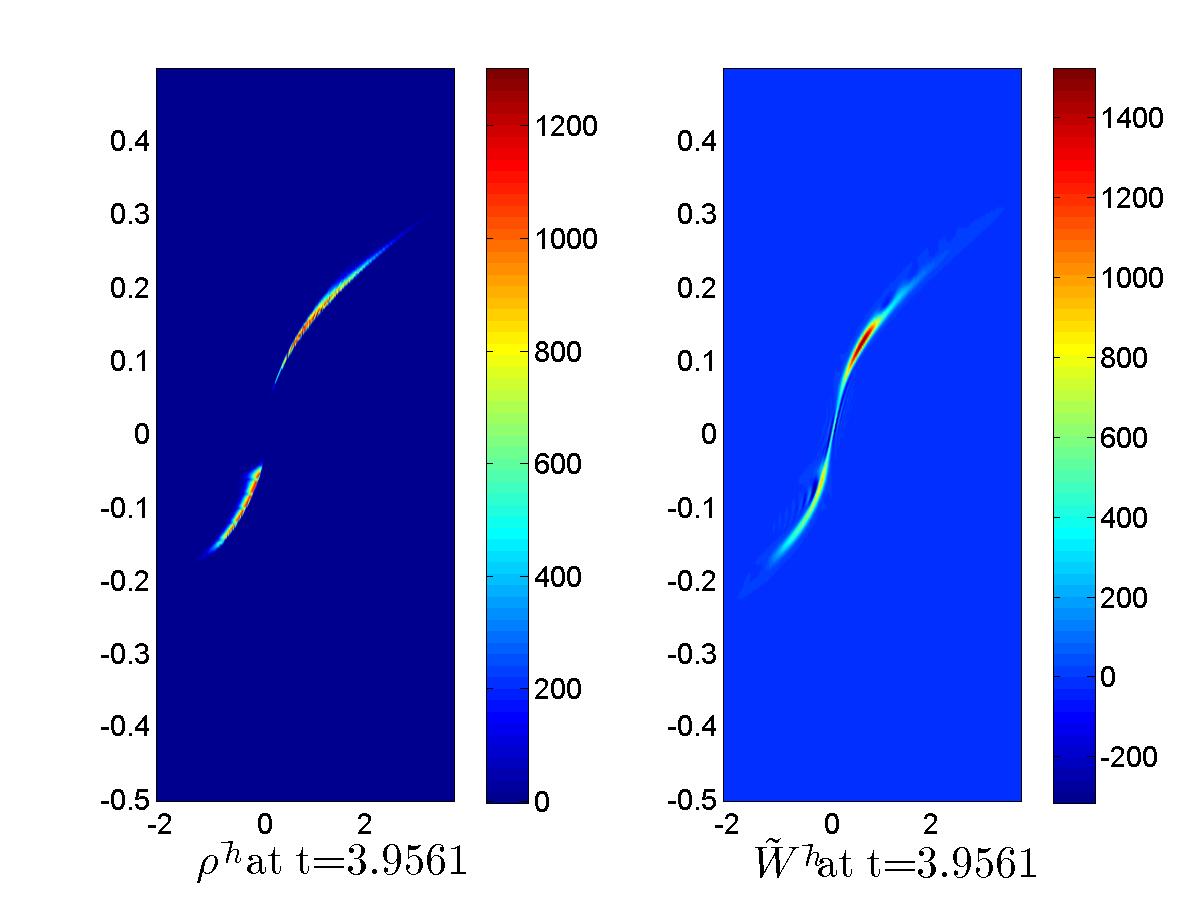} &
\includegraphics[width=68mm,height=60mm]{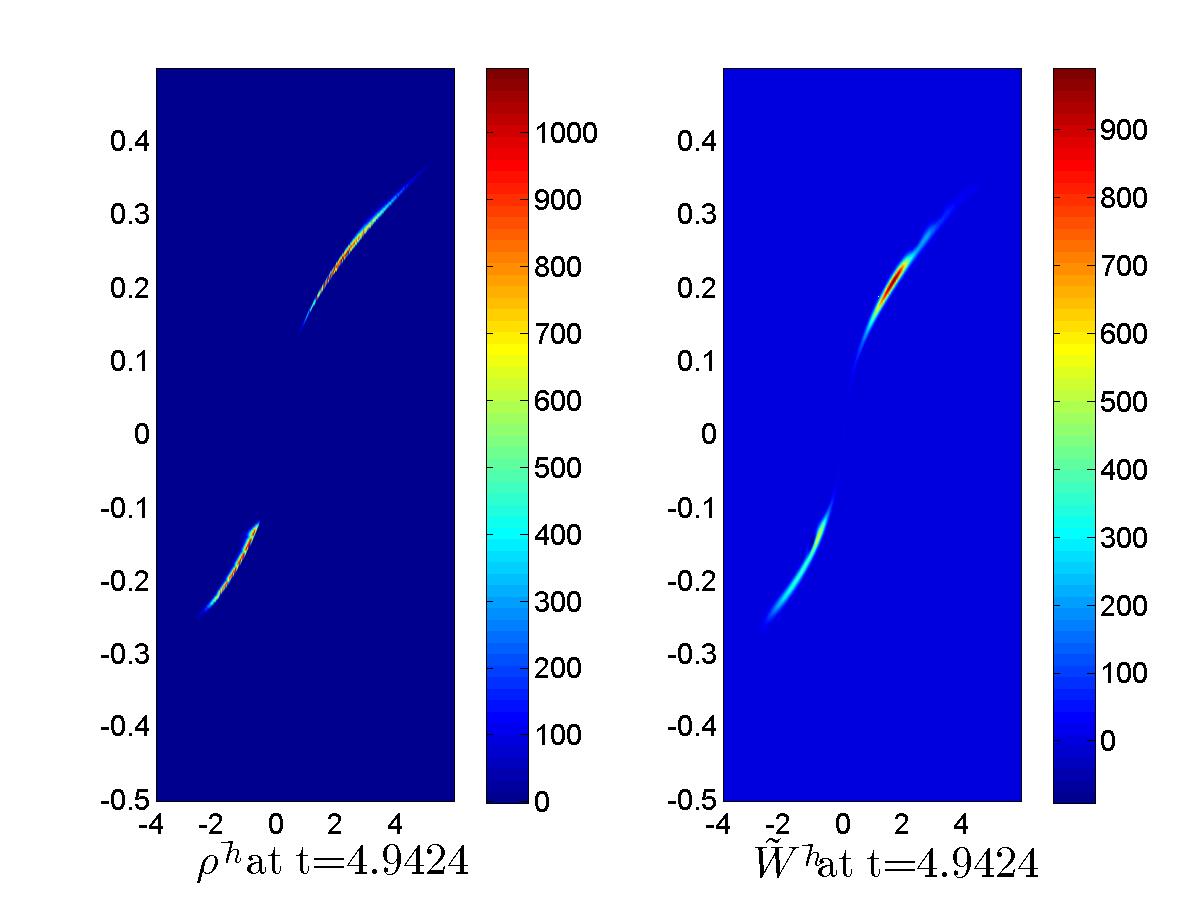} 
\end{tabular}
\caption{ The quantum phase-space density $\widetilde{W}^\hbar(t)$ versus its proposed semiclassical approximation $\rho^\hbar(t)$, for the data of eq. \eqref{eq:A3___1}, $\hbar=10^{-2}$ and various times before, during, and after interaction with the singularity.
}
\label{figA3___1}
\end{figure}

A more singular non-interference problem example is given by the initial data
\beq
\bac \label{eq:A3___1}
 u_0(x) = a_{0}(x) \mathrm{e}^{i  \frac{S_{0}(x)}{\hbar}}, \quad
 a_{0}(x)=  (1+tanh(7(x+3)) )\cdot(1+tanh(7(-x+1)) ), \quad
 S_{0}(x) =  \frac{-2}3  |x|^{\frac{3}{2}}
\ea
\ee
and
\begin{equation}
V(x) = 1 + (1+\tanh(4(x+4)))(1+\tanh(-4(x-4)))\frac{(-|x|+4)}{8}.
\end{equation}

The initial WM of this problem is 
 a line supported measure.  The concentration limit of $\rho^\hbar$ predicts that this measure would be ``sliced'' into two lines.  We see clear qualitative agreement between $\rho^\hbar$ and $\widetilde{W}^\hbar$, see Figure \ref{figA3___1}. The quantum observables \eqref{eqobsset} (including mass scattered to the left / right) are within around $4\%$ of their semiclassical prediction. Overall it seems that quantitative convergence as $\hbar \rightarrow 0$ is taking place, albeit somewhat more slowly for this type of initial data than for the data of eq. \eqref{u0nonint}.
 

\end{document}